\documentclass[a4paper,12pt]{amsart}

\usepackage[top=3cm,bottom=2.85cm,outer=3cm,inner=2cm,marginpar=2.45cm]{geometry}
\usepackage[destlabel,final,colorlinks=true,urlcolor=RoyalBlue]{hyperref}
\usepackage[abbrev]{amsrefs}

\usepackage[french,ngerman,english]{babel}
\usepackage[utf8]{inputenc}
\usepackage[T1]{fontenc}

\usepackage[all,pdf]{xy}

\renewcommand{\objectstyle}{\displaystyle}

\usepackage[inline]{enumitem}
\usepackage{mathtools}  
\usepackage[usenames,dvipsnames]{xcolor}
\usepackage{calc} %% for doing length computations

\babeltags{de = ngerman}
\babeltags{fr = french}

\usepackage{upref}  %% for uprighting texts generated by \ref
\usepackage{embrac} %% for uprighting parentheses in \emph and
\usepackage[
]{marionotations}
\usepackage{bbm}     %% to use \mathbbm (like \mathbb but works also
\usepackage[Smaller]{cancel} %%%%% for crossing out argument in math mode via
\usepackage{breakurl}  %% used so that line breaks for contents in
\usepackage{subfiles}

\newcommand{\defaultDimension}{n}

\newcommand{\defaultAmbientSpace}{X}

\newcommand{\defaultlcIndex}{\sigma}
\newcommand{\setDefaultlcIndex}[1]{\renewcommand{\defaultlcIndex}{#1}}

\newcommand{\defaultcohDegree}{q}

\newcommand{\defaultlclocus}{S}
\newcommand{\setDefaultlclocus}[1]{\renewcommand{\defaultlclocus}{#1}}

\newcommand{\defaultvphi}{\vphi_L}
\newcommand{\setDefaultvphi}[1]{\renewcommand{\defaultvphi}{#1}}

\newcommand{\defaultpsi}{\psi}
\newcommand{\setDefaultpsi}[1]{\renewcommand{\defaultpsi}{#1}}

\newcommand{\defaultMetric}{\omega}
\newcommand{\setDefaultMetric}[1]{\renewcommand{\defaultMetric}{#1}}

\newcommand{\defaultMap}{\pi}

\newcommand{\alert}[2][RoyalBlue]{{\color{#1}#2}}

\NewDocumentCommand{\logKX}{
  D<>{\defaultAmbientSpace}  %% #1 base space
  t{.}                       %% #2 include lc locus if .-ed
  D(){\defaultlclocus}       %% #3 divisorial lc locus
  O{F}                       %% #4 replace F by the argument when provided
  t{M}                       %% #5 include M in the tensor product if present
}{K_{#1} \IfBooleanT{#2}{\otimes #3} \otimes #4 \IfBooleanT{#5}{\otimes M}}

\NewDocumentCommand{\vphilist}{
  t{.}                 %% #1 include "\phi_D" if .-ed
  G{\defaultvphi}      %% #2 the main potential "\vphi_F"
  t{M}                 %% #3 add "+\vphi_M" if present
  o                    %% #4 add any other potentials
  e{,}                 %% #5 replace the K\"ahler metric "\omega" if present
}{\IfBooleanT{#1}{\phi_{\defaultlclocus} +} #2
  \IfBooleanT{#3}{+\vphi_M} \IfNoValueF{#4}{+#4}, \IfNoValueTF{#5}{\defaultMetric}{#5}}

\NewDocumentCommand{\Ltwo}{ %% the space of L2 sections
  t{.}                      %% #1 with "loc" when .-ed
  D//{\bullet,\bullet}      %% #2 the order of forms
  D<>{\defaultAmbientSpace} %% #3 base space
  s                         %% #4 base space is hidden if * is present
  m                         %% #5 coefficient
}{L^{#2}_{(2)\IfBooleanT{#1}{\:\text{loc}}}\paren{\IfBooleanF{#4}{#3;} #5}}

\newcommand{\HarmSym}{\mathcal{H}}
\NewDocumentCommand{\Harm}{ %% the space of harmonic forms
  t{'}                      %% #1 with preassigned hol degree if present
  D//{\defaultcohDegree}    %% #2 anti-holomorphic degree
  D<>{\defaultAmbientSpace} %% #3 the base space
  g                         %% #4 the coefficient; will be hidden
  t{,}                      %% #5 separator
  G{\defaultvphi}           %% #6 potential on line bundle
  s                         %% #7 hide the ambient space metric if *-ed
  E{_}{\defaultMetric}      %% #8 metric on the base space or other subscripts
}{\HarmSym^{\IfBooleanT{#1}{\defaultDimension,}#2}\IfNoValueF{#4}{\paren{#3;#4}}_{
    \IfBooleanF{#1}{\IfNoValueT{#4}{#3,}} #6 \IfBooleanF{#7}{,#8}}}

\NewDocumentCommand{\lcIndex}{ %% for displaying the lc index,
  m  %% #1 the basic lc index (\sigma)
  m  %% #2 amount added to the index
  m  %% #3 amount substracted from the index
}{#1\IfNoValueF{#2}{+#2}\IfNoValueF{#3}{-#3}}

\NewDocumentCommand{\lcData}{ %% for displaying lc data in the format
  G{\defaultvphi}  %% #1 potential or q-psh function
  O{\defaultpsi}   %% #2 lc locus psi function
  e{.}             %% #3 jumping number
}{\paren{#1; \IfNoValueF{#3}{#3 \cdot} #2}}

\NewDocumentCommand{\lcdata}{ %% for displaying lc data in the
  s                %% #1 no parentheses if starred 
  d<>              %% #2 base space
  G{\defaultvphi}  %% #3 potential or q-psh function
  O{\defaultpsi}   %% #4 lc locus psi function
  e{.,}            %% #5 jumping number
}{{\newcommand{\datalist}{\IfNoValueF{#2}{#2,}#3,#4\IfNoValueF{#5}{,#5}\IfNoValueF{#6}{,#6}}
\IfBooleanTF{#1}{\datalist}{\paren{\datalist}}}}

\NewDocumentCommand{\drR}{ %% higher direct image (d*er*ived functor R*ight)
  m              %% #1 degree of the direct image
  O{\defaultMap} %% #2 the map
  d()            %% #3 coefficient
}{R^{#1}{#2}_*\IfNoValueF{#3}{\paren{#3}}}

\newcommand{\spHsym}{\mathbb{H}}
\NewDocumentCommand{\spH}{ %% cohomology group with coefficients 
  O{\spHsym}              %% #1 the symbol 
  D//{\defaultcohDegree}  %% #2 degree of anti-hol form
  t{M}                    %% #3 with 'M' to display M in the coefficient
  m                       %% #4 the coefficient 
}{#1^{#2}\paren{\IfBooleanT{#3}{M\otimes}#4}}

\def\spHarm{\spH[\HarmSym]}
\def\spR{\spH[\mathbb{R}]}

\newcommand{\deltaH}{\delta_{\HarmSym}}

\DeclareMathOperator{\lc}{lc} %% lc centre
\NewDocumentCommand{\lcc}{ %% union of lc centres
  D||{\defaultlcIndex}       %% #1 lc index \sigma
  e{+-}                      %% #2,#3
  D<>{\IfBooleanF{#7}{\defaultAmbientSpace}}  %% #4 base space
  t{'}                       %% #5 '-ed to show lc locus instead of
  D(){\defaultlclocus}       %% #6 lc locus
  s                          %% #7 *-ed to hide the parentheses for lc data
}{\lc_{#4}^{\lcIndex{#1}{#2}{#3}}\IfBooleanF{#7}{\IfBooleanTF{#5}{\paren{#6}}{\lcData}}}

\NewDocumentCommand{\lcS}{  %% a local lc centre
  s                       %% #1 symbol with \mathring when starred
  D(){\defaultlclocus}    %% #2 symbol for the subvariety
  D||{\defaultlcIndex}    %% #3 codimension
  e{+-}                   %% #4,#5
  d<>                     %% #6 open set where the lc centre lives
  O{p}                    %% #7 index among the \sigma-lc centres
}{\IfBooleanTF{#1}{\mathring{\mathtt{#2}}}{\mathtt{#2}}^{\lcIndex{#3}{#4}{#5}}_{\IfNoValueF{#6}{#6,}#7}}

\NewDocumentCommand{\PRes}{ %% Poincare Residue map
  O{}      %% subvariety
  d()      %% forms from the domain
}{\mathcal R_{#1}\IfNoValueF{#2}{\paren{#2}}}

\NewDocumentCommand{\HRes}{ %% Harmonic residue
  d()   %% #1 harmonic form
}{\mathfrak{R}\IfNoValueF{#1}{\paren{#1}}}

\newcommand{\defidlof}[1]{\mathcal{I}_{#1}}  %% defining ideal of (a set)
\NewDocumentCommand{\mtidlof}{   %% multiplier ideal of (a potential)
  O{}      %% #1 base space (for compatibility)
  D<>{#1}  %% #2 base space
  m        %% #3 potential / psh function
}{\multidl_{#2}\paren{#3}} 

\NewDocumentCommand{\residlof}{  %% multiplier ideal sheaf on the
  D||{\defaultlcIndex}   %% #1 codim of lc centres or supporting lc
  e{+-}                  %% #2,#3
  d<>                    %% #4 base space
  s                      %% #5 display the symbol without arguments when starred
  e{^_}                  %% #6,#7 extra super- and sub-scripts before \lcData
  t{,}                   %% #8 separator
}{\sheaf R_{\IfNoValueTF{#4}{}{#4,} \lcIndex{#1}{#2}{#3}\IfNoValueF{#7}{;#7}}\IfNoValueF{#6}{^{#6}}\IfBooleanF{#5}{\lcData}}

\NewDocumentCommand{\Adjidlof}{
  D||{\defaultlcIndex}       %% #1 codim of lc centres under concern
  D<>{\defaultAmbientSpace}  %% #2 base space
  D(){\defaultlclocus}       %% #3 lc locus
  m                          %% #4 potential or ideal
}{\operatorname{\mathit{Adj}}^{#1}_{\paren{#2,#3}}\paren{#4}}

\NewDocumentCommand{\aidlof}{
  D||{\defaultlcIndex}   %% #1 codim of lc centres under concern
  e{+-}                  %% #2,#3
  d<>                    %% #4 base space
  s                      %% #5 display the symbol without arguments when starred
  e{^_}                  %% #6,#7 extra super- and sub-scripts before \lcData
  t{,}                   %% #8 separator
}{\sheaf{J}_{\!\IfNoValueTF{#4}{}{#4,} \lcIndex{#1}{#2}{#3}\IfNoValueF{#7}{;#7}}\IfNoValueF{#6}{^{#6}}\IfBooleanF{#5}{\lcData}}

\NewDocumentCommand{\faidlof}{
  D||{\defaultlcIndex}   %% #1 codim of lc centres in numerator
  e{+-}                  %% #2,#3
  t{/}                   %% #4 a separator for arguments
  D||{\defaultlcIndex}   %% #5 codim of lc centres in denominator
  e{+-}                  %% #6,#7
}{\fracAidlof{\lcIndex{#1}{#2}{#3}}{\lcIndex{#5}{#6}{#7}}}

\NewDocumentCommand{\fracAidlof}{
  m                  %% #1 lcIndex in numerator
  m                  %% #2 lcIndex in denominator
  d<>                %% #3 base space
  s                  %% #4 display the symbol without arguments when starred
  e{^_}              %% #5,#6 extra super- and sub-scripts before \lcData
  G{\defaultvphi}    %% #7 potential or ideal
  O{\defaultpsi}     %% #8 defining function of the lc locus
  e{.}               %% #9 jumping number
}{\frac{
    \aidlof|#1|<#3>*^{#5}_{#6}\IfBooleanF{#4}{\lcData{#7}[#8].{#9}}
  }{
    \aidlof|#2|<#3>*^{#5}_{#6}\IfBooleanF{#4}{\lcData{#7}[#8].{#9}}
  }}

\NewDocumentCommand{\lcV}{ %% measure on lc centres
  D||{\defaultlcIndex}    %% #1 codim of supporting lc centres
  D//{\defaultvphi}       %% #2 potential for bundle valued section
  d()                     %% #3 metric on the ambient space
  e{^}                    %% #4 jumping number
  O{\defaultpsi}          %% #5 defining function (potential) of subvariety 
}{\:d\operatorname{lcv}^{#1\IfNoValueF{#4}{,\paren{#4}}}_{\IfNoValueF{#3}{#3,}#2}\left[#5\right]}

\NewDocumentCommand{\Ohvol}{ %% Ohsawa measure %%
  D//{\defaultvphi} %% #1 potential for bundle valued section
  d()               %% #2 metric on the ambient space
  O{\defaultpsi}    %% #3 defining function of subvariety
}{\dvol_{\IfNoValueF{#2}{#2,}#1}\left[#3\right]}

\newcommand{\dvol}{\:d\vol}

\NewDocumentCommand{\lcDataNormSubscript}{
  d<>                   %% #1 base space
  s                     %% #2 no potential and psi function when starred
  G{\defaultvphi}       %% #3 potential or q-psh function
  O{\defaultpsi}        %% #4 lc locus psi function
  e{.}                  %% #5 jumping number
  D||{\defaultlcIndex}  %% #6 lc Index
  e{+-}                 %% #7,#8
}{\IfNoValueF{#1}{#1,}
  \IfBooleanF{#2}{#3, \IfNoValueF{#5}{#5 \cdot} #4,}
  \lcIndex{#6}{#7}{#8}}

\newcommand{\RTFsym}{\mathfrak{F}} 
\NewDocumentCommand{\RTF}{ %% residue transform function
  s          %% #1 adding \smash[t] when starred
  G{\RTFsym} %% #2 symbol body
  o          %% #3 general superscript
  >{\SplitArgument{1}{,}} d<> %% #4 superscript in inner product
  d||        %% #5 superscript in \abs{}^2
  D(){\eps}  %% #6 for adding variable (\eps)
  t{,}       %% #7 separator
}{%
  \begingroup%
    \newif\ifsmasht%
    \IfBooleanTF{#1}{\smashttrue}{\smashtfalse}%
    \newif\ifboolup%
    \booluptrue%
    \IfNoValueT{#3}{\IfNoValueT{#4}{\IfNoValueT{#5}{\boolupfalse}}}%
    \newcommand{\supsrptstr}{\IfNoValueF{#3}{#3}\IfNoValueF{#4}{\inner#4}\IfNoValueF{#5}{\abs{#5}^2}}
    \newcommand{\RTFvar}{#6}
    #2\RTFprocess
}

\NewDocumentCommand{\RTFprocess}{
  o                     %% #1 overwrite subscript if given
  d<>                   %% #2 base space
  t{,}                  %% #3 with potential and psi function when ,-ed
  G{\defaultvphi}       %% #4 potential or q-psh function
  O{\defaultpsi}        %% #5 lc locus psi function
  e{.}                  %% #6 jumping number
  D||{\defaultlcIndex}  %% #7 lc Index
  e{+-}                 %% #8,#9
}{\newcommand{\subsrptstr}{%
    \IfNoValueTF{#1}{
    \IfNoValueF{#2}{#2,}
    \IfBooleanT{#3}{#4,#5,\IfNoValueF{#6}{#6,}}
    \lcIndex{#7}{#8}{#9}}{#1}}%
  \newcommand{\srptstr}{\cramped{{}^{\supsrptstr}%
      \ifboolup _
      \fi{\ifboolup\displaystyle\fi\paren{\RTFvar}%
          \ifboolup {\scriptstyle \subsrptstr} \else _{\subsrptstr} \fi%
        }}}%
  \ifboolup%
    \ifsmasht%
      \smash[t]{
        \raisebox{\depthof{$\srptstr$} * \real{0.3}}{$\srptstr$}%
      }%
    \else%
      \raisebox{\depthof{$\srptstr$} * \real{0.3}}{$\srptstr$}%
    \fi%
  \else%
    \srptstr%
  \fi%
  \endgroup%
}
\NewDocumentCommand{\mtlog}{O{e} d() D||{\defaultpsi}}{\log\!#1^{\paren{#2}}\abs{#3}}
\NewDocumentCommand{\slog}{O{e} D||{\defaultpsi}}{\log\abs{#1 #2}}
\NewDocumentCommand{\dlog}{O{e} D||{\defaultpsi}}{\mtlog[#1](2)|#2|}

\NewDocumentCommand{\logpole}{ %% log-pole in the residue transform
  D||{\defaultpsi}       %% #1 log singularity defining function
  t{,}                   %% #2 separator
  D||{\defaultlcIndex}   %% #3 codim of lc centres in question
  e{+-}                  %% #4,#5
  E{.^}{{e}{1+\eps}}     %% #6 multiplicative constant in logarithm 
  s                      %% #8 no parentheses and exponent on log|\psi| when starred
}{\abs{#1}^{\lcIndex{#3}{#4}{#5}} \IfBooleanTF{#8}{\slog[#6]|#1|}{\paren{\slog[#6]|#1|}^{#7}}}

\DeclareFontFamily{OMX}{MnSymbolE}{}
\DeclareSymbolFont{MnLargeSymbols}{OMX}{MnSymbolE}{m}{n}
\SetSymbolFont{MnLargeSymbols}{bold}{OMX}{MnSymbolE}{b}{n}
\DeclareFontShape{OMX}{MnSymbolE}{m}{n}{
    <-6>  MnSymbolE5
   <6-7>  MnSymbolE6
   <7-8>  MnSymbolE7
   <8-9>  MnSymbolE8
   <9-10> MnSymbolE9
  <10-12> MnSymbolE10
  <12->   MnSymbolE12
}{}
\DeclareFontShape{OMX}{MnSymbolE}{b}{n}{
    <-6>  MnSymbolE-Bold5
   <6-7>  MnSymbolE-Bold6
   <7-8>  MnSymbolE-Bold7
   <8-9>  MnSymbolE-Bold8
   <9-10> MnSymbolE-Bold9
  <10-12> MnSymbolE-Bold10
  <12->   MnSymbolE-Bold12
}{}
\DeclareMathDelimiter{\llangle}{\mathopen}%
{MnLargeSymbols}{'164}{MnLargeSymbols}{'164}
\DeclareMathDelimiter{\rrangle}{\mathclose}%
{MnLargeSymbols}{'171}{MnLargeSymbols}{'171}

\newcommand{\iinner}[2]{\left\llangle#1,#2\right\rrangle}
\newcommand{\eqcls}[1]{\left[#1\right]}

\NewDocumentCommand{\idxup}{ %% operator for raising indices via a
  m                  %% #1 the differential form whose indices to be raised
  O{\defaultMetric}  %% #2 the hermitian metric on X
  t{,}               %% #3 separator
  o                  %% #4 extra superscripts
  s                  %% #5 smash the vertical spacing on the top of the metric if present
  t{.}               %% #6 with contraction operator \ctrt if '.'-ed
}{\paren{#1}^{
    \!\IfBooleanTF{#5}{\smash[t]{#2}}{#2}\IfNoValueF{#4}{, #4}
  }\IfBooleanT{#6}{\!\!\ctrt}}
\newcommand{\dbadj}{\dbar^{\smash{\mathrlap{*}\;\:}}}

\NewDocumentCommand{\sm}{s m}{{#2}\IfBooleanTF{#1}{_}{^}\text{sm}}

\newcommand{\tloc}{{\text{loc}}}

\NewDocumentCommand{\idx}{ %% multi-indices
  O{i} %% #1 symbol of the indices
  m    %% #2 starting subscript
  o    %% #3 additional stuff to add before \dotsm
  t{.} %% #4 display "\dotsm" if '.'-ed
  t{,} %% #5 display ",\dots," if ','-ed
  o    %% #6 additional stuff to add after \dotsm
  m    %% #7 ending subscript
}{{#1}_{#2} \IfNoValueF{#3}{#3}
  \IfBooleanT{#4}{\dotsm} \IfBooleanT{#5}{,\dots,}
  \IfNoValueF{#6}{#6} {#1}_{#7}}

\newcommand{\defaultPscript}{p}
\newcommand{\setDefaultPscript}[1]{\renewcommand{\defaultPscript}{#1}}

\NewDocumentCommand{\ps}{ %% sub-/super-script with parentheses
  E{_^}{{\defaultPscript}}
}{\IfNoValueTF{#2}{_{(#1)}}{^{(#2)}}}

\DeclareFontFamily{U} {MnSymbolC}{}
\DeclareSymbolFont{MnSyC} {U} {MnSymbolC}{m}{n}
\SetSymbolFont{MnSyC} {bold}{U} {MnSymbolC}{b}{n}
\DeclareFontShape{U}{MnSymbolC}{m}{n}{
  <-6> MnSymbolC5
  <6-7> MnSymbolC6
  <7-8> MnSymbolC7
  <8-9> MnSymbolC8
  <9-10> MnSymbolC9
  <10-12> MnSymbolC10
  <12-> MnSymbolC12}{}
\DeclareFontShape{U}{MnSymbolC}{b}{n}{
  <-6> MnSymbolC-Bold5
  <6-7> MnSymbolC-Bold6
  <7-8> MnSymbolC-Bold7
  <8-9> MnSymbolC-Bold8
  <9-10> MnSymbolC-Bold9
  <10-12> MnSymbolC-Bold10
  <12-> MnSymbolC-Bold12}{}
\DeclareMathSymbol{\smallstar}{\mathbin}{MnSyC}{"80}
\DeclareMathSymbol{\medstar}{\mathbin}{MnSyC}{"82}
\DeclareMathSymbol{\largestar}{\mathbin}{MnSyC}{"83}

\newcommand{\charfct}{\mathbbm 1}
\NewDocumentCommand{\lelong}{m O{x}}{\operatorname{\boldsymbol{\nu}}\paren{#1,#2}}

\newcommand{\cvr}[1]{\mathfrak{#1}} %% set of covering subsets
\NewDocumentCommand{\rs}{ %% putting ~ on objects on the
  s  %% when * is given, \smash[t] is applied
  m  %% the main object 
}{\IfBooleanTF{#1}{\smash[t]{\widetilde{#2}}}{\widetilde{#2}}}

\DeclareMathOperator{\mlc}{mlc} %% minimal lc centre
\newcommand{\Diff}{\operatorname{Diff}^*} %% general different (adjunction formula)

\newcommand{\sect}[1][s]{\mathtt{#1}} %% canonical section

\NewDocumentCommand{\cbn}{  %% group of combinations
  D//{\defaultlcIndex_V}
  D||{\defaultlcIndex}
}{\mathfrak{C}^{#1}_{#2}} 
\NewDocumentCommand{\Iset}{  %% index set for lc centres on log-resolution
  D||{\defaultlcIndex}    %% #1
  e{+-}                   %% #2,#3
  O{\defaultlclocus}      %% #4
  d()                     %% #5 open set on which the index set is valid
}{I^{\lcIndex{#1}{#2}{#3}}_{#4}\IfNoValueF{#5}{\paren{#5}}} 
\ifcsname defineNoThmInMarionotations\endcsname
  \relax
\else 

  \newtheorem{THMprop}{Proposition}[subsection]
  \newtheorem{THMlemma}[THMprop]{Lemma}
  \newtheorem{THMthm}[THMprop]{Theorem}
  \newtheorem{THMcor}[THMprop]{Corollary}

  \newtheorem{THMconjecture}[THMprop]{Conjecture}
  \newtheorem*{THMclaim}{Claim}

  \def\makeparenother{\catcode`\(=12 \catcode`\)=12 }
  \def\makeparenactive{\catcode`\(=\active\catcode`\)=\active}

  \makeparenactive
  \NewDocumentEnvironment{textupparenenvir}{}{
    \everymath\expandafter{\makeparenother}
    \everydisplay\expandafter{\makeparenother}

    \def({\textup{\char`\(}}
    \def){\textup{\char`\)}}

    \makeparenactive
  }{\makeparenother}
  \makeparenother

  \NewDocumentEnvironment{prop}{ +o }{
    \IfNoValueTF{#1}{\begin{THMprop}}{\begin{THMprop}[{#1}]}
      \begin{textupparenenvir}
  }{
      \end{textupparenenvir}
    \end{THMprop}
  }

  \NewDocumentEnvironment{lemma}{ +o }{
    \IfNoValueTF{#1}{\begin{THMlemma}}{\begin{THMlemma}[{#1}]}
      \begin{textupparenenvir}
  }{
      \end{textupparenenvir}
    \end{THMlemma}
  }

  \NewDocumentEnvironment{thm}{ +o }{
    \IfNoValueTF{#1}{\begin{THMthm}}{\begin{THMthm}[{#1}]}
      \begin{textupparenenvir}
  }{
      \end{textupparenenvir}
    \end{THMthm}
  }

  \NewDocumentEnvironment{cor}{ +o }{
    \IfNoValueTF{#1}{\begin{THMcor}}{\begin{THMcor}[{#1}]}
      \begin{textupparenenvir}
  }{
      \end{textupparenenvir}
    \end{THMcor}
  }

  \NewDocumentEnvironment{conjecture}{ +o }{
    \IfNoValueTF{#1}{\begin{THMconjecture}}{\begin{THMconjecture}[{#1}]}
      \begin{textupparenenvir}
  }{
      \end{textupparenenvir}
    \end{THMconjecture}
  }

  \NewDocumentEnvironment{claim}{ +o }{
    \IfNoValueTF{#1}{\begin{THMclaim}}{\begin{THMclaim}[{#1}]}
      \begin{textupparenenvir}
  }{
      \end{textupparenenvir}
    \end{THMclaim}
  }

  \theoremstyle{remark}
  \newtheorem{remark}[THMprop]{Remark}

  \theoremstyle{definition}
  \newtheorem{definition}[THMprop]{Definition}

  \numberwithin{equation}{subsection}

\fi

\allowdisplaybreaks  %% allow multi-line equations to spread across
\ifx\pdftexversion\undefined
  \renewcommand{\ibar}{{\raisebox{-0.9ex}{$\mathchar'26$}\mkern-6.7mu i}}
\fi

\setDefaultlclocus{D}
\setDefaultvphi{\vphi_F}
\setDefaultpsi{\psi_D}

\newtheorem{step}{Step}
\newcommand{\xb}[1]{\textcolor{blue}{#1}}

\newif\ifshowNewtextnpara %%% a flag to determine whether to highlight
\NewDocumentCommand{\newtext}{
  O{blue}  %% #1 colour for highlighting the text
  s        %% #2 
  +m       %% #3 the text
}{{\ifshowNewtextnpara\IfBooleanT{#2}{\color{#1}}\fi%
    #3}}

\NewDocumentEnvironment{newpara}{
  O{blue}
  s
}{\ifshowNewtextnpara\IfBooleanT{#2}{\color{#1}}\fi}{}

\begin{document}

%%%%% Aliases for the citation reference names;
%%%%% uncomment the needed one to use or add new ones.
%%%
\citealias{Amb03}{Ambro_quasi-log-var}
\citealias{Amb14}{Ambro_injectivity}
\citealias{Eno90}{Enoki}
\citealias{EV92}{Esnault&Viehweg_book}
\citealias{Fuj11}{Fujino_log-MMP}
\citealias{Fuj12b}{Fujino_vanishing-thms}
\citealias{Fuj13a}{Fujino_injectivity-II}
\citealias{Fuj13b}{Fujino_injectivity-hodge-theoretic}
\citealias{Fuj15b}{Fujino_survey}
\citealias{CCM}{Chan&Choi&Matsumura_injectivity}
\citealias{Mat22}{Matsumura_injectivity-Kaehler}
\citealias{Tak96}{Takegoshi_higher-direct-images}

%%%%%
%%%%% File name  : titleinfo.tex
%%%%% Creator    : Mario Chan
%%%%% Date       : 25th February, 2024
%%%%% Description: This file contains the info needed for maketitle
%%%%%              for the project "Relative-Fujino".
%%%%%
%%
%%%
\newcommand{\titlestr}{%
  Injectivity theorems for higher direct images \\ under proper K\"ahler morphisms on snc spaces%
  % : \\
  % an application of the theory of
  % harmonic integrals on log-canonical centers via adjoint ideal
  % sheaves%
}

\newcommand{\shorttitlestr}{%
  Injectivity theorems for higher direct images%
}

\newcommand{\MCname}{Tsz On Mario Chan}
\newcommand{\MCnameshort}{Mario Chan}
\newcommand{\MCemail}{mariochan@ntu.edu.tw}

\newcommand{\NTUAddressstr}{%
  Department of Mathematics, National Taiwan
  University, Taipei 106319, R.O.C.~Taiwan%
}

\newcommand{\YJname}{Young-Jun Choi}
\newcommand{\YJnameshort}{Young-Jun Choi}
\newcommand{\YJemail}{youngjun.choi@pusan.ac.kr}

\newcommand{\PNUAddressstr}{%
  Department of Mathematics and Institute of Mathematical Science, Pusan National
  University, 2, Busandaehak-ro 63beon-gil, Geumjeong-gu, Busan,
  46241, Republic of Korea%
}

\newcommand{\ShMname}{Shin-ichi Matsumura}
\newcommand{\ShMnameshort}{Shin-ichi Matsumura}
\newcommand{\ShMemail}{mshinichi0@gmail.com, mshinichi-math@tohoku.ac.jp}

\newcommand{\TohokuAddressstr}{%
  Mathematical Institute $\&$ Division for the Establishment of
  Frontier Science of Organization for Advanced Studies,
  Tohoku University, 6-3, Aramaki Aza-Aoba,
  Aoba-ku, Sendai 980-8578, Japan%
}

\newcommand{\subjclassstr}[1][,]{%
  32J25 (primary)#1  %% Transcendental methods of algebraic geometry (complex-analytic aspects) 
  32Q15#1   %% 	Kähler manifolds
  14B05 (secondary)%   %% Singularities in algebraic geometry
  % 14E30 (secondary)%   %% Minimal model program (Mori theory, extremal rays)
}

\newcommand{\keywordstr}[1][,]{%
  $L^2$ injectivity#1
  adjoint ideal sheaf#1
  multiplier ideal sheaf#1
  log-canonical center%
}

\newcommand{\dedicatorystr}{%
}

\newcommand{\thankstr}{%
}

%%% Local Variables:
%%% mode: latex
%%% TeX-master: "Relative-Fujino"
%%% coding: utf-8
%%% End:

\title[\shorttitlestr]{\titlestr}
 
\author[\MCnameshort]{\MCname}
\email{\MCemail}
\address{\NTUAddressstr}

\author{\YJname}
\email{\YJemail}
\address{\PNUAddressstr}

\author{\ShMname}
\email{\ShMemail}
\address{\TohokuAddressstr}

% \thanks{\thankstr}
 
\subjclass[2020]{\subjclassstr}


% \dedicatory{\dedicatorystr}

\begin{abstract}
  
% Let $X$ be a complex manifold, $Y$ and $D$ be two reduced snc divisors
% on $X$ with no common irreducible components.
% Given a proper locally K\"ahler morphism $\pi \colon X \to \Delta$
% from $X$ to a complex analytic space $\Delta$, we prove a more general
% form of Fujino's conjecture in the relative setting, i.e.~an
% injectivity theorem for higher direct images under $\pi$, for the lc
% pair $(X,D)$ as well as $(Y,D_Y)$ (where $D_Y := D \cap Y$).
% The main technique in the proof consists of the theory of harmonic
% integrals together with the residue formulae associated to the adjoint
% ideal sheaves, which we have used already in our previous work to
% handle the case in the absolute setting (where $\Delta$ is a point and
% $X$ is compact).
% We also make use of the Takegoshi harmonic forms in order to deal with
% the non-compactness of $X$.

% \bigskip 

% \xb{
Let $X$ be a complex manifold, and let $Y$ and $D$ be two reduced simple-normal-crossing (snc) divisors
on $X$ with no common irreducible components.
Given a proper locally K\"ahler morphism $\pi \colon X \to \Delta$
from $X$ to a complex analytic space $\Delta$, 
we prove Fujino's conjecture on the injectivity theorem in the relative setting in a generalized form. 
Specifically, we establish an injectivity result for the higher direct images under $\pi$ 
for the lc pairs $(X, D)$  as well as  $(Y, D_Y)$, where $D_Y := D \cap Y$. 
As an application, this result immediately implies the injectivity theorem on holomorphically convex K\"ahler manifolds with reduced snc divisors.
The main technique in the proof consists of the theory of harmonic
integrals together with  residue formulae associated with adjoint
ideal sheaves, which are developed from our previous work for the absolute case (where $\Delta$ is a point and $X$ is compact). Additionally, we make use of the Takegoshi harmonic forms to deal with the non-compactness of $X$.
% }

\end{abstract} 

\date{\today} 

\maketitle

%%%%% End of Top matter %%%%%%%%%%

\section{Introduction}\label{sec:intro}

%\xb{SM: the texts in blue are personal memos; do not remove them.}
%\xb{We use $Y$ for singular snc space and use $X$ for a manifold.}
%\xb{24/07:All the statements of theorems and definitions are changed to the case of global embedding. }
%\xb{24/07:I think that the formulation with singular weights should not be written in Section 1, but in Section 3. 
%I am a little confused about the restriction of multipliear ideals in
%the formulation of the theorem; do we restrict the multipliear ideals
%to $Y$, or do we restrict the metric to $Y$ and consider the
%multipliear ideals (in the latter case, it needs to be explained that
%$Y$ is a snc.)}

%\alert{MARIO: I guess you mean we should assume $\vphi_F$ and
% $\vphi_M$ to be smooth when we consider the injectivity on $Y$.
%  I agree on that, especially if we don't want
%  to clarify what a multiplier ideal sheaf on an snc space is (we have
%  to clarify what the measure is around the singularites).
%  There might be a canonical choice when the space $Y$ is snc, but I
 % believe people want the multiplier ideal sheaf to have certain
 % ``birationally invariant properties'' and the choice of the measure
 % should not be too arbitrary.  
%}

{
  \let\thesubsection\thesection
  
  This paper establishes injectivity theorems for higher direct images
  under proper \emph{locally K\"ahler morphisms} on simple-normal-crossing (snc)
  K\"ahler spaces by further developing the analytic techniques for handling
  log-canonical (lc) strata, including the theory of harmonic integrals,
  analytic adjoint ideal sheaves and the associated residue techniques,
  introduced in \cite{CCM}.
  The main results of this paper (Theorems \ref{thm:main-log-smooth} and
  \ref{main-thm} below) generalize our previous work \cite{CCM} on the
  injectivity theorems in \textit{the absolute setting} (for cohomology
  groups on compact spaces) to \textit{the relative setting} (for higher
  direct image sheaves under proper morphisms).
  In fact, we resolve a more general form of Fujino's conjecture (see
  \cite[Problem 1.8]{Fuj13a}; also cf.~\cite[Conj.~2.21]{Fuj15b} for
  the absolute setting), which is stated only for projective morphisms
  on complex manifolds.
  
    \mmark{}{\xb{Yes, \cite[Conj.~2.21]{Fuj15b} is the conjecture on compact spaces, 
    but I cited this here. You may remove this. 
    \cite[Problem 1.8]{Fuj13a} is the relative version.}}

  We remark that, %\mmark{
    % the proof delves into the complex analytic aspects of 
    % Ambro and Fujino's injectivity theorem using the theory of mixed Hodge
    % structure
    our proof delves into the complex analytic aspects of the problem,
    which is in parallel with the mixed-Hodge-structure techniques used
    in Ambro's and Fujino's injectivity theorems
 % }{Not sure what the original sentence means, so this amendment may
 %   have distorted the meaning.}
    (see \cite{Amb03, Amb14, Fuj11, Fuj12b,
    Fuj13b}). 
    \mmark{}{
    \xb{No problem with Mario's rewrite: I just wanted to say here that
    our proof studies an analytic aspect  of Ambro and Fujino's injectivity theorem, 
     which was proved using mixed-Hodge-structure techniques.}}%
  The overlapping of the two different approaches demands for further
  study and comparison, which we are regretfully not able to pursue in
  this paper.

  % \mmark{}{\xb{I found that our requirement for $\pi$ should be called {\bf{locally}} K\"ahler morphisms.}}
  \newtext{We use the notation $\ibar := \ibardefn \;$\ibarfootnote\
    in this paper.}
  Our main results are stated as follows. 
  \begin{thm}[{Injectivity for higher direct images -- the
      smooth case; Theorems \ref{T:prime_divisor_case} and
      \ref{thm:main-thm-in-section}}]
    % \label{main-cor}
    \label{thm:main-log-smooth}
    Let $(X, D)$ be a log smooth and lc pair (i.e.~a pair of a complex
    manifold $X$ and a reduced snc divisor $D$ on $X$)
    and $\pi \colon X \to \Delta$ be a proper locally K\"ahler morphism 
    to a (not necessarily irreducible or reduced) analytic space $\Delta$. 
    Let $F$ (resp.~$M$) be a line bundle on $X$ 
    with a smooth Hermitian metric $h_{F}$  (resp.~$h_{M}$) 
    such that 
    \begin{equation*}
      \ibar\Theta_{h_F}(F)\geq 0   \quad\text{and}\quad
      \ibar\Theta_{h_M}(M) \leq C \ibar\Theta_{h_F}(F)
      \quad\text{ for some } C>0 \; .
    \end{equation*}
    Consider a section $s \in H^{0}(X, M)$ whose zero locus $s^{-1}(0)$
    contains no lc \mmark{centers}{Since ``lc centers'' is used throughout
      the proof, it seems better to use it also in the main statements to
      avoid unnecessary confusion.} of the pair $(X,D)$.  
    Then, the multiplication map 
    induced by the tensor product with $s$ between the higher direct image sheaves
    \begin{equation*}
      R^{q}\pi_{*}\paren{K_{X} \otimes D \otimes F }
      \xrightarrow{\otimes s} 
      R^{q}\pi_{*}(K_{X} \otimes D \otimes F \otimes M )
    \end{equation*}
    is injective for every $q \geq 0$.
  \end{thm}

  \begin{thm}[{Injectivity for higher direct images
      -- the snc case; Theorem \ref{thm:reduction-to-log-smooth}}]\label{main-thm}
    % Let $(Y, D_Y)$ be an snc  pair globally embedded into a complex manifold $X$ and 
    % $\pi \colon X \to \Delta$ be a proper locally K\"ahler morphism to 
    % % \mmark{}{Where do we need $\pi$ be onto?}
    % % \xb{You are right; the subjectivity is not necessarily. I changed 'onto' to 'to'}
    % a  analytic space $\Delta$. 
    % Let $F$ (resp.~$M$) be a  line bundle on $X$ 
    % with a smooth Hermitian metric $h_{F}$  (resp.~$h_{M}$) 
    % such that 
    % \begin{equation*}
    %   \ibar\Theta_{h_F}(F)\geq 0   \quad\text{and}\quad
    %   \ibar\Theta_{h_M}(M) \leq C \ibar\Theta_{h_F}(F) \quad\text{ for
    %   some } C>0 \; .  
    % \end{equation*}
    % % satisying the following$:$
    % % \begin{itemize}
    % % \item[$\bullet$] $ \ibar\Theta_{h_F}(F)\geq 0   \text{ and } 
    % %   \sqrt{-1}(\Theta_{h_F}(F)-t \Theta 
    % %   _{h_M}(M))\geq 0
    % %   -C\omega \leq
    % %   \ibar\Theta_{h_M}(M) \leq C \ibar\Theta_{h_F}(F) \text{ holds for some } C>0 $
    % % \item[$\bullet$] $h_{F}$ and $h_{M}$ have only neat analytic singularities 
    % %   such that their singular loci $P_F$ and $P_M$ (with reduced structure) are divisors on $X$ 
    % %   and $(Y, D+P_F+P_M)$ is an snc  pair globally embedded into  $X$. 
    % % \end{itemize}
    Under the same notation and assumptions in Theorem
    \ref{thm:main-log-smooth}, consider a reduced snc divisor $Y$ in $X$
    such that $Y$ and $D$ \newtext{have} no common components and $Y+D$ has only snc.
    Suppose that the zero locus $s^{-1}(0)$ of the section $s \in H^{0}(X, M)$
    contains no lc centers of the pair $(X, Y +D)$.
    Let $D_Y := D \cap Y$, $F_Y := \res F_Y$, $M_Y := \res M_Y$, $s_Y :=
    \res s_Y$ and $\pi_Y := \res\pi_Y$.
    % and the norm $|s|_{h_M}$ is locally bounded. 
    Then, the multiplication map 
    induced by the tensor product with $s_Y$ between the higher direct image sheaves
    \begin{equation*}
      {R^{q}\pi_Y}_{*}\paren{ K_{Y} \otimes D_Y \otimes F_Y }
      \xrightarrow{\otimes s_Y} 
      {R^{q}\pi_Y}_{*}(K_{Y} \otimes D_Y \otimes F_Y \otimes M_Y )
    \end{equation*}
    is injective for every $q \geq 0$, where $K_Y := \parres{K_X \otimes Y}_Y$. 
    % where $\pi_Y\colon Y \to \Delta$ is the restriction of $\pi \colon X \to \Delta$ to $Y$. 
  \end{thm}

  \mmark{}{Following previous choice of Shin-ichi, I refer most of
    the assumptions of Thm.~\ref{main-thm} (also Cor.~\ref{cor:main}) to
    Thm.~\ref{thm:main-log-smooth}. I also rewrite it in such a way that
    we don't have to specifically define ``snc pair globally embedded in
    a manifold''.}

  % \mmark{}{I put a mark here to remind ourselves that this should be
  % changed to a ``Theorem'' instead of a ``Corollary'', and everywhere
  % in the text referring to this statement should also be changed.}
  % \mmark{}{\xb{I changed Cor to Thm, but the label is still using cor.}}

  % \begin{thm}\label{main-thm}
  %   Let $(X, D)$ be an snc  pair and 
  %   let $\pi \colon X \to \Delta$ be a proper locally K\"ahler morphism to 
  %   a (not necessarily irreducible or reduced) analytic space $\Delta$. 
  %   Let $F$ (resp.~$M$) be a line bundle on $X$ 
  %   with a smooth Hermitian metric $h_{F}$  (resp.~$h_{M}$) 
  %   such that 
  %   \begin{equation*}
  %     \ibar\Theta_{h_F}(F)\geq 0 \quad  \text{ and } \quad
  %     \sqrt{-1}(\Theta_{h_F}(F)-t \Theta 
  %     _{h_M}(M))\geq 0
  %     -C\omega \leq
  %     \ibar\Theta_{h_M}(M) \leq C \ibar\Theta_{h_F}(F)
  %     \quad \text{ for some } C>0 \; . 
  %   \end{equation*}
  %   Consider a section $s \in H^{0}(X, M)$  whose zero locus $s^{-1}(0)$ contains no lc centers of the pair $(X,D)$. 
  %   Then, the multiplication map  induced by the tensor product with $s$ between the higher direct image sheaves
  %   \begin{equation*}
  %     R^{q}\pi_{*}\paren{X, K_{X} \otimes D \otimes F}
  %     \xrightarrow{\otimes s} 
  %     R^{q}\pi_{*}(X, K_{X} \otimes D \otimes F \otimes M )
  %   \end{equation*}
  %   is injective for every $q$.
  % \end{thm}

  % \mmark{}{\xb{I have added a new corollary below.}}

  The above results can also be seen as a generalization of results in
  \cite{CCM} on compact spaces to holomorphically convex spaces.
  \mmark[blue]{%
    Any holomorphically convex manifold $X$ 
    admits a proper surjective map $\pi \colon X \to \Delta$ to a Stein space.
    By the Leray spectral sequence and Cartan's Theorem B for Stein spaces,
    % \mmark{}{While this name is aptly referring
    % to the statement below according to Hirzebruch (Topo.~Method ...),
    % it refers to the thm. that \v Cech cohomology can be computed with
    % acyclic cover according to Hartshorne, Grauert--Remmert (Theory of
    % Stein spaces) and Wikipedia.
    % Perhaps have to call this the result of the Leray spectral sequence
    % and Cartan(--Serre)'s Thm. B. 
    % }
    we have the isomorphism 
    \begin{equation*}
      H^q(X, \mathcal{F}) \cong H^0(\Delta, R^q \pi_*\mathcal{F})
    \end{equation*}
    for any coherent sheaf $\mathcal{F}$.  
    The functor $H^0(\Delta, \bullet)$ of global sections on $\Delta$ is left exact, 
    so we infer the following injectivity theorem on holomorphically
    convex manifolds, as well as their snc divisors, 
    from Theorems \ref{thm:main-log-smooth} and \ref{main-thm}.
  }{I have added a new corollary below.}
  \begin{cor}[Injectivity on holomorphically convex K\"ahler spaces]\label{cor:main}
    \mhlight[blue]{
      Let $X$ be a holomorphically convex K\"ahler manifold. 
      % Let $F$ (resp.~$M$) be a  line bundle on $X$ 
      % with a smooth Hermitian metric $h_{F}$  (resp.~$h_{M}$) 
      % such that 
      % \begin{equation*}
      %   \ibar\Theta_{h_F}(F)\geq 0   \quad\text{and}\quad
      %   \ibar\Theta_{h_M}(M) \leq C \ibar\Theta_{h_F}(F) \quad\text{ for
      %   some } C>0 \; .  
      % \end{equation*}
      % Consider a section $s \in H^{0}(X, M)$ whose zero locus $s^{-1}(0)$ contains no lc centers of the pair $(X,D)$.  
      % Then,
      Under the same notation and assumptions on $Y$, $D$, $F$, $M$ and $s$ in Theorems
      \ref{thm:main-log-smooth} and \ref{main-thm},
      the multiplication maps induced by the tensor products with $s$ and $s_Y$ between
      the cohomology groups
      \begin{align*}
        H^q(X, K_{X} \otimes D \otimes F )
        &\xrightarrow{\otimes s} 
          H^q(X, K_{X} \otimes D \otimes F \otimes M )
          \mathrlap{\quad\text{ and}}
        \\
        \cohgp q[Y]{\logKX<Y>[D_Y \otimes F_Y]}
        &\xrightarrow{\otimes s_Y}
          \cohgp q[Y]{\logKX<Y>[D_Y \otimes F_Y \otimes M_Y]}
      \end{align*}
      are injective for every $q \geq 0$.
    }
  \end{cor}

  All these results can be extended to the case where $h_F$ and $h_M$
  are singular (with restricted assumptions on the singularities;
  see Section \ref{subsec:notation}), by handling the singularities
  using the techniques established already in \cite{Chan&Choi_injectivity-I}.
  The statements will then incorporate the multiplier ideal sheaves of
  the singular Hermitian metrics. 
  % However, these extensions are somewhat technical and are therefore discussed 
  These more general statements are stated and proved in Section
  \ref{sec:main_results_proofs}.
  Note that we have to use Theorem \ref{thm:main-log-smooth} with
  singular Hermitian metrics to prove even Theorem \ref{main-thm} with
  smooth Hermitian metrics.
  See the proof of Theorem \ref{thm:reduction-to-log-smooth} for details.

  % The definitions of \mmark{snc pairs, lc strata,}{Maybe not needed?} and locally K\"ahler morphisms, which appear in the theorems, are given as follows:
  The definition of a locally K\"ahler morphism on a complex manifold
  $X$, compatible with the more general \cite{Takegoshi_higher-direct-images}*{Def.~6.1}, is given as follows.

  \mmark{
    % \begin{definition}[snc pairs and lc strata]\label{def-snc}
    %   Let $Y$ be a (not necessarily irreducible) reduced analytic space and $D=\sum_{i \in I_D}D_{i}$ be a Weil divisor on $X$. 
    %   The pair $(Y,D)$ is called a \textit{simple-normal-crossing (snc) pair} globally embedded into a complex manifold $X$ 
    %   if there exists an snc divisor $\bar{Y} + \bar{D}$ on $X$ such that $\bar{Y}|_{X}=Y$ and $\bar{D}|_{Y}=D$. 
    %   An {\textit{lc stratum}} of the snc pair $(Y,D)$ is the restriction to $Y$ of an lc stratum of $(X, \bar{Y} + \bar{D})$ 
    %   (i.e.~a connected component of the non-empty intersections of the irreducible components of $\bar{Y} + \bar{D}$).
    % \end{definition}
  }{The rewritten Thm.~\ref{main-thm} saves us from the need of the
    definition of ``snc globally embedded''. Also, the lc centers are defined in \S
    \ref{subsec:notation}. I've removed them, but can bring them back
    easily if needed.}

  \begin{definition}[locally K\"ahler morphisms on smooth $X$;
    cf.~\cite{Takegoshi_higher-direct-images}*{Def.~6.1}] \label{def-Kahler}
    A holomorphic map $\pi\colon X \to \Delta$ from a complex manifold $X$ to an analytic space $\Delta$ 
    is said to be a \mmark{\textit{locally K\"ahler morphism}}{
      I've added back ``locally'' here. Let me know if I shouldn't
      have to.
      \\
      And, we are following Takegoshi's definition of ``locally
      K\"ahler'', right? I've added a reference to that.
    } if every point in $\Delta$ 
    admits an open neighborhood $U$ whose inverse image $\pi^{-1}(U)$ is a K\"ahler manifold.
    % for a sufficiently small neighborhood $U$ of every point in $\Delta$, 
    % the inverse image $\pi^{-1}(U)$ admits a K\"ahler form. 
  \end{definition}

  We briefly describe the development of the injectivity theorems in the relative setting prior to our study. 
  % compare the results in this paper with previous studies,
  % focusing specifically on the relative setting.
  For the absolute setting, see
  \cite[Sec.~1]{Chan&Choi&Matsumura_injectivity} for a concise
  history. 
  Consider the special case of Theorem \ref{thm:main-log-smooth}
  formulated in the framework of algebraic geometry: $\pi\colon X \to
  \Delta$ being a projective morphism, $F$ a semi-ample line bundle,
  and $M$ a positive multiple of $F$. 
  Under the assumption $D=0$, this special case corresponds to
  Koll\'ar's injectivity theorem \cite{Kollar_injectivity}. 
  Ambro and Fujino generalize Koll\'ar's theorem to the case $D \neq
  0$ using the theory of mixed Hodge structures.
  (See \cite{Takumi_Murayama} for a recent advancement in this
  direction and \cite{Junchao_Chen} for a slightly different
  approach.) 
  On the other hand, building upon Enoki's work \cite{Enoki} which
  makes use of the theory of harmonic integrals in the absolute
  setting, Takegoshi extends Koll\'ar's theorem (with $D = 0$) to the
  complex analytic setup, where $\pi\colon X \to \Delta$ is a locally
  K\"ahler morphism and $F$ is a semi-positive line bundle.
  Fujino proposes a conjecture (in both the absolute
  \cite[Conj.~2.21]{Fuj15b} and relative \cite[Problem 1.8]{Fuj13a}
  settings) that generalizes Enoki's and Takegoshi's results in the
  analytic setup to the case $D \neq 0$ which contains the results of
  Ambro and Fujino in the algebraic setup.
  In the absolute setting, partial results on Fujino's conjecture can
  be found in \cite{Matsumura_harmonic, Matsumura_injectivity-survey,
    Fujino&Matsumura, Cao&Demailly&Matsumura, Gongyo&Matsumura} for
  the klt case and \cite{Matsumura_injectivity-lc,
    Chan&Choi_injectivity-I, Chan&Choi_injectivity-proceedings} for
  the plt case.
  % This development is particularly significant as it elucidates the analytical aspects of mixed Hodge theory.
  Recently, solutions to Fujino's conjecture in the absolute setting
  are announced first by Cao--P\u aun \cite{Cao&Paun_LC-inj},
  and shortly afterward by us independently in \cite{Chan&Choi&Matsumura_injectivity}.
  The method used in \cite{Chan&Choi&Matsumura_injectivity} has the
  advantage of being applicable to establish the injectivity theorem
  on (singular) compact K\"ahler snc spaces.
  In this paper, we further develop our method to resolve Fujino's
  conjecture in the relative setting completely.

}

\subsection*{Outline of the proof}
\mmark{}{\xb{I have not used the blue color, but I have revised this
    subsection.I have tried to be as faithful as possible to Section
    3.1.}}

Theorem \ref{main-thm} follows from Theorem \ref{thm:main-log-smooth}
with a simple reduction argument via adjoint ideal sheaves (Theorem
\ref{thm:reduction-to-log-smooth}; see Section
\ref{sec:adjoint-ideal-n-residue} for a brief review of the adjoint
ideal sheaves). 
We present here an outline of the proof of Theorem
\ref{thm:main-log-smooth} with a comparison with the proof in
\cite{Chan&Choi&Matsumura_injectivity} (i.e.~the absolute setting,
where $X$ is compact).
For simplicity, we focus on the proof in the case where $D$ is a prime divisor,
% Although this case appears considerably simpler than the general case, 
which effectively illustrates the essential difficulties that we
encountered when studying the relative setting.
Once these difficulties are resolved, the case of general $D$ can be
handled with adjoint ideal sheaves as in \cite{Chan&Choi&Matsumura_injectivity}.
% The difficulties arising from $D$ not being a prime divisor 
% has been resolved by \cite{Chan&Choi&Matsumura_injectivity}. 
% Therefore, the general case can be proved by combining their methods.

The claim in Theorem \ref{thm:main-log-smooth} is a local statement on
$\Delta$, so we can assume that $\Delta$ is a relatively compact Stein
domain and $X = \pi^{-1}(\Delta)$ is a relatively compact
holomorphically convex K\"ahler domain by shrinking $\Delta$.
The Leray spectral sequence implies 
\begin{equation*}
  H^q(X, \mathcal{F}) \cong H^0(\Delta, R^q \pi_*\mathcal{F})
\end{equation*}
for any coherent sheaf $\mathcal{F}$.  
\mmark{}{Here I try to emphasise the additional technicalities from
  the absolute setting in the current one, at the beginning.}%
The key technical difficulty in the relative setting comes from the
non-compactness of $X$.
The Dolbeault theory on $X$ only guarantees that
cohomology classes with values in a vector bundle twisted by a
multiplier ideal sheaf can in general be represented by \emph{locally}
$L^2$ $\dbar$-closed forms (see Section \ref{sec:L2-theory} or
\cite{Matsumura_injectivity-Kaehler}*{Prop.~2.16}), while we want to
make use of harmonic forms which are \emph{globally} $L^2$.
Readers will notice that most of the modifications to the proof in
\cite{Chan&Choi&Matsumura_injectivity} made here is to ensure certain
differential forms (namely, $u$ and $\delta w$ in this paper) being
globally $L^2$.
The trick here is that, since the question at hand is a local problem,
once the locally $L^2$ differential form is fixed, we can
shrink $\Delta$, and therefore $X$, so that the form becomes globally $L^2$.
 \newtext{The Bochner-type formulae are needed to apply to harmonic forms. 
However, after such a shrinking, the initially fixed K\"ahler metric cannot remain complete. 
This naturally leads us to consider the \emph{Takegoshi harmonic spaces} \eqref{T-space} developed in
\cite{Takegoshi_higher-direct-images} (see Section~\ref{sec:Takegoshi-argument}), which fit our purposes.}

%Moreover, since the Bochner-type formulae are needed to apply to
%harmonic forms while the K\"ahler metric cannot remain complete after
%the restriction to a relatively compact subset, 
%we are naturally led
%to consider the \emph{Takegoshi harmonic spaces} \eqref{T-space} developed in
%\cite{Takegoshi_higher-direct-images} (see Section
%\ref{sec:Takegoshi-argument}), which fulfill our needs.

The proof is divided into four steps.
In Step \ref{item:pf-reduction-simple}, we consider the long exact
sequence of direct images of $\pi$ induced by the standard exact sequence 
$
0 \to K_X \to K_X \otimes D \to K_D \to 0
$: 
\begin{equation*}
\cdots \to \drR {q-1}[\pi_D](K_D \otimes F) \xrightarrow{\delta} 
\drR q(K_X\otimes F) \xrightarrow{\tau} \drR q(K_X\otimes D \otimes F) \to \drR q[\pi_D](K_D \otimes F) \to \cdots
\end{equation*}
(where $\pi_D := \res\pi_D$).
Our goal is to prove that, for any germ $\beta_t \in \drR q(K_X\otimes
D \otimes F)_t$ at an arbitrary $t \in \Delta$ with $s \beta_t = 0 \in
\drR q(K_X\otimes D \otimes F \otimes M)_t$, we have $\beta_t=0$.
From the injectivity theorem in \cite{Takegoshi_higher-direct-images}
or the more general form in \cite{Matsumura_injectivity-Kaehler} (with
multiplier ideal sheaves) on the manifold $D$, the map $\otimes \res s_D \colon \drR
q[\pi_D](K_D \otimes F) \to \drR q[\pi_D](K_D \otimes F \otimes M)$ is
injective (note that this result is included in Theorem
\ref{thm:main-thm-in-section}).
Then, by a diagram-chasing argument via the above exact sequence
together with the morphism $\otimes s$, we can find a germ $\alpha_t \in \drR
q(K_X \otimes F)_t$ such that $\beta_t = \tau(\alpha_t)$.
\mmark{}{I try to emphasise here we fix $\alpha$ before shrinking $\Delta$.}%
At this point, we fix a choice of $\alpha_t$ and shrink $\Delta$
sufficiently such that there is a section $\alpha \in \cohgp
0[\Delta]{\drR q(\logKX)} \isom \cohgp q[X]{\logKX}$ whose germ at $t$
is $\alpha_t$ and can be viewed as a Dolbeault class which is
represented by some \emph{globally} $L^2$ $\dbar$-closed form.
The problem is thus reduced from handling sections of higher direct
image sheaves to handling cohomology classes.
We can also consider the harmonic representative of $\alpha$ (by
taking the harmonic projection of a globally $L^2$ representative)
even though it does not exist for a general Dolbeault class in
$\cohgp q[X]{\logKX}$.
It remains to prove that $\alpha \in \ker\tau$ (for technical reasons,
only) on some smaller subset $X_c \Subset X$, which is some
neighborhood of $\pi^{-1}(t)$.

In Step \ref{item:pf-Take-repr-simple}, following
\cite{Chan&Choi&Matsumura_injectivity}*{Step 1 of proof of Thm.~1.2},
we seek for an \emph{``optimal''} representative of $\beta =
\tau(\alpha)$ by applying the theory of harmonic integrals to obtain a
harmonic representative $u$ of $\alpha$ and taking an orthogonal
projection to get rid of the component in $\ker\tau$, resulting in a
harmonic form $u^\perp$ ``representing'' $\beta$ (as $\beta =
\tau(u^\perp)$ by a slight abuse of notation).
The goal is then to show that $u^\perp = 0$.
Recall that $\delta$ is the connecting morphism from the long exact
sequence above.
For the absolute setting in \cite{Chan&Choi&Matsumura_injectivity},
it is the fact $u^\perp$ being orthogonal to $\im\delta
\:\paren{=\ker\tau}$ that is actually being used to prove that
$u^\perp = 0$ (more precisely, $\HRes(u^\perp) = 0$; see below).
\mmark{}{Although the original decomposition $\cohgp q[\cdot]{\cdot}
  =\im \deltaH \oplus (\im \deltaH )^\bot$ is for the compact case and
  is correct, I would like to avoid that statement, especially in the
  display math mode, to avoid confusion.}%
% $\alpha \in H^{q-1}(X, K_X \otimes F)$ on the smooth variety $X$.  
% In the absolute setting, this is achieved by 
% the $L^2$ Dolbeault isomorphism, the $L^2$ harmonic representation, 
% and the orthogonal decomposition: 
% $$
% H^q(X, K_X \otimes F) \cong \cohgp{n,q}<\dbar,\L{L^2}>[X]{F}_{h_F}
%     \cong \mathcal{H}^{n,q}(X;F)_{h_F}=\im \deltaH \oplus (\im \deltaH )^\bot. 
% $$
% where \mmark{$\delta_\mathcal{H}$ is the composition of the connecting
%   morphism $\delta$ and the harmonic projection.}{Need revision.} 
However, in the relative setting, even the images of the (classes of)
harmonic forms under $\delta$ are not guaranteed to have harmonic
representatives (as the representatives are only \emph{locally} $L^2$).
This causes some troubles even to define $u^\perp$ properly, not to
say to claim the vanishing of $u^\perp$ (or $\HRes(u^\perp)$) using
orthogonality argument in the harmonic (Hilbert) space.

For this reason, we consider the sublevel sets $X_c := \{\Phi < c\}$
of an appropriately chosen exhaustion psh function $\Phi \geq 0$ on
$X$, for $c \in (0, \infty)$ (and arrange such that $X = X_{\infty}$).
Then, for any given $c > 0$, there is a subspace $\Gamma_c$ of harmonic forms on $D_c := D
\cap X_c$ which are extendable to harmonic forms on some larger spaces
$D_{c'}$ for some $c' > c$ (see \eqref{eq:Gamma_c-simple} or
\eqref{eq:Gamma_c-L2-delta-image-subsp}), such that $\delta\Gamma_c$
contains only (classes of) \emph{globally $L^2$ forms on $X_c$} (for a
harmonic form $w$ on $D_{c'}$ with $c' > c$, the class $\delta w$ is
represented by \emph{locally $L^2$ forms on $X_{c'}$}, but their
restrictions to $X_c$ ($\Subset X_{c'}$) are \emph{globally $L^2$ on
  $X_c$}).
\mmark{}{I think to ensure restrictions to different $X_c$'s being compatible with one
  another is not a good reason for using the Takegoshi harmonic
  spaces, because restrictions of forms are always compatible among
  themselves. That's why I made a substantial changes here.}%
As mentioned before, in order to use the Bochner-type
identities, all the harmonic forms that we are dealing with here should
come from the \emph{Takegoshi harmonic spaces}
$\Harm'/\bullet/<X_c>{F;\Phi},{h_F}* = \Harm/\bullet/<X_c>{\logKX ;\Phi},{h_F}*$ \eqref{T-space} (note that the
ordinary and the Takegoshi harmonic forms coincide on $X = X_\infty$ by Theorem
\ref{thm:Takegoshi-argument} and Remark \ref{rem-T-property}; see
Proposition \ref{prop:Takegoshi-harmonic-forms} for a proof that a
Bochner-type formula is valid for Takegoshi harmonic forms on $X_c$
for every $c > 0$).
Let $\deltaH\Gamma_c$ be the image of the harmonic projection of the
space $(-1)^{q-1} \delta\Gamma_c$ (see Theorem \ref{thm:HRes-duality}
for the reason of the factor $(-1)^{q-1}$).
We then replace $\im\delta$ by $\deltaH\Gamma_c$ in the definition of
$u^\perp$, that is, we consider the orthogonal decomposition
\begin{equation*}
  \res u_{X_c} =: u^\perp + \mu \in \paren{\deltaH\Gamma_c}^\perp
  \oplus \cl{\deltaH\Gamma_c} = \Harm'<X_c>{F;\Phi},{h_F}* \; ,
\end{equation*}
where $\paren{\deltaH\Gamma_c}^\perp$ is the orthogonal complement of
$\deltaH\Gamma_c$ and $\cl{\deltaH\Gamma_c}$ is the closure in the
Takegoshi harmonic space $\Harm'<X_c>{F;\Phi},{h_F}*$.
Since $\cl{\deltaH\Gamma_c} \subset \ker\tau$ on $X_c$, it suffices to
show that $u^\perp = 0$ (thus $\beta = \tau(\alpha) = \tau(u) =
\tau(\mu) = 0$) for the rest of the proof.

Step \ref{item:expression-of-su-simple} consists of the same
computation in \cite{Chan&Choi&Matsumura_injectivity}*{Step 2 of proof
  of Thm.~1.2} or \cite{Chan&Choi&Matsumura_injectivity}*{Step II of
  proof of Thm.~3.4.1}.
Using the assumption $0 = s \beta = s\tau(u^\perp)$ and an explicit
formula between the \v Cech and Dolbeault representatives (see Section
\ref{sec:Dolbeault}), the squared norm of $su^\perp$ on $X_c$
is rewritten as an inner product on $D_c$ via a residue formula
(Proposition \ref{prop:res-formula-dbar-exact-dot-harmonic}), namely,
\begin{equation*}
  \norm{su^\perp}_{X_c}^2 = \iinner{\bullet}{s\:\HRes(u^\perp)}_{D_c}
  \; .
\end{equation*}
(Indeed, to be precise, the right-hand side should be a limit of inner products.
See Step \ref{item:express-su-in-residue-norm} of the proof of Theorem
\ref{thm:main-thm-in-section} for details.)
The form $\HRes(u^\perp)$ is computed from the Poincar\'e residue on
$D_c$ of a form derived from $u^\perp$ (see
\eqref{eq:definition-of-w}).
The proof is complete if we show that $\HRes(u^\perp) = 0$ (thus
$u^\perp = 0$).

\mmark{}{Put some emphasise to harmonic residues.}%
We name $\HRes(u^\perp)$ as the \emph{harmonic residue} of $u^\perp$,
which has appeared already in \cite{Chan&Choi&Matsumura_injectivity}
(denoted by ``$w$'' there)
and is shown to be harmonic on $D_c$ in
\cite{Chan&Choi&Matsumura_injectivity}*{Sec.~2.4} (in particular, it
is \emph{globally} $L^2$ on $D_c$) under the curvature assumption on
$h_F$.
It is also a Takegoshi harmonic form in
$\Harm/q-1/<D_c>{\logKX<D>;\Phi},{h_F}*$ by Theorem
\ref{thm:Takegoshi-harmonicity-of-HRes}.
We are giving this harmonic residue $\HRes$ a slightly more systematic
treatment in Section \ref{sec:residue-element}, as it turns out to
have a crucial role in Step \ref{item:pf:use_u-ortho-w}.
Indeed, $\HRes$ is an adjoint of $\deltaH$ (with domain $\Gamma_c$) by
Theorem \ref{thm:HRes-duality}.
Therefore, for all $w \in \Gamma_c$,
\begin{equation*}
  \iinner{w}{\HRes(u^\perp)}_{D_c} = \iinner{\deltaH w}{u^\perp}_{X_c}
  = 0 \; .
\end{equation*}
The right-hand side vanishes for $u^\perp \in
\paren{\deltaH\Gamma_c}^\perp$.
We will obtain $\HRes(u^\perp) = 0$ if we can show that
$\HRes(u^\perp) \in \cl{\Gamma_c}$ (closure of $\Gamma_c$ in
$\Harm/q-1/<D_c>{\dotsm},{h_F}*$).
This last claim is proved by using the properties of $\HRes$. 

% , we consider the \emph{harmonic residue} $w := \HRes(u^{\bot})$
% introduced in Section \ref{sec:adjoint-ideal-n-residue}. 
% Furthermore, we prove that $w = \HRes(u^{\bot})$ is an $L^2$ harmonic form with the Takegoshi property 
% and $w$ is an obstruction to $u^{\bot}$ being zero 
% by using the residue technique and the assumption $s \beta = 0$. 

% In Step 4, by noting that the harmonic residue behaves like an adjoint operator of $\delta_\mathcal{H}$ 
% (see Theorem \ref{thm:HRes-duality}),
% we  show that the $L^2$-norm $\|w\|_{D_c}$ on $D_c $ 
% can be written as 
% $$
% \|w\|^2_{D_c} = \iinner{\HRes(u^{\bot})}{\HRes(u^{\bot})}
% =\iinner{\delta_\mathcal{H} (\HRes(u^{\bot}))}{u^{\bot}}.
% $$
% Then, the right-hand side is zero by the choice of $u^{\bot}$, 
% which shows that  the desired conclusion $u^{\bot}=0$. 

\mmark{}{\xb{We must change this part if we split Section 2.}}
\mmark{}{I remove the Section ``Organization...'' to see how it
  looks. Can bring it back easily.}
% \subsection*{Organization of this paper}
This paper is organized as follows:
% In Section \ref{sec:preliminaries}, to prove the main results, 
% we develop theories,  including harmonic integrals, analytic adjoint ideal sheaves, associated residue techniques, 
% and the Takegoshi property as discussed in \cite{Tak96}. 
% The technical core of this paper is particularly concentrated 
% in Sections \ref{sec:Takegoshi-argument}, \ref{sec:adjoint-ideal-n-residue}, and \ref{sec:residue-element}.
% Section \ref{sec:main_results_proofs} is dedicated 
% to proving Theorems \ref{thm:main-log-smooth} and \ref{main-thm} 
%  in a form involving singular Hermitian metrics. 
%  Section \ref{sec:main_results_proofs} consists of three subsections. 
%  In Section \ref{sec:proof-on-X}, 
%  we focus on the case where $D$ is a prime divisor in Theorem \ref{thm:main-log-smooth}. 
%  This allows us to illustrate the essential difficulties encountered in studying the relative setting 
%  and to explain our ideas for resolving them in detail.
% Section \ref{sec:proof-on-D} provides a proof for Theorem \ref{thm:main-log-smooth}, 
% and Section \ref{sec:reduction-to-log-smooth} provides a proof for Theorem \ref{main-thm}.

\tableofcontents

\subsection*{Acknowledgment}\label{subsec-ack}
This project started to take shape during our visit to Kagoshima
University and we would like to thank our host Masaaki Murakami for
his invitation and his hospitality during our stay.
S.\,M.\,wrote this paper during his stay at Westlake University. 
He would like to thank the faculty members and his host  Xin Fu 
for providing a wonderful environment.
We would like to thank the anonymous referee for the careful reading
and the long list of suggestions which helps to improve the
readability of the paper.
S.\,M.\,is supported 
by Grant-in-Aid for Scientific Research (B) $\#$21H00976 from JSPS 
and JST FOREST Program $\#$JPMJFR2368 from JST.
Y.C.~and M.C.~were supported by the National Research Foundation
of Korea (NRF) Grant funded by the Korean government
(No.~2023R1A2C1007227).
Y.C.~was partially supported by Samsung Science and Technology Foundation under Project
Number SSTF-BA2201-01.
M.C.~was partially supported by the National Science and Technology
Council of Taiwan (NSTC) under the project number NSTC 115-2115-M-002-002-MY3.

\section{Preliminary Results}
\label{sec:preliminaries}

% \mmark{
  In the following subsections, we fix the notation and lay down the basic
  assumptions that we use throughout the paper.
  We also state the known results needed for the proofs.
% }{This is about the style of the presentation.
%   I put a headline statement here, so that we can avoid stating 
%   ``In this section, we...'' after the title of each subsection.
%   That looks to me redundant, because the title is supposed to have
%   explained already what we are doing in that section and we should
%   not repeat the same thing again in the same line.}%

\subsection{Notation, conventions and assumptions}
\label{subsec:notation}

% In this subsection, we explain the notation and conventions used throughout the paper. 
Let $\pi \colon X \to \Delta$ be a proper locally K\"ahler morphism to an analytic space $\Delta$. 
The desired conclusions in the main results are local on $\Delta$. 
Thus, by shrinking $\Delta$ from the initial space and replacing $X$ with $X := \pi^{-1}(\Delta)$, 
we assume that 
\begin{itemize}
\item $\Delta$ is a Stein space, 
\item $X$ is a holomorphically convex K\"ahler manifold, and 
\item $\Delta$ and $X$ are relatively compact domains in the initial spaces. 
\end{itemize}
For convenience, we write the boundaries of $\Delta$ and $X := \pi^{-1}(\Delta)$ 
in the initial spaces as $\bdry \Delta$ and $\bdry X := \pi^{-1}(\bdry \Delta)$. 
The following definitions and conventions will be used consistently throughout this paper.

\begin{itemize}
\item \mmark{$n$ is the dimension of $X$.}{Do we seriously need a
    separated item for this?}
\mmark{}{\xb{ This may be removed. It is my style to keep important
    information in bullet points. (It may not look good as a paper,
    but I think it is easy to read). I wrote it clearly because many
    papers write “n” for the fiber dimension. } I see. Let's keep it
  as is then.}

\item $\omega$ is \emph{complete} K\"ahler metric  on $X$. 

\item $\Phi$ is a smooth, lower-bounded (hence assuming $\Phi \geq 0$), exhaustion psh 
function on $X$  such that 
\begin{equation*}
  \sup_X \Phi=\infty  \quad \text{ and } \quad \sup_X
  \abs{d\Phi}_\omega < \infty \; .
\end{equation*}
The above-mentioned $\Phi$ and $\omega$ are constructed  as follows: 
Take a smooth exhaustion strictly psh function $\Phi_\Delta \geq 0$ on the
Stein space $\Delta$.
Replace it with $\frac{1}{C -\Phi_\Delta}$ for $C := \sup_\Delta
\Phi_\Delta$ if $C < \infty$, so that we can assume that $\sup_\Delta
\Phi_\Delta = \infty$.
% \mmark{}{I avoid the use of the symbol ``$\chi$'' here as it is used
%   for the cut-off function. Also, the choice can be made explicitly
%   easily, and so I do it.}%
% Choose $\chi \colon \mathbb{R} \to \mathbb{R}$ is a (rapidly) increasing convex function 
% such that the composition $\Phi:=\chi \circ \pi^* \Phi_\Delta$
Then $\Phi := \pi^*\Phi_\Delta$ satisfies 
the desired properties except possibly for $\sup_X \abs{d\Phi}_\omega < \infty$. 
This last property % $\sup_X \abs{d\Phi}_\omega < \infty$ 
is fulfilled by replacing $\omega $ with the new complete K\"aher
metric $\omega + \ibddbar \Phi^2$.

%and $\chi \colon \mathbb{R} \to \mathbb{R}$ is a convex  increasing function to ensure $\sup_X \Phi=\infty$. 
%The property $\sup_X \abs{d\Phi}_\omega < \infty$ is ensure 

%\mmark{Such a function $\Phi$ is constructed by pulling back a smooth strictly plurisubharmonic function on the Stein space $\Delta$ and then shrinking $X$ to a smaller subdomain to ensure that $\abs{d\Phi}_\omega < \infty$. 
%(for example, see \cite[Step 1 in Section 3.1]{Mat22} for the
%precise argument.)}{But shrinking $X$ after fixing $\omega$ would destroy the
%completeness of $\omega$. Perhaps saying ``$d\Phi$ is smooth up to
%$\bdry X$'' rather than ``shrinking $X$'' is clearer and safer.}
%\xb{Sorry; I made a mistake here. When we shrink $X$, we can not expect  $\sup_X \Phi=\infty$. 
%I revised the sentence in blue. }

\item $X_c := \set{\Phi < c}$ and $X_\infty := X$. 
Note that $X_c \Subset X_{c'}  \Subset X$ for any $0 < c< c' <
\infty$, and each $X_c$ is holomorphically convex. 
% \item $\omega$ is a K\"ahler form on $X$. 

\item $h_F  = e^{-\vphi_F}$ and $h_M = e^{-\vphi_M}$ 
are singular Hermitian metrics on $F$ and $M$ respectively, 
where $\vphi_F$ and $\vphi_M$  denote their quasi-plurisubharmonic
(quasi-psh) potentials (of their respective curvature currents). 
Assume that $\vphi_F$ and $\vphi_M$ have at worst neat analytic
singularities\footnote{\label{fn:an-sing}\newtext{%
    A quasi-psh function $\vphi$ is said to have \emph{neat analytic
      singularities} if it is locally of the form $\vphi \equiv c \:\log
    \paren{\sum_{j=1}^N \abs{g_j}^2} \mod \smooth$ for some
    holomorphic functions $g_j$ and some constant $c \in \fieldR_{\geq 0}\;$.
  }} such that the (reduced) varieties $P_F
:=\vphi_F^{-1}(-\infty)$ and $P_M :=\vphi_M^{-1}(-\infty)$ are snc
divisors, and so is $P_F + P_M$.

\item $s$ is a holomorphic section of $M$ on $X$ such that $\sup_X
  \abs s_{\vphi_M} < +\infty$.

\item $D=\sum_{i \in \Iset||}D_{i}$ is a (reduced) snc divisor on $X$ such that 
$D$ and $P_F +P_M$ have no common irreducible component and
$D+P_F+P_M$ is an snc divisor.
Note that the index set $\Iset||$ is finite.

%\xb{In my definition, the snc divisor is a divisor with coefficients $1$; 
%hence, if $D$ has connected component with $P_F\cup P_M$, 
%then $D+P_F+P_M$  is not snc. 
%In my definition, the divisor is a ''finite'' sum of hypersurfaces. 
%I'm not that picky, so you can revert back to the previous version.}
%\mmark{}{Thanks! I've been using ``snc'' for non-reduced divisors
%  elsewhere, so let me be pedantic here. Let me also emphasise the
 % finiteness of $\Iset||$ just for extra caution in the non-compact case.} 

\item $\sect_i$ is a canonical section  of the irreducible component $D_{i}$. 

\item $\sect_D := \prod_{i\in \Iset||} \sect_i$ is the canonical section of $D$. 

\item $\sigma \in \{0,1,2,\cdots, n\}$.

% \item  $\Iset$ is the set of $p:=\{i_{1}, i_{2}, \cdots, i_{\sigma}\}$ such that  
% \mmark{$\lcS:=\cap_{k=1}^{\sigma} D_{i_{k}} $}{$\cap D_{i_k}$ may have more than
% one component.} is of codimension $\sigma$. 

% \item   $\lcc' := \cup_{p \in \Iset} \lcS$ is the union of $\sigma$-lc centers $\lcS$ of  $(X,D)$

\item Let $\lcc'$ be the union of \emph{$\sigma$-lc centers} of
  $(X,D)$ indexed by $\Iset$, i.e. 
  \begin{equation*}
    \lcc' :=\bigcup_{p \in \Iset} \lcS \; ,
  \end{equation*}
where, under the assumption $(X,D)$ being log-smooth and lc, each
$\sigma$-lc center $\lcS$ is
a $\sigma$-codimensional irreducible component of an intersection
of some irreducible components $D_i$ of $D$.
Set $\lcc|0|' := X$ and let $\Iset|0|$ be a singleton for convenience.
Note that $\Iset|1| = \Iset||$ and $P_F \cup P_M$ does not contain any lc centers of $(X,D)$.

\item $\Diff_{p}D$ is the effective divisor on $\lcS$ defined by the 
adjunction formula 
\begin{equation*}
  K_{\lcS} \otimes \Diff_{p}D = \parres{K_X \otimes D}_{\lcS}. 
\end{equation*}

\item  $\sect_{(p)}:=
\smashoperator{\prod\limits_{i \in \Iset|| \colon D_i
    \not\supset \lcS}} \sect_i $. 
Note that the restriction  $\sect_{(p)} |_{\lcS}$ is a canonical section of
$\Diff_{p}D$.

\item $\phi_D :=\log\abs{\sect_D}^2$ and $\phi_{(p)}
  :=\log\abs{\sect_{(p)}}^2$ are the (psh) potentials on $D$ and
  $\Diff_p D$ induced from their respective canonical sections.
%\xb{global potential?? local potential??} \alert{MARIO: it is a
%  collection of local potentials glued together via a specific
%  transformation rule. It is used in \cite{Chan&Choi&Matsumura_injectivity} and is explained in
% \cite{Chan&Choi_injectivity-I}.}
%\xb{Thank you!}

\item $\mtidlof{\vphi} := \mtidlof<X>{\vphi}$ is the multiplier ideal
  sheaf of the potential $\vphi$ on $X$ given at each $x \in X$ by
  \mmark{}{Definition of $\mtidlof{\vphi_L}$ added.}%
  \begin{equation*}
    \mtidlof{\vphi}_x  % := \mtidlof[X]{\vphi}_x
    :=\setd{f \in \holo_{X,x}}{
      \exists~\text{open set } V_x \ni x \vphantom{f^{f^f}} \; , \;
      \int_{V_x} \abs f^2 e^{-\vphi} \dvol_{V_x} < +\infty
    } \; ,
  \end{equation*}
  \newtext{where $\dvol_{V_x}$ denotes a smooth volume form.}

%  where $\dvol_{V_x}$ denotes 
%  the Lebesgue measure on the local open
%  set $V_x \Subset X$.
  Multiplier ideal sheaves $\mtidlof<S>{\vphi}$ on any submanifolds $S
  \subset X$ are defined similarly, but the subscripts will not be omitted.

\item $\cvr V := \{V_{i}\}_{i \in I}$ is a finite open cover of $X$ by 
 Stein admissible open sets induced from a locally finite cover of the initial ambient space
(see the below for the definition of admissible open sets). 

\item $\{\rho^{i}\}_{i\in I}$ is a partition of unity subordinate to
  $\cvr V =\{V_{i}\}_{i \in I}$ such that $\supp\paren{\eta \rho^i}
  \Subset V_i$ for any cut-off function $\eta$ on $X$ with $\supp \eta
  \Subset X$ (for example, the functions $\eta_\nu$ given in \eqref{eq:cut-off-functions}). 

\item $\rs\omega$ is a  complete K\"ahler metric on $X^\circ := X \setminus
\paren{P_F\cup P_M}$ defined by 
\begin{equation*}
  \rs\omega  :=2\omega +\ibddbar \frac 1{\log\abs{\ell\psi_{P_F \cup
        P_M}}} \geq \omega
  \quad\text{ so that }
\end{equation*}
\begin{itemize}[label=$\circ$]
\item $\rs\omega \geq \omega$ holds on $X^\circ$ (after choosing the
  constant $\ell \gg e$ suitably), 
\item $\rs\omega$ admits a \textit{bounded} potential locally on  $X$ (not only  $X^\circ$). 
\end{itemize}
Here $\psi_{P_F \cup P_M} \leq -1 $ is a global potential function 
defined by a canonical section of the divisor $P_F\cup P_M$.
% and   is a constant such that $\rs\omega \geq \omega$. 
%Here  $\psi_{P_F \cup P_M} := \phi_{P_F \cup P_M} -\sm\vphi_{P_F \cup P_M} \leq
%-1$ is a function on $X$ constructed from potentials of $P_F\cup P_M$, 
%where $\sm\vphi_{P_F\cup P_M}$ is any chosen smooth potential and 
%$\phi_{P_F\cup P_M}$ is the potential induced from a canonical section,
%and $\ell \gg e$ is a constant such that $\rs\omega \geq \omega$. 
Note that $\frac 1{\log\abs{\ell\psi_{P_F \cup P_M}}} $ is locally bounded on $X$ (not only $X^\circ$) 
and $\abs{d\log\paren{e \log    \abs{\ell\psi_{P_F\cup P_M}}}}_{\rs\omega} < +\infty$ 
(see, for example, \cite{Chan&Choi_injectivity-I}*{\S 2.2, item (4)}). 
Notice that $\rs\omega$ is a complete metric on $X^\circ$ but not on
$X^\circ \setminus D$ (and the same happens for its restriction to
each $\lcS* := \lcS \cap X^\circ$).

%Suppose that $X $ is a relatively compact subset in the initial ambient space $\bar X$. 
%Take open covers $ \{ \bar{U}_{i} \}_{i \in I}$ and $\{ \bar{V}_{i} \}_{i \in I}$ of $\bar X$ 
%such that $\bar{V}_{i} \Subset \bar{U}_{i} $ for any $i \in I$. 
%By $X \Subset X$, 
%we choose a finite set $I$ (which we use the same notation) such that 
%$\{ \bar{V}_{i} \}_{i \in I}$ covers $X$, 
%where $V_{i}:=\bar{V}_{i} \cap X$. 

%Fix  a partition of unity $\{\rho^{i}\}_{i\in I}$ subordinate to $ \{ \bar{U}_{i} \}_{i \in I}$
%$\{\rho^{i}\}_{i\in I}$ i
%We may assume that $I$ is a finite set by $X \Subset X$. 
%Fix a par

%For an open cover ${U_{i}}_{i \in J}$ of $X$ 
\end{itemize}
Recall from \cite{Chan&Choi&Matsumura_injectivity}*{\S 2.1} that an
open set $V \subset X$ is said to be \emph{admissible} with 
respect to $D$ if $V$ is biholomorphic to a polydisc centered at the
origin under a holomorphic coordinate system $(z_{1}, z_{2}, \cdots,
z_{n})$ such that
\begin{equation*} % \label{eq:local-expression-bphi-psi}
  D =\set{z_1 \dotsm z_{\sigma_V} =0}, \quad 
  \log r_{j}^2 < 0, \quad \text{and }
  r_j \fdiff{r_j} \psi_D >0 \text{ on } V \; , 
  % \res{\vphi_\bullet}_V = \smashoperator{\sum_{k=\sigma_V+1}^n} b_{\bullet,k}
  % \log\abs{z_k}^2 +\beta_\bullet \;\;\text{ for } \bullet= F, M \; ,
\end{equation*} 
where  $r_j := \abs{z_j}$  and $\res{\psi_D}_V := \parres{\phi_D
  -\sm\vphi_D}_V =\sum_{j=1}^{\sigma_V} \log\abs{z_j}^2
-\res{\sm\vphi_D}_V$.
(The assumption is made so that the residue computation in
\cite{Chan&Choi&Matsumura_injectivity}*{Prop.~2.3.2} is valid on any
admissible open set.)
An index $p \in \Iset$ such
that $\lcS \cap V \neq \emptyset$ is identified with a permutation
representing a choice of $\sigma$ elements from the set
$\set{1,2,\dots,\sigma_V}$. Under this identification, we have 
\begin{equation*}
  \lcS \cap V = \set{z_{p(1)} = z_{p(2)} = \dotsm = z_{p(\sigma)} = 0}
  \quad\text{ and }\quad
  \res{\sect_{(p)}}_V = z_{p(\sigma+1)} \dotsm z_{p(\sigma_V)}
\end{equation*}
(cf.~the definition of the set $\cbn$ in \cite{Chan_adjoint-ideal-nas}*{\S 3.1}).

\subsection{$L^2_\tloc$ Dolbeault cohomology and $L^2$ harmonic spaces}
\label{sec:L2-theory}
%\input{L2-Dolbeault-cohomology}

%%%%%
%%%%% File name  : L2-Dolbeault-cohomology.tex
%%%%% Author     : Mario Chan
%%%%% Date       : 4th May, 2024
%%%%% Description: 
%%%%%
%%
%%%

{
  \setDefaultvphi{\vphi_L}
  \setDefaultMetric{\rs*\omega}
  % \setDefaultAmbientSpace{X^\circ}

  \NewDocumentCommand{\decor}{
    O{\dbar}               %% #1 invisible height reference for the
                           %% super- and sub-scripts
    D//{\bullet,\bullet}   %% #2 degrees (superscript)
    E{_}{{(2)}}            %% #3 L_p type (subscript)
    t{c}                   %% #4 show X_c in subscript 
  }{{\vphantom{#1}}^{#2}_{#3\IfBooleanT{#4}{,X_c}}}
  
  % In this subsection, we review the $L^2_\tloc$ Dolbeault cohomology and the $L^2$ harmonic spaces.

  Let $(L, \vphi_L)$ denote either $(F, \vphi_F)$ or $(F\otimes M,
  \vphi_F+\vphi_M)$ for the remainder of Section \ref{sec:preliminaries}. 
  We first consider the Fr\'echet space of $L$-valued $(n,q)$-forms 
  that are locally $L^2$ on $X$ (not only on $X^\circ$) 
  with respect to $\vphi_L$ and $\rs\omega$: 
  \begin{equation*}
    \Ltwo./n,q/<X>{L} := \Ltwo./n,q/<X>{L}_{\vphilist} \; .
  \end{equation*}
  The $\dbar$-operator determines the densely defined closed operator $\dbar$ 
  with domain 
  \begin{equation*}
    \paren{\Dom\dbar}\decor/n,q/_{\tloc} :=\paren{\Dom\dbar}\decor/n,q/_{\tloc, \vphilist} 
    := \setd{
      \zeta \in \Ltwo./n,q/{L}_{\vphilist}
    }{
      \dbar\zeta \in \Ltwo./n,q+1/{L}_{\vphilist}
    } \; .
  \end{equation*}
  Let $\paren{\ker\dbar}\decor/n,q/_{\tloc}$ and
  $\paren{\im\dbar}\decor/n,q/_\tloc$ denote the kernel and image of $\dbar$, respectively, 
  where the superscript and subscript indicate that these are subspaces of $\Ltwo./n,q/<X>{L}$.
  Since $\rs\omega$ admits a bounded potential locally on $X$ (not only $X^\circ$), 
  we obtain the $L^2_\tloc$ Dolbeault isomorphism on $X$: 
  \begin{equation} \label{eq:L2-Dolbeault-isom}
    \cohgp q[X]{K_X \otimes L \otimes \mtidlof{\vphi_L}}
    \isom
    \frac{\paren{\ker\dbar}\decor/n,q/_{\tloc}}{\paren{\im\dbar}\decor/n,q/_\tloc}
    =: \cohgp{n,q}<\dbar,\Lloc[2]>[X]{L}_{\vphilist} \; ,
  \end{equation}
  where the left-hand side is treated as the \v Cech cohomology group,
  \newtext{and the completeness of $\rs\omega$ on $X^\circ$ is not
    needed to obtain the isomorphism}.
  Since $X$ is holomorphically convex, the left-hand side is Hausdorff. 
  This implies that $\paren{\im\dbar}\decor/n,q/_\tloc$ is closed in $\Ltwo./n,q/<X>{L}$. 
  See \cite{Matsumura_injectivity-Kaehler}*{\S 2.6} for details of
  the above isomorphism \newtext{and the claims on the cohomology
    groups.
    Note also that the isomorphism and claims also apply to $X_c$ for
    any $c \in (0,\infty]$.}
  % , i.e.~
  % \begin{equation*}
  %   \paren{\im\dbar}\decor/n,q/_\tloc = \cl{\paren{\im\dbar}}\decor/n,q/_\tloc
  %   \subset \Ltwo./n,q/*{L}. \; 
  % \end{equation*}

  We now consider the Hilbert space of $L$-valued $(n,q)$-forms 
  that are globally $L^2$ on $X$  with respect to $\vphi_L$ and $\rs\omega$: 
  \begin{equation*}
    \Ltwo/n,q/<X>{L} := \Ltwo/n,q/{L}_{\vphilist} \:\footnotemark
    \;\;\text{ with the $L^2$ norm }\norm{\cdot}_{X^\circ} :=
    \norm\cdot_{X^\circ,\vphilist} \; . 
  \end{equation*}%
  \footnotetext{
    The notation ``$\Ltwo/n,q/<X>{L}_{\vphilist}$'' is not
    consistent with the one in \cite{Chan&Choi_injectivity-I}.
    In \cite{Chan&Choi_injectivity-I}, this space is denoted by
    ``$\Ltwo/n,q/<X^\circ>{L}_{\vphilist}$'', which emphasizes that $\vphi_L$ and $\rs\omega$ are smooth on $X^\circ$ and
    the forms in this space can be approximated by smooth forms with compact
    support in $X^\circ$.
    Although these conditions remain unchanged in the current context,
    the forms in ``$\Ltwo./n,q/<X>{L}_{\vphilist}$'' have locally
    $L^2$ coefficients not only on $X^\circ$ but also on $X$. 
    To maintain consistency with other notations in this paper, we use
    the notation ``$\Ltwo/n,q/<X>{L}_{\vphilist}$.'' 
  }%
  Then, by \cite{Demailly}*{Ch.~VIII, \S 1},
  we have the following orthogonal decomposition:
  \begin{equation} \label{decom}
    \Ltwo/n,q/{L}= \Ltwo/n,q/{L}_{\vphilist}
    = \overbrace{\newtext{\underbrace{\ker \dbar \cap \ker \dbadj}_{=: \:
        \Harm'}}
    \oplus \cl{\paren{\im\dbar}}\decor/n,q/}^{\newtext{\ker\dbar\: =}}
    \oplus \cl{\paren{\im\dbadj}}\decor/n,q/ \; ,
  \end{equation}
  where $\dbadj$ denotes the Hilbert space adjoint
  ($\dfadj$ denotes the formal adjoint in this paper), the spaces
  $\paren{\im\dbar}\decor/n,q/$ and $\paren{\im\dbadj}\decor/n,q/$ are the images
  of $\dbar$ and $\dbadj$ in $\Ltwo/n,q/{L}$ respectively.
  The space $\Harm'$ \newtext{defined above, namely,}
  \begin{align*}
    \Harm'
    := \Harm'{L},_{\rs\omega}
    &:=\newtext{\Harm<X>} \;\footnotemark
    := \Harm{\logKX<X>[L]}
    \\
    &:= \setd{u \in
      \Ltwo/n,q/<X>{L}_{\vphilist}}{
      \begin{gathered}
        \dbar u = 0 \; , \; \dbadj u = 0
      \end{gathered}
      \;\;\text{ on } X^\circ } \; ,
  \end{align*}%
  \footnotetext{\newtext{The notation is introduced so that we don't have to
    write something like $\Harm/n-\sigma-1,q/<\lcS+1[b]>$ on a codimension-$(\sigma+1)$
    space $\lcS+1[b]$, a symbol with redundant decorations, in the latter part of this paper.}}%
  \newtext{is the space of $L^2$ harmonic $(n,q)$-forms.}
  % Thus, by \eqref{eq:L2-Dolbeault-isom}, 
  % we obtain the monomorphism
  % \begin{equation} \label{eq:map-for-harmonic-representatives}
  %   \Harm' \xhookrightarrow{\:\jmath\:}
  %   \cohgp q[X]{K_X \otimes L \otimes \mtidlof{\vphi_L}} \; .
  % \end{equation}
  % The injectivity follows from the Takegoshi property (see Proposition \ref{prop:injective-jmath}),  
  % which is consistent with \cite{Takegoshi_higher-direct-images}*{Thm.~4.3 (iv)}.
  % % (See Proposition \ref{prop:injective-jmath} for a brief proof.)
  \newtext{Note that there is the corresponding orthogonal
    decomposition of $\Ltwo/n,q/<X_c>{L}_{\vphilist}$ for any $c > 0$ as in
    \eqref{decom} with the image spaces of $\dbar$ and $\dbadj$
    denoted by $\paren{\im\dbar}\decor/n,q/c$ and
    $\paren{\im\dbadj}\decor/n,q/c$ respectively.}

  Following \cite{Takegoshi_higher-direct-images} \newtext{(in which
    the case $\vphi_L$ being smooth is handled)}, 
  we consider  the \emph{Takegoshi harmonic space} for every $c \in (0,\infty]$: 
  \begin{equation}\label{T-space}
    \begin{aligned}
      \Harm'(c)
      &:= \Harm'<X_c>{L;\Phi}
      :=\newtext{\Harm<X>(c)
      \;\addtocounter{footnote}{-1}\footnotemark
      :=\Harm<X_c>{K_X \otimes L; \Phi}}
      \\
      &:= \setd{u \in
        \Ltwo/n,q/<X_c>{L}_{\vphilist}}{
        \begin{gathered}
          \dbar u = 0 \; , \; \dfadj u = 0 \; , \\
          \idxup{\diff\Phi}.u = 0
        \end{gathered}
        \;\;\text{ on } X_c^\circ } \; ,%\footnotemark
    \end{aligned}
  \end{equation}%
  % \footnotetext{In \cite{Takegoshi_higher-direct-images}*{Thm.~4.3}, elements in
  %   $\Harm'<X_c>{L;\Phi}$ are not assumed to be $L^2$ on $X_c^\circ$,}%
  where $X_c^\circ := X_c \cap X^\circ$ and $\dfadj$ is the formal adjoint of $\dbar$ 
  with respect to the $L^2$ norm $\norm{\cdot}_{X^\circ} := \norm{\cdot}_{X^\circ, \vphilist}$.
  \newtext{Note that elements in $\Harm'(c)$ are smooth $X_c^\circ$ by the
  regularity of $\dbar$ and $\dfadj$ operators (see, for example,
  \cite{Hormander}*{Proof of Thm.~4.2.5}).
  When $c = \infty$, since} $\rs\omega$ is a complete metric on $X^\circ$, 
  the formal adjoint $\dfadj$ coincides with the Hilbert space adjoint
  $\dbadj$ \newtext{(see \cite{Demailly}*{Ch.~VIII, Thm.~(3.2)})}. 
  Furthermore, the Takegoshi property  $\idxup{\diff\Phi}.u = 0$ is automatically satisfied 
  (see \newtext{Theorem \ref{thm:Takegoshi-argument}}). 
  Thus, we have 
  $$ \newtext{\Harm'(\infty) =\Harm'<X>{L;\Phi} = \Harm'{L} = \Harm' \; .} $$ 
  \newtext{When $c < \infty$, although} $\rs\omega$ is not complete on $X_c$, 
  % the formal adjoint need not coincide with the Hilbert space adjoint a priori. 
  % Nevertheless, thanks to
  the Takegoshi property  $\idxup{\diff\Phi}.u = 0$ \newtext{helps to show that}
  the forms in $\Harm'<X_c>{L;\Phi}$ are genuine \newtext{$L^2$} harmonic forms in
  $\Ltwo/n,q/<X_c>{L}_{\vphilist}$ \newtext{such that $\Harm'(c) \subset
  \paren{\ker\dbar}\decor/n,q/_{(2), X_c} \cap
  \paren{\ker\dbadj}\decor/n,q/_{(2), X_c}$} (see \newtext{Proposition
    \ref{prop:Takegoshi-harmonic-forms}}). 

  \newtext{For every $c \in (0,\infty]$, since
    $\paren{\im\dbar}\decor/n,q/_\tloc c$ is closed in
    $\Ltwo./n,q/<X_c>{L}$, we have 
    \begin{equation*}
      \cl{\paren{\im\dbar}}\decor/n,q/c \subset
      \cl{\paren{\im\dbar}}\decor/n,q/_\tloc c
      =\paren{\im\dbar}\decor/n,q/_\tloc c \; .
    \end{equation*}
    It follows from the $\Lloc[2]$ Dolbeault isomorphism
    \eqref{eq:L2-Dolbeault-isom} and the orthogonal decomposition
    \eqref{decom} that we have a homomorphism}
  % the same way as in \eqref{eq:map-for-harmonic-representatives}, we obtain the monomorphism
  \begin{equation} \label{eq:map-for-harmonic-representatives-on-X_c}
    \Harm'(c) \xhookrightarrow{\:\jmath^c\:} \cohgp q[X_c]{\logKX[L]
      \otimes \mtidlof{\vphi_L}} \; ,
  \end{equation}
  \newtext{which is injective thanks to the Takegoshi property
  $\idxup{\diff\Phi} . u = 0$ (see Proposition \ref{prop:injective-jmath}).}
  Furthermore, the restriction from $X_{c'}$ to
  $X_c$ also induces the commutative diagram
  \begin{equation} \label{eq:restriction-maps-between-Takegoshi-harm-sp}
    \begin{aligned}
      \xymatrix{ {\Harm'(c')} \ar[d]_-{\jmath^{c'}_c}
        \ar@{^(->}[r]^-{\jmath^{c'}} \ar@{}[dr]|-*+{\circlearrowleft}
        & {\cohgp q[X_{c'}]{\logKX[L] \otimes \mtidlof{\vphi_L}}}
        \ar[d]
        \\
        {\Harm'(c)} \ar@{^(->}[r]^-{\jmath^c} & {\cohgp
          q[X_c]{\logKX[L] \otimes \mtidlof{\vphi_L}} \; }  }
    \end{aligned}
  \end{equation}
  for any $0 < c < c' \leq \infty$. 
  \newtext{We note again that these structures of $\Harm'(c)$ are shown already
  in \cite{Takegoshi_higher-direct-images}*{Thm.~4.3} for the case with
  smooth $\vphi_L$.}

  \subsection{A \v Cech--Dolbeault map with respect to a partition of unity}
  \label{sec:Dolbeault}

  % Note that the Leray theorem assures the isomorphism
  % $\cohgp q[\cvr V]{\mathcal{F}} \isom \cohgp q[X]{\mathcal{F}}$ 
  % for any coherent sheaf $\mathcal{F}$ on $X$. 

  Recall that we have the finite Stein cover $\cvr V = \set{V_i}_{i \in I}$ and the
  partition of unity $\set{\rho^i}_{i \in I}$ given in Section
  \ref{subsec:notation}.
  Consider a $\dbar$-closed form 
  \begin{equation*}
    u \in
    % \Harm' \xhookrightarrow{\jmath}
    % \cohgp q[X]{\logKX[L] \otimes \mtidlof{\vphi_L}}
    \paren{\ker{\dbar}}\decor/n,q/ \subset 
    % \Ltwo/n,q/{L}:=
    \Ltwo/n,q/{L}_{\vphilist}
    \; .
  \end{equation*}
  We can solve $\dbar$-equations \newtext{%
    \begin{alignat*}{2}
      \dbar\set{\beta_{i_0}}
      &=\set{\res{u}_{V_{i_0} \cap X^\circ}}
      &&\text{for } \beta_{i_0} \in
      \Ltwo/n,q-1/<V_{i_0} \cap X^\circ>{L}_{\vphilist}
      \quad\text{and}
      \\
      \dbar\set{\beta_{i_0 \dots i_\nu}}
      &=\underbrace{\delta\set{\beta_{i_0 \dots i_{\nu-1}}}}_{=:
        \:\set{\alpha_{i_0 \dots i_\nu}}}
      &\quad&\text{for } \beta_{i_0 \dots i_\nu} \in
      \Ltwo/n,q-\nu-1/<V_{i_0 \dots i_\nu} \cap X^\circ>{L}_{\vphilist}
    \end{alignat*}
    for $\nu = 0, \dots, q-1$ via $L^2$ method on relatively compact
    Stein subsets (see, for example, \cite{Demailly}*{Ch.~VIII, \S 6} for the
    procedures for dealing with the metric $\rs*\omega$ which is not
    complete on $V_{i_0 \dots i_\nu} \cap X^\circ$),
  }%
  % derived from $u$ with $L^2$ estimates
  % successively on various intersections of the Stein open sets in
  % $\cvr V \cap X^\circ := \set{V_i \cap X^\circ}_{i \in I}$
  to obtain a \v Cech cocycle
  \begin{equation*}
    \seq{\alpha_{\idx 0.q}}_{\idx 0,q \in I} \in 
    Z^q(\cvr V, {K_X \otimes L \otimes \mtidlof{\vphi_L}}) 
  \end{equation*}
  (see \cite{Matsumura_injectivity}*{Prop.~5.5}
  % or \cite{Chan&Choi_injectivity-I}*{Lemma 3.2.1 and Remark 3.2.2}
  for the details of the construction, \newtext{in particular for the argument
    which lets the $L^2$ holomorphic $(n,0)$-form $\alpha_{\idx 0.q}$
    extends from $V_{\idx 0.q} \cap X^\circ$ to $V_{\idx 0.q}$}).
  From the isomorphism $\cohgp q[\cvr V]{K_X\otimes L \otimes
    \mtidlof{\vphi_L}} \isom \cohgp q[X]{K_X\otimes L \otimes
    \mtidlof{\vphi_L}}$ given by Leray's theorem, together with
  the $L^2_\tloc$ Dolbeault isomorphism \eqref{eq:L2-Dolbeault-isom},
  such cocycle represents the cohomology class which corresponds to
  the $L^2_{\tloc}$ Dolbeault class of $u$ in $\cohgp q[X]{\logKX[L] \otimes
    \mtidlof{\vphi_L}}$.

  \newcommand{\cubrace}[2]{
    \begingroup
    \colorlet{currcolor}{.}
    \color{Gray}
    \underbrace{\color{currcolor}#1}_{\mathclap{#2}}
    \endgroup
  }
  
  From the above construction, we have, under the Einstein summation convention,
  \begin{alignat}{1} 
    % \begin{aligned}
    \notag
    u
    &\newtext{{}=\dbar\beta_{i_0} \qquad\paren{\text{on } V_{i_0} \cap X^\circ \;\;\;\forall~ i_0 \in I}}
    \\ \notag
    &\newtext{{}=\rho^{i_0} \dbar\beta_{i_0}} \\ \notag
    &\newtext{{}=\dbar\paren{\rho^{i_0} \beta_{i_0}}
      \cubrace{{}-\dbar\rho^{i_0} \wedge
        \beta_{i_0}}{{}=\:+\dbar\rho^{i_0} \wedge %\rho^{i_1}
        \paren{\beta_{i_1} -\beta_{i_0}} \mathrlap{\newtext{\quad\paren{\forall~ i_1 \in I}}}} 
      =\dbar\paren{\rho^{i_0} \beta_{i_0}}
      +\dbar\rho^{i_0} \wedge \rho^{i_1} \dbar\beta_{i_0 i_1}} \\  \notag
    &\newtext{{}=\dbar\paren{\rho^{i_0} \beta_{i_0}}
      +\dbar\rho^{i_0} \wedge \dbar\paren{\rho^{i_1}\beta_{i_0 i_1}}
      \cubrace{{}-\dbar\rho^{i_0} \wedge \dbar\rho^{i_1} \wedge
        \beta_{i_0 i_1}}{{}=\: +\dbar\rho^{i_1} \wedge
        \dbar\rho^{i_0} \wedge %\rho^{i_2}
        \paren{\beta_{i_0 i_1} -\beta_{i_0 i_2} +\beta_{i_1 i_2}}
        \mathrlap{\newtext{\quad\paren{\forall~ i_2 \in I}}}}}
    \\  \notag
    &\newtext{{}=\dbar\paren{\rho^{i_0} \beta_{i_0}
        -\dbar\rho^{i_0} \wedge
        \rho^{i_1}\beta_{i_0 i_1}}
      +\dbar\rho^{i_1} \wedge
      \dbar\rho^{i_0} \wedge \rho^{i_2} \dbar\beta_{i_0 i_1 i_2}} \\  \notag
    &\newtext{{}=
      \begin{aligned}[t]
        &\dbar\paren{ \rho^{i_0} \beta_{i_0}
          -\dbar\rho^{i_0} \wedge
          \rho^{i_1}\beta_{i_0 i_1}
          +\dbar\rho^{i_1} \wedge
          \dbar\rho^{i_0} \wedge  \rho^{i_2} \beta_{i_0
            i_1 i_2}} \\  \notag
        &\cubrace{{}-\dbar\rho^{i_2} \wedge \dbar\rho^{i_1}
          \wedge \dbar\rho^{i_0} \wedge \beta_{i_0 i_1 i_2}}{=
          \:+\dbar\rho^{i_2} \wedge \dbar\rho^{i_1}
          \wedge \dbar\rho^{i_0} \wedge %\rho^{i_3} \:
          \alpha_{i_0 i_1 i_2 i_3} \mathrlap{\newtext{\quad\paren{\forall~ i_3 \in I}}}}
      \end{aligned}
    }
    \\  \notag
    &=
    \begin{aligned}[t]
      &
      \begin{multlined}[t]
        \dbar\left(
          \alert{\rho^{i_0} \beta_{i_0}
            -\dbar\rho^{i_0} \wedge
            \rho^{i_1}\beta_{i_0 i_1}
            +\dbar\rho^{i_1} \wedge
            \dbar\rho^{i_0} \wedge \rho^{i_2} \beta_{i_0 i_1 i_2}}
        \right. \\ 
        \left. 
          \alert{-\dots
            +(-1)^{q-1} \:\dbar\rho^{i_{q-2}} \wedge \dots \wedge
            \dbar\rho^{i_0} 
            \cdot \rho^{i_{q-1}} \beta_{i_0\dots i_{q-1}}}
        \right)
      \end{multlined}
      \\ 
      &+\dbar\rho^{i_{q-1}} \wedge \dbar\rho^{i_{q-2}}
      \wedge \dots \wedge \dbar\rho^{i_0} \cdot %\rho^{i_q}
      \alpha_{i_0 \dots i_{q}} \newtext{\qquad\paren{\text{on }
          V_{i_q} \cap X^\circ \;\;\;\forall~ i_q \in I}}
    \end{aligned}
    \\  \notag
    &=: \dbar \alert{v_{(2)}} +\dbar \rho^{i_{q-1}} \wedge \dotsm \wedge
    \dbar\rho^{i_0} \cdot \rho^{i_q} \:\alpha_{\idx 0.q}
    \\ \label{eq:Cech-Dolbeault-isom}
    &=\dbar v_{(2)} +(-1)^q \:\underbrace{\dbar \rho^{i_q} \wedge
      \dotsm \wedge \dbar\rho^{i_1} \cdot \rho^{i_0} }_{=: \:
      \paren{\dbar\rho}^{\idx q.0}} \alpha_{\idx 0.q} \; , 
    % \end{aligned}
  \end{alignat}
  where $v_{(2)}$ is an $L$-valued $(n,q-1)$-form that is a priori
  locally $L^2$ on $X^\circ$.
  Since the cover $\cvr V$ is finite and \newtext{all $\beta$'s are globally
  $L^2$ with respect to $\vphi_L$ and $\rs*\omega$ on their respective domains}, both forms
  \begin{equation*}
    \text{
      $v_{(2)}$ and $(-1)^q \paren{\dbar\rho}^{\idx q.0}
      \alpha_{\idx 0.q}$ can be chosen to be \emph{globally $L^2$ on $X^\circ$}.
    }
  \end{equation*}
  Also note that these conclusions still hold true when $X$ is replaced by $X_c$,
  $\lcS$ or $\lcS<c>$.

  % we can trace through 
  % to obtain 

  % From the construction of , we have

  % Note that the components of are obtained by .

  % The \v Cech cocycle $\seq{\alpha_{\idx 0.q}}_{\idx 0,q \in I}$ determines the cohomology class in 
  % $\cohgp q[X]{\logKX[L] \otimes \mtidlof{\vphi_L}} $ . 

  % Notice that $\seq{\alpha_{\idx 0.q}}_{\idx 0,q \in I}$ represents
  % the class $\jmath(u)$ in $\cohgp q[X]{K_X \otimes L \otimes
  %   \mtidlof{\vphi_L}}$ for any $u \in \Harm'$. 

  % % For the application in the current article, the line bundle $L$ is
  % % either $F$ or $F \otimes M$ (or their restrictions to an lc centre
  % % of $(X,D)$ when $X$ is restricted to the lc centre).
  % Sometimes $L^2$ forms with $K_X \otimes D \otimes L$-values (or
  % $K_{\lcS} \otimes \Diff_{p} D \otimes \res L_{\lcS}$-values on an
  % lc centre $\lcS$) are considered and the (twisted)
  % Bochner--Kodaira--Nakano formula (see, for example,
  % \cite{Matsumura_injectivity-lc}*{Prop.~2.5} or
  % \cite{Chan&Choi_injectivity-I}*{Lemma 2.4.2}) for these forms is
  % needed.
  % In this case, the line bundle $D$ is equipped with the potential
  % $\phi_D$ induced from a canonical section $\sect_D$
  % % , and $X :=
  % % X \setminus Z$, with $Z := \supp \paren{P_F +
  % % P_M}$, and $\rs\omega$ retain the same meaning as before, that is,

  \section{Harmonic forms with the Takegoshi property}
  
% \xb{In the following subsections, we introduce several cut-off functions and the Takegoshi property to handle boundaries that appear when performing integration by parts. 
% These are necessary to ensure that a Bochner-type formula holds 
% for harmonic forms with the Takegoshi property.}

%   \mmark{}{Need some descriptive text for this section here.}

  \mmark{
    In order to apply the Bochner-type identities to harmonic forms on
    $X_c^\circ$ ($c < \infty$) without the completeness of the
    K\"ahler metric, we put on an extra condition on the harmonic
    forms following Takegoshi \cite{Takegoshi_higher-direct-images}.
    We provide a self-contained treatment to the use of such ``Takegoshi
    harmonic forms'', with some minor generalizations compared to
    \cite{Takegoshi_higher-direct-images} (see Remark
    \ref{rem-T-property}).
  }{Sorry that I keep changing your writing, Shin-ichi...
    \\
    Btw, I start to doubt whether ``Bochner-type identities/formulae''
    is a good way to call them. I'm imagining some people caring
    history may criticise it.}

  \subsection{Various cut-off functions}
  \label{sec:cutoff}

  We introduce several cut-off functions to handle the boundaries
  $\bdry X$ and $\bdry X_c$, and the singular loci $P_{F} \cup P_{M}$
  and $D$ in the computations of the relevant integrals. 
\newtext{This enables us to apply the Bochner--Kodaira--Nakano formula
to forms that are not necessarily smooth and do not necessarily have compact support.}

  Take a non-increasing  smooth function $\rho \colon [0,+\infty) \to
  [0,1]$ such that $\res{\rho}_{[0,\frac 12]} \equiv 1$,
  $\res\rho_{[1,+\infty)} \equiv 0$, and $\abs{\rho'} \lesssim 1$ on
  its domain (where $\rho'$ denotes the derivative).
  For $\nu \in \Nnum$ and $\eps > 0$, 
  we define $\eta_\nu$, $\eta_{c,\nu}$, $\chi_\nu$, and $\theta_\eps$ as follows:
  \begin{equation} \label{eq:cut-off-functions}
    \begin{aligned}
      \eta_\nu
      &:= \rho\paren{\frac{\Phi}{\nu}} \; ,
      \quad\;
      \eta_{c,\nu}
      :=
      \begin{cases}
        \rho\paren{\frac{1}{\nu\paren{c - \Phi}}}
        & \text{on } X_c 
        \\
        0 & \text{on } X \setminus X_c 
      \end{cases}
      \;\;\text{ for } c \in (0,\infty) \; ,
      \\
      \chi_\nu
      &:= \rho\paren{\frac{\log\paren{e \log\abs{\ell\psi_{P_F\cup P_M}}}}\nu}
      \quad\text{ and }\quad
      \theta_\eps := 1 -\rho\paren{\frac{1}{\abs{\psi_D}^\eps}} \; .
    \end{aligned}
  \end{equation}
  By the properties of $\Phi$, $\omega$, and $\tilde \omega$ given in Section \ref{subsec:notation}, 
  we can easily verify that  
  \begin{equation*}
    \begin{aligned}
      \supp \eta_\nu \cap \bdry X &= \emptyset \\
      \supp \eta_{c,\nu} \cap \bdry X_c &= \emptyset \\
      \supp \chi_\nu \cap \paren{P_F\cup P_M} &= \emptyset \\
      \supp \theta_\eps \cap D &= \emptyset
    \end{aligned}
    \quad\text{ and }\quad
    \begin{aligned}
      \eta_\nu &\ascendsto 1 &&\text{ on } X
      \\
      \eta_{c,\nu} &\ascendsto 1 &&\text{ on } X_c
      \\
      \chi_\nu &\ascendsto 1 &&\text{ on } X^\circ
      \\
      \theta_\eps &\ascendsto 1 &&\text{ on } X\setminus D
    \end{aligned}
    \quad\text{ as }\quad
    \begin{aligned}
      \nu &\ascendsto +\infty
      \\
      \eps &\descendsto 0
    \end{aligned}
  \end{equation*}
  and 
  \begin{equation*}
    \abs{d\eta_\nu}_{\rs\omega} \leq \abs{d\eta_\nu}_\omega
    \lesssim \frac 1\nu \; , \quad
    \abs{d\chi_\nu}_{\rs\omega} \lesssim \frac 1\nu
    \quad\text{ and }\quad
    d\theta_\eps = \eps \:\frac{\theta'_\eps
      \:d\psi_D}{\abs{\psi_D}^{1+\eps}} \; ,
  \end{equation*}
  where $\theta'_\eps := -\rho' \circ \frac{1}{\abs{\psi_D}^\eps} \geq
  0$ and the constants involved in $\lesssim$ are independent of
  $\nu$.
  Note that $\abs{d\eta_{c,\nu}}_{\rs\omega}$ is not uniformly bounded
  on $X_c$.

  As an example of their applications, we can make the Bochner--Kodaira--Nakano formula
  (see \cite[Prop.~2.5]{Matsumura_injectivity-lc} or \cite[Lemma 2.4.2]{Chan&Choi_injectivity-I}) 
  applicable to (a dense set of) forms in $\Ltwo/n,q/{L}_{\vphilist}$ via Friedrichs' lemma. 
   \newtext{This is often used in this paper without explicit mention. }
  Indeed, following the standard proof (for example, see \cite[Ch.VIII, Thm.~(3.2)]{Demailly}),
  we can use the cut-off functions $\eta_\nu \chi_\nu$
  to construct a sequence of compactly supported approximations of a
  given $L^2$ form before applying the smoothing kernels.
  
  With the valid Bochner--Kodaira--Nakano formula on the non-compact $X$,
  the same proof as in \cite{Enoki}, \cite{Matsumura_injectivity-lc}, or \cite{Chan&Choi_injectivity-I}
  guarantees the following result.

  \begin{prop}[\cite{Chan&Choi&Matsumura_injectivity}*{Prop.~2.2.2}]
    \label{prop:consequence-of-positivity}
    
    Suppose that
    $\ibddbar\vphi_F \geq 0$ and $u \in \Harm'{F},{\vphi_F}$.
    Then, we have
    \begin{equation*}
      \nabla^{(0,1)}u = 0 \quad\text{ and }\quad
      \idxup{\ibddbar\vphi_F} \ptinner{u}{u}_{\rs\omega} = 0
      \quad\text{ on } X^\circ \; .\footnotemark
    \end{equation*}
    \footnotetext{
      \label{fn:contract-notation}%
      Given an $(n,q)$-form $u$ and a function $\vphi$, we define
      $\idxup{\diff\vphi}.u$, $\idxup{\ibddbar\vphi}.u$ and
      $\idxup{\ibddbar\vphi}\ptinner u u_{\rs\omega}$ as raising the
      (holomorphic) indices of the coefficients (of $\diff\vphi$ or $\ibddbar\vphi$) via
      $\rs\omega$ and then contracting with the anti-holomorphic
      indices of coefficients of $u$ (and we have $\idxup{\ibddbar\vphi} \ptinner u
      u_{\rs\omega} = \inner{\idxup{\ibddbar\vphi}.u}{u}_{\rs\omega}$).
      The operator $\idxup{\diff\vphi}.\cdot$ happens to be the adjoint
      of $\dbar\vphi \wedge \cdot$ with respect to the inner product
      $\inner\cdot\cdot_{\rs\omega}$, and is sometimes denoted by
      $(\dbar\vphi)^*$ or $e(\dbar\vphi)^*$ (in
      \cite{Takegoshi_higher-direct-images}) and can be computed by
      $\pm * \diff\conj{\vphi} \wedge * \cdot$ (where $*$ is the
      Hodge-$*$ operator with respect to $\rs\omega$).
      The function $\idxup{\ibddbar\vphi}\ptinner u u_{\rs\omega}$ can
      also be denoted by $\inner{[\ibddbar\vphi, \Lambda_{\rs\omega}]
        u}{u}_{\rs\omega} = \inner{\ibddbar\vphi \Lambda_{\rs\omega}
        u}{u}_{\rs\omega}$ (for $(n,q)$-forms) (in \cite{Demailly}).
      See \cite{Chan&Choi_injectivity-I}*{Remark 2.4.3} for more
      details.
      For most of the computations in this paper, the forms of
      $(1,0)$- and $(0,1)$-types in a differential form $u$ are
      handled separately (e.g.~$u$ is mostly treated as a
      ``$K_X$-valued $(0,q)$-form'' rather than an ``$(n,q)$-form'').
      Our choice of notation is intended to make the computations more
      intuitive and avoid the unnecessary interaction between the two
      types of forms.
    }%
    Furthermore, if $\vphi_M$ 
    satisfies $\ibddbar\vphi_M \leq C\ibddbar\vphi_F$ for some constant $C > 0$,
    then we also have $su \in \Harm'{F\otimes
      M},{\vphi_F+\vphi_M}_{\rs\omega}$. 
  \end{prop}

  By using $\eta_\nu \newtext{\chi_\nu} \theta_\eps$ in place of $\theta_\eps$ in the
  proof of \cite{Chan&Choi&Matsumura_injectivity}*{Lemma 2.2.1} and
  with a little care of the order of taking the limits $\eps \tendsto
  0^+$ and $\nu \tendsto +\infty$, we also obtain the following.
  
  \begin{lemma}[\cite{Chan&Choi&Matsumura_injectivity}*{Lemma 2.2.1}]
    \label{lem:su-harmonicity}
    If $u \in \Harm'{L},{\vphi_L}_{\rs\omega}$, then $\sect_D u \in
    \Harm'{D\otimes L},{\phi_D+\vphi_L}_{\rs\omega}$. 
  \end{lemma}

  % To assure that the Bochner--Kodaira--Nakano formula is valid for
  % general $L^2$ forms (sitting inside the domains $\Dom \dbar$ and
  % $Dom \dbadj$ of $\dbar$ and its Hilbert space adjoint $\dbadj$), one
  % should in principle check that the smoothing via the convolution
  % with a smoothing kernel in Friedrichs' lemma commutes with the
  % approximation via the cut-off functions in
  % \eqref{eq:cut-off-functions} (especially the one via $\theta_\eps$)
  % under the graph norms of the relevant differential operators.
  % While this can be done, for the purpose of this paper, the
  % formula is needed only for harmonic forms 
  
  % \mmark{Under construction... Validity of Bochner--Kodaira formula}{}

}

%%% Local Variables:
%%% mode: latex
%%% TeX-master: "Relative-Fujino"
%%% coding: utf-8
%%% End:

\subsection{Takegoshi property for harmonic forms}
% \footnote{Under the assumption that the psh potential has 
%   only analytic singularities.}%
\label{sec:Takegoshi-argument}

%\input{takegoshi-argument}

%%%%%
%%%%% File name  : takegoshi-argument.tex
%%%%% Author     : Mario Chan
%%%%% Date       : 26th February, 2024
%%%%% Description: This is aiming to generalise a result of Takegoshi,
%%%%%              and therefore shows that harmonicity is a property
%%%%%              of the (Dolbeault) cohomology class of a
%%%%%              pseudoeffective line bundle.
%%%%%
%%
%%%

{ %%%%%% DefaultMetric redefined %%%%%%%%%%%%%
  % \setDefaultvphi{\vphi_L}
  \setDefaultMetric{\rs\omega}
  
  % Assume that $X$ is only a weakly $1$-complete manifold (i.e.~there
  % exists a smooth psh exhaustion function on $X$) in this section.
  % Let $L$ be a pseudoeffective holomorphic line bundle over $X$ equipped
  % with a singular hermitian metric $e^{-\vphi_L}$ such that
  % $\ibddbar\vphi_L \geq 0$ and \emph{the polar set $\vphi_L^{-1}(-\infty)$
  % lies inside some closed analytic subset $Z \subset X$}.
  % Let $\omega$ be a K\UTF{00E4}hler metric which is complete on $X^\circ :=X \setminus Z$.
  % Let $\Harm' :=\Harm'<X^\circ>{L}$, where $n :=\dim_\fieldC X$, be
  % the space of $L$-valued $L^2$ harmonic $(n,q)$-forms with respect to
  % $\vphi_L$ and $\omega$ on $X^\circ$.
  
  In this section, we introduce \textit{the Takegoshi property} and present several of its applications.
  We begin with a slight generalization of  \cite{Takegoshi_higher-direct-images}*{Thm.~3.4 (ii) and 4.3 (ii)}.

  \begin{thm}[cf.~\cite{Takegoshi_higher-direct-images}*{Thm.~3.4 (ii)
      and 4.3 (ii)}]
    \label{thm:Takegoshi-argument}
    Let $\Psi$ be a smooth function on $X$ such that
    \begin{equation*}
      \begin{gathered}
        \ibddbar\paren{C\vphi_F +\Psi} \geq 0
        \;\;\text{ and }\;\;
        \Psi > -C
      \end{gathered}
      \;\;\text{ on $X$ for some constant } C > 0 \; .
    \end{equation*}
    Suppose that $\ibddbar\vphi_F \geq 0$ holds and  
    $u \in \Harm'$ satisfies that $\norm{\idxup{\diff\Psi}.u}_{X^\circ,\vphi_F,\rs\omega}^2 < \infty$. 
    Then, the form $u$ is a harmonic form in $\Harm',{\vphi_F+\Psi}$ and satisfies that 
    \begin{equation*}
      \idxup{\diff\Psi}.u =0
      \quad\text{ and }\quad
      \idxup{\ibddbar\Psi}. u = 0
      \quad\text{ on } X^\circ \; .\addtocounter{footnote}{-1}\footnotemark
    \end{equation*}
  \end{thm}
  % \footnotetext{See footnote \ref{fn:contract-notation}.}%
  
  \begin{remark}\label{rem-T-property}
    The equations in Theorem \ref{thm:Takegoshi-argument} are referred to 
    as \emph{the Takegoshi property} in this paper.
    Note that the Takegoshi property holds for all $u \in \Harm'$ with $\Psi := \Phi$, 
    where $\Phi$ is the psh exhaustion function on $X$ given in Section \ref{subsec:notation}.
    Furthermore,
    Theorem \ref{thm:Takegoshi-argument} relaxes
    the conditions on $\Psi$ compared to 
    \cite{Takegoshi_higher-direct-images}*{Thm.~3.4 (ii) and 4.3 (ii)}.
    In our case, the function $\Psi$ is neither necessarily psh nor bounded on $X^\circ$. 
    Although such a generalization is not strictly required for the results of this paper, 
    it is of interest to investigate to what extent the statement can be generalized to any smooth function $\Psi$
    such that the harmonicity of an $F$-valued form depends only on the class $c_1(F)$.
  \end{remark}

  \begin{proof}
    Without loss of generality, assume that $C = 1$.
    Let $\dfadj$ and $\dfadj_\Psi$ be the formal adjoints of $\dbar$
    with respect to $\vphi_F$ and $\vphi_F + \Psi$, respectively. 
    By definition, we have $\dfadj_\Psi = \dfadj + \idxup{\diff\Psi} . \cdot$.
    From $u \in \Harm'$ (i.e.~$\dbar u = 0$ and $\dfadj u = 0$)
    and the assumption $\ibddbar\vphi_F \geq 0$, Proposition
    \ref{prop:consequence-of-positivity} guarantees that
    \begin{equation*} \tag{$*$} \label{eq:pf:u-1st-properties}
      \nabla^{(0,1)} u = 0
      \;\;\;\text{ and }\;\;\;
      \idxup{\ibddbar\vphi_F}\ptinner u u_{\rs\omega} = 0
      \quad\text{ on } X^\circ \; .
    \end{equation*}
    Given the assumptions on $\Psi$ and $\diff \Psi$, we also have
    \begin{equation*}
      \int_{X^\circ} \abs u_{\rs\omega}^2 \:e^{-\vphi_F-\Psi} < \infty
      \;\;\;\text{ and }\;\;\;
      \int_{X^\circ} \abs{\idxup{\diff\Psi} . u}_{\rs\omega}^2
      \:e^{-\vphi_F-\Psi}
      \lesssim \int_{X^\circ} \abs{\idxup{\diff\Psi} . u}_{\rs\omega}^2
      \:e^{-\vphi_F} < \infty \; . 
    \end{equation*}
    In particular, we have $\norm{u}_{\vphilist[\Psi]} < \infty$ and $\norm{\dfadj_\Psi
      u}_{\vphilist[\Psi]} < \infty$.

    Recall from \cite{Chan&Choi&Matsumura_injectivity}*{Lemma 2.4.2} that the formula
    \begin{equation}  \label{eq:pf:partial-BK-formula}
      \dbar\paren{\idxup{\diff\Psi}* . u}
      = \idxup{\ibddbar\Psi}. u
      - \idxup{\diff\Psi} . \cancelto{0}{\paren{\dbar u}}
      + \idxup{\diff\Psi} \cdot \cancelto{0}{\nabla^{(0,1)}_\bullet u}
    \end{equation}
    holds on $X^\circ$. Using the cut-off functions given in \eqref{eq:cut-off-functions}, we apply $\iinner{\:\cdot\:}{\eta_\nu \chi_\nu
      u}_{\vphi_F,\rs\omega}$ to both sides, integrate by parts on the
    left-hand side, and then take the limit $\nu \to
    \infty$ \newtext*{to} obtain
    \begin{equation*}
      0 = \iinner{\idxup{\diff\Psi} . u}{\dfadj u}_{\vphi_F,
        \rs\omega}
      =\int_{X^\circ} \idxup{\ibddbar\Psi}\ptinner u
      u_{\vphi_F,\rs\omega}
    \end{equation*}
    \newtext*{(notice that we use $\idxup{\ibddbar\Psi} \ptinner u
      u_{\vphilist} \overset{\text{\eqref{eq:pf:u-1st-properties}}}= \idxup{\ibddbar\paren{\vphi_F+\Psi}} \ptinner u
      u_{\vphilist} \geq 0$ from the assumption when considering the
      limit on the right-hand side)}.
    Using the sequence $\seq{\eta_\nu \chi_\nu}_{\nu\in\Nnum}$, we can verify that 
    the Bochner--Kodaira formula with respect to $\vphi_F+\Psi$ and $\rs\omega$ 
    is also valid for $u$. Thus, we have 
    \begin{align*}
      \cancelto{0}{\norm{\dbar u}_{\vphi_F+\Psi\mathrlap{,\rs\omega}}^2} \;\;\;
      +\norm{\dfadj_\Psi u}_{\vphi_F+\Psi,\rs\omega}^2
      &= \cancelto{0}{\norm{\nabla^{(0,1)} u}_{\vphi_F+\Psi\mathrlap{,\rs\omega}}^2}
        \;\;\;
        +\int_{X^\circ}
        \idxup{\ibddbar\paren{\vphi_F+\Psi}} \ptinner u
        u_{\vphi_F, \rs\omega} \:e^{-\Psi} \\
      &\lesssim \int_{X^\circ}
        \idxup{\ibddbar\paren{\vphi_F+\Psi}}\ptinner u u_{\vphi_F,
        \rs\omega}
        =
        \int_{X^\circ}
        \idxup{\ibddbar\Psi}\ptinner u u_{\vphi_F, \rs\omega} = 0 \; . 
    \end{align*}
    Therefore, we conclude that $\dfadj_\Psi u = 0$ and hence $\idxup{\diff\Psi} . u
    = 0$ on $X^\circ$ (cf.~\cite{Takegoshi_higher-direct-images}*{Thm.~4.3}).
    This further implies that $\idxup{\ibddbar\Psi}. u = 0$ on
    $X^\circ$ by \eqref{eq:pf:partial-BK-formula}.
    We can see that $u$
    lies in the domain of the Hilbert space adjoint of $\dbar$ with
    respect to $\iinner\cdot\cdot_{\vphi_F+\Psi,\rs\omega}$ 
    by the standard argument using $\seq{\eta_\nu \chi_\nu}_{\nu \in \Nnum}$
    together with Friedrichs' lemma (which ensures that compactly supported
    smooth forms are dense in the graph norm of $\dbar$; see, for example, \cite{Demailly}*{Ch.~VIII, Lemma (3.3)}). 
    Consequently,  we obtain $u \in \Harm',{\vphi_F+\Psi}$ as desired.
  \end{proof}

  Theorem \ref{thm:Takegoshi-argument} ensures that
  \begin{equation*}
    \Harm' =\Harm'{F} = \Harm'<X_\infty>{F;\Phi} = \Harm'(\infty) \; .
  \end{equation*}
  Although $\rs\omega$ is not complete on $X_c^\circ$ for $c < \infty$, 
  the Takegoshi property works as an effective substitute for the
  completeness of $\rs\omega$ in our context \mmark{(indeed, the property
  works like the boundary condition for elements in
  $\paren{\Dom\dbadj}^{n,q}_{(2), X_c}$ (\newtext{domain of the
    Hilbert space adjoint of $\dbar$}, see below), but without requiring the
  smoothness of the boundary $\bdry X_c$)}{I added this.}. 
  Specifically, it implies the following result.

  \begin{prop}[cf.~\cite{Takegoshi_higher-direct-images}*{Thm.~4.3 (i)}]
    \label{prop:Takegoshi-harmonic-forms}
    Suppose that $u \in \Harm'(c) = \Harm'<X_c>{F;\Phi}$. 
    Then,
    \begin{newpara}*%
      \begin{enumerate}[label=\roman*., ref=\emph{\roman*}, left=0pt]
      \item \label{item:prop:u-harmonic}
        the form $u $ is a genuine $L^2$ harmonic form in
        $\Ltwo/n,q/<X_c>{F}_{\vphilist}$ for every $c > 0$, that
        is,\mmark{}{Turn $\Dom \dbadj$ to $\ker \dbadj$}
        \begin{equation*}
          u \in \newtext{\paren{\ker \dbar}^{n,q}_{(2), X_c} \cap \paren{\ker \dbadj}^{n,q}_{(2), X_c}} \; ,
          % \quad\text{ and }\quad \dbar u=0 \;\;\text{ and }\;\; \dfadj
          % u=0 \; ,
        \end{equation*}
        where \newtext{$\dbadj$} is the Hilbert space adjoint of
        $\dbar \colon \Ltwo/n,q-1/<X_c>{F}_{\vphilist} \birat
        \Ltwo/n,q/<X_c>{F}_{\vphilist}\;$;

      \item \label{item:prop:u-harmonic-properties}
        $u$ satisfies the full Takegoshi property
        with $\Phi$ and the conclusions of Proposition
        \ref{prop:consequence-of-positivity},
        % In particular, we have
        \newtext{namely,}
        \begin{equation*}
          \nabla^{(0,1)}u = 0 \; , \;\;
          \idxup{\ibddbar\vphi_F} \ptinner{u}{u}_{\rs\omega} = 0
          \;\;\text{ and }\;\;
          \idxup{\ibddbar\Phi} \ptinner{u}{u}_{\rs\omega} = 0
          \;\;\text{ on } X_c^\circ = X_c \cap X^\circ \; ;
        \end{equation*}

      \item \label{item:prop:su-properties}
        the images of $u$ under the multiplications by $\sect_D$
        and $s$ satisfy the conclusions of Proposition
        \ref{prop:consequence-of-positivity} and Lemma
        \ref{lem:su-harmonicity}, namely,
        \begin{alignat*}{1}
          \sect_D u &\in \paren{\ker \dbadj_{\phi_D}}^{n,q}_{(2), X_c}
          \cap \Harm'<X_c>{D\otimes F; \Phi},{\phi_D+\vphi_F}
          \quad\text{ and}
          \\
          s u &\in \paren{\ker \dbadj_{\vphi_M}}^{n,q}_{(2), X_c} \cap
          \Harm'<X_c>{F\otimes M; \Phi},{\vphi_F +\vphi_M} \; ,
        \end{alignat*}
        where $\dbadj_{\phi_D}$ and $\dbadj_{\vphi_M}$ are the Hilbert
        space adjoints of $\dbar$ with respect to the inner products
        $\iinner\cdot\cdot_{X_c^\circ, \vphilist.}$ and
        $\iinner\cdot\cdot_{X_c^\circ, \vphilist M}$ respectively.
      \end{enumerate}
    \end{newpara}
  \end{prop}

  \begin{proof}
    % Recall that $\eta_{c,\nu} = \rho\paren{\frac{1}{\nu\paren{c -\Phi}}}$ 
    % is the cut-off function defined in    .
    Notice that the Takegoshi property (in the definition of $\Harm'(c)$) ensures
    that $\idxup{\diff\eta_{c,\nu}}. u = 0$ on $X_c^\circ$ by the
    definition of $\eta_{c,\nu}$ in \eqref{eq:cut-off-functions}.
    % Then,
    \begin{enumerate}[label=(\emph{\roman*}), ref=\emph{\roman*},
      left=0pt..0pt,
      % leftmargin=0pt,
      itemindent=*,
      % align=left,
      % labelsep=*,
      %widest=iii,
      align=left, itemsep=1.5ex]
    \item[(\ref{item:prop:u-harmonic})] \newtext{With the cut-off function $\chi_{\nu}$ in
        \eqref{eq:cut-off-functions},} for any
      $\xi \in \paren{\Dom\dbar}^{n,q-1}_{(2), X_c}$, we have
      \begin{equation} \label{eq:in-dom-dbadj}
        \begin{aligned}
          \iinner{u}{\dbar\xi}_{X_c^\circ, \vphilist} \xleftarrow{\nu
            \tendsto +\infty} &~\iinner{\eta_{c,\nu}
            \newtext{\chi_\nu} u}{\dbar\xi}_{X_c^\circ, \vphilist}
          \\
          =&
          \begin{aligned}[t]
            &-\cancelto{ 0 \;\;(\because \text{ Takegoshi property})
            }{ \iinner{ \newtext{\chi_\nu} \idxup{\diff\eta_{c,\nu}}
                . u }{\xi} }_{X_c^\circ, \vphilist}
            \newtext{~-\iinner{\eta_{c,\nu}
                \idxup{\diff\chi_\nu}. u}{\xi}_{X_c^\circ, \vphilist}}
            \\
            &+\iinner{ \eta_{c,\nu} \newtext{\chi_\nu} \dfadj u
            }{\xi}_{X_c^\circ, \vphilist}
          \end{aligned}
          \\
          \xrightarrow{\nu \tendsto +\infty} &~\iinner{\dfadj
            u}{\xi}_{X_c^\circ, \vphilist} \qquad\text{(a bounded
            functional in $\xi$)} \; .
        \end{aligned}
      \end{equation}
      \newtext{Note that the integration by parts above is valid as
        $u$ is smooth on $X_c^\circ$ (by the regularity of $\dbar$ and
        $\dfadj$) and thus $\eta_{c,\nu} \chi_\nu u$ is smooth with
        compact support in $X_c^\circ$ for any $\nu \in \Nnum$, and
        the convergence on the last step is making use of the estimate
        $\abs{\diff\chi_\nu}_{\rs\omega} \lesssim \frac 1\nu$.}  This
      implies that
      $u \in \newtext{\paren{\ker \dbadj}^{n,q}_{(2),X_c} \subset}
      \paren{\Dom\dbadj}^{n,q}_{(2),X_c}$ \newtext{(note that
        $\dfadj u = 0$ by assumption).  The first claim thus follows.}

    \item[(\ref{item:prop:u-harmonic-properties})]
      \setDefaultPscript{\nu,\nu'}

      \newtext{Other formulae on $u$} follow since the
      Bochner--Kodaira--Nakano formula on $X_c$ with respect to
      $\vphi_F+\Phi$ and $\rs\omega$ is valid for harmonic forms in
      $\Harm'(c)$.  The proof presented here is essentially the same
      as the proof in \cite{Takegoshi_higher-direct-images}*{Thm.~4.3
        (i)}.  For any $u \in \Harm'(c)$, set
      $u\ps := \eta_{c,\nu} \chi_{\nu'} u$ and
      $\dfadj_\Phi := \dfadj +\idxup{\diff\Phi} .\cdot$.  Since $u\ps$
      is smooth and compactly supported in $X_c^\circ$ and that
      $e^{-\Phi} \leq 1$, the Bochner--Kodaira--Nakano formula is
      valid for $u\ps$, which yields
      \begin{equation} \label{eq:BK-F+Phi}
        \begin{aligned}
          &~\int_{X_c^\circ} \abs{\dbar u\ps}_{\vphilist[\Phi]}^2
          +\int_{X_c^\circ} \abs{\dfadj_\Phi u\ps}_{\vphilist[\Phi]}^2
          \\
          =&~\int_{X_c^\circ} \abs{\nabla^{(0,1)}
            u\ps}_{\vphilist[\Phi]}^2 +\int_{X_c^\circ}
          \idxup{\ibddbar\paren{\vphi_F +\Phi}}
          \ptinner{u\ps}{u\ps}_{\vphilist[\Phi]}
        \end{aligned}
      \end{equation}
      (see, for example, \cite{Chan&Choi_injectivity-I}*{Lemma
        2.4.2}).  The fact $u \in \Harm'(c)$ implies that
      $\dbar u\ps = \dbar\paren{\eta_{c,\nu} \chi_{\nu'}} \wedge u$
      and
      $\dfadj_\Phi u\ps = -\idxup{\diff\paren{\eta_{c,\nu}
          \chi_{\nu'}}} . u$ (note that $\idxup{\diff\Phi}.u = 0$ is
      used here).  Also, we have the identities
      \begin{equation*}
        \nabla^{(0,1)} u\ps
        = \eta_{c,\nu} \chi_{\nu'} \nabla^{(0,1)} u
        +\dbar\paren{\eta_{c,\nu} \chi_{\nu'}} \otimes u
      \end{equation*}
      and
      \begin{equation*}
        \abs{\dbar \paren{\eta_{c,\nu} \chi_{\nu'}} \otimes u}_{\rs\omega}^2
        =\abs{\dbar \paren{\eta_{c,\nu} \chi_{\nu'}} \wedge u}_{\rs\omega}^2
        +\abs{\idxup{\diff \paren{\eta_{c,\nu} \chi_{\nu'}}} . u}_{\rs\omega}^2
        \quad\text{ on } X_c^\circ
      \end{equation*}
      (see \cite{Chan&Choi_injectivity-I}*{footnote 9 on p.33 (arXiv
        version)} or \cite{Federer}*{1.5.3}).  The
      \newtext{formula~\eqref{eq:BK-F+Phi}} is then reduced to
      \begin{align*}
        0 =
        &~\int_{X_c^\circ} \abs{\eta_{c,\nu} \chi_{\nu'} \nabla^{(0,1)} u}_{\vphilist[\Phi]}^2
          +\int_{X_c^\circ} \idxup{\ibddbar\paren{\vphi_F+\Phi}}
          \ptinner{u\ps}{u\ps}_{\vphilist[\Phi]}  \\
        &+2\Re\int_{X_c^\circ} \inner{
          \eta_{c,\nu} \chi_{\nu'} \nabla^{(0,1)}u
          }{
          \dbar\paren{\eta_{c,\nu} \chi_{\nu'}} \otimes u
          }_{\vphilist[\Phi]} \; .
      \end{align*}
      Note that
      $\dbar\paren{\eta_{c,\nu} \chi_{\nu'}} = \chi_{\nu'}
      \dbar\eta_{c,\nu} + \eta_{c,\nu} \dbar\chi_{\nu'}$ and
      $\dbar\chi_{\nu'} \to 0$ uniformly.  By taking the limit as
      $\nu' \to \infty$, the above formula yields
      \begin{equation*}
        \tag{$\dagger$} \label{eq:pf:BK-on-X_c-reduced}
        \begin{aligned}
          0 = &~\int_{X_c^\circ} \abs{\eta_{c,\nu} \nabla^{(0,1)}
            u}_{\vphilist[\Phi]}^2 +\int_{X_c^\circ}
          \idxup{\ibddbar\paren{\vphi_F+\Phi}}
          \ptinner{\eta_{c,\nu} u}{\eta_{c,\nu} u}_{\vphilist[\Phi]}  \\
          &+2\Re\int_{X_c^\circ} \underbrace{\inner{ \eta_{c,\nu}
              \nabla^{(0,1)}u }{ \dbar\eta_{c,\nu} \otimes u
            }_{\vphilist[\Phi]}}_{= \: \inner{ \eta_{c,\nu} \:
              \idxup{\diff\eta_{c,\nu}} \cdot \nabla_\bullet^{(0,1)}u
            }{\:u}_{\vphilist[\Phi]} \;\footnotemark} \; .
        \end{aligned}
      \end{equation*}%
      \footnotetext{\label{fn:explain_ctrt-nabla}%
        \newtext{For the sake of clarity, the notation
          ``$\idxup{\diff\eta_{c,\nu}} \cdot \nabla^{(0,1)}_\bullet
          u$'' means
          ``$\sum_j \paren{\diff\eta_{c,\nu}}^{\conj j} \nabla_{\conj
            j} u$'' in local coordinates.}  }%
      From the definition of $\eta_{c,\nu}$ in
      \eqref{eq:cut-off-functions}, we have
      $\dbar\eta_{c,\nu} = \frac{\eta_{c,\nu}' }{\nu (c - \Phi)^2}
      \dbar\Phi$, where
      \begin{equation*}
        \eta_{c,\nu}' := \rho'\paren{\frac{1}{\nu (c - \Phi)}} \leq 0
        \quad\text{ on } X_c \; .
      \end{equation*}
      The formula in \cite{Chan&Choi&Matsumura_injectivity}*{Lemma
        2.4.2} yields
      \begin{equation*}
        \dbar\cancelto{0}{\paren{\idxup{\diff\Phi}* . u}}
        =\idxup{\ibddbar\Phi}. u
        -\idxup{\diff\Phi} . \cancelto{0}{\paren{\dbar u}}
        +\idxup{\diff\Phi} \cdot \nabla^{(0,1)}_\bullet u \; .
      \end{equation*}
      Since $\Phi$ is psh, it follows that the term
      $2\Re\int_{X_c^\circ} \dotsm$ in \eqref{eq:pf:BK-on-X_c-reduced}
      is
      \begin{equation*}
        - 2 \int_{X_c^\circ} \eta_{c,\nu} \frac{\eta_{c,\nu}'}{\nu (c -
          \Phi)^2} \idxup{\ibddbar\Phi} \ptinner u u_{\vphilist[\Phi]}
        \geq 0 \; . 
      \end{equation*}
      Thus, we can see that each of the integrals in
      \eqref{eq:pf:BK-on-X_c-reduced} is $0$ (even without taking
      $\nu \to +\infty$) by noting the fact that $\vphi_F$ is psh.
      Considering the integrals involving $\nabla^{(0,1)}u$ and
      $\ibddbar\paren{\vphi_F +\Phi}$ and taking the limit
      $\nu \to +\infty$ yield $\nabla^{(0,1)} u = 0$,
      $\idxup{\ibddbar\vphi_F} \ptinner u u_{\rs\omega} = 0$ and
      $\idxup{\ibddbar\Phi} \ptinner u u_{\rs\omega} = 0$ on
      $X_c^\circ$.

    \item[(\ref{item:prop:su-properties})$_{\sect_D}$] 
    
      The \newtext{claim}
      % corresponding statement to Lemma \ref{lem:su-harmonicity}: if
      % $u \in \Harm'<X_c>{F;\Phi}$, then
      $\sect_D u \in \newtext{\paren{\ker
          \dbadj_{\phi_D}}^{n,q}_{(2),X_c} \cap~} \Harm'<X_c>{D\otimes
        F;\Phi},{\phi_D+\vphi_F}$ can be proved by the argument in
      \cite{Chan&Choi&Matsumura_injectivity}*{Lemma 2.2.1},
      % with the additional use of the cut-off function $\eta_{c,\nu}$
      % explained as above.  The proof is left to the reader.
      \begin{newpara}%
        as sketched below.  Write
        $\dfadj_{\phi_D} := \dfadj +\idxup{\diff\phi_D}. \cdot$ as the
        formal adjoint of $\dbar$ with respect to $\phi_D+\vphi_F$ and
        $\rs\omega$.  Since $\idxup{\diff\Phi} . \sect_D u = 0$,
        $\dbar\paren{\sect_D u} = \sect_D \dbar u = 0$ and
        $\dfadj_{\phi_D}\paren{\sect_D u} = \sect_D \dfadj u = 0$ on
        $X_c^\circ$ (see \cite{Chan&Choi&Matsumura_injectivity}*{Lemma
          2.2.1}), we see that
        $\sect_D u \in \Harm'<X_c>{D\otimes F;\Phi},{\phi_D+\vphi_F}$.
        To show that $\sect_D u$ lies in $\ker \dbadj_{\phi_D}$, put
        $\eta_{c,\nu} \chi_\nu \theta_\eps$ for $\nu \in \Nnum$ and
        $\eps >0$ in place of $\eta_{c,\nu} \chi_\nu$ (and $\sect_D u$
        in place of $u$) in the computation \eqref{eq:in-dom-dbadj}
        above.  More explicitly,
        \begin{align*}
          &~\iinner{\sect_D u}{\dbar\xi}_{X_c^\circ, \vphilist.}
          \\
          \xleftarrow{\nu \tendsto +\infty, \; \eps \tendsto 0^+}
          &
            \begin{aligned}[t]
              &-\cancelto{ 0 \;\;(\because \text{ Takegoshi property})
              }{ \iinner{ \chi_\nu \theta_\eps
                  \idxup{\diff\eta_{c,\nu}} . \sect_D u }{\xi}
              }_{X_c^\circ, \vphilist.}  -\iinner{ \eta_{c,\nu}
                \theta_\eps \idxup{\diff\chi_\nu}. \sect_D u
              }{\xi}_{X_c^\circ, \vphilist.}
              \\
              &-\eps \iinner{ \eta_{c,\nu} \chi_\nu \theta'_\eps
                \idxup{\frac{\diff\psi_D}{\abs{\psi_D}^{1+\eps}}}
                . \sect_D u }{\xi}_{X_c^\circ, \mathrlap{\vphilist.}}
              +\iinner{ \eta_{c,\nu} \chi_\nu \theta_\eps
                \dfadj_{\phi_D} \paren{\sect_D u} }{\xi}_{X_c^\circ,
                \vphilist.}
            \end{aligned}
          \\
          \xrightarrow{\eps \tendsto 0^+}
          &-\iinner{
            \eta_{c,\nu} \idxup{\diff\chi_\nu}. \sect_D u
            }{\xi}_{X_c^\circ, \vphilist.}            +\iinner{
            \eta_{c,\nu} \chi_\nu
            \sect_D \dfadj u
            }{\xi}_{X_c^\circ, \vphilist.}
          \\
          \xrightarrow{\nu \tendsto +\infty}
          &~\iinner{\sect_D \dfadj u}{\xi}_{X_c^\circ, \vphilist.}
            \qquad\text{(a bounded functional in $\xi$)}
        \end{align*}
        for any
        $\xi \in \Dom\dbar \subset \Ltwo/n,q-1/<X_c>{D\otimes
          F}_{\vphilist.}\;$.  Note that
        $\eps\eta_{c,\nu} \chi_\nu \theta'_\eps
        \idxup{\frac{\diff\psi_D}{\abs{\psi_D}^{1+\eps}}} . \sect_D u
        \xrightarrow{\eps \tendsto 0^+} 0$ in the norm
        $\norm\cdot_{X_c^\circ, \vphilist.}$ by the residue
        computation \cite{Chan&Choi_injectivity-I}*{Prop.~3.2.3 and
          Remark 3.2.4} (in which the smoothness of $u$ is needed).
        The claim is thus concluded.
      \end{newpara}

    \item[(\ref{item:prop:su-properties})$_{s}$] 
    \begin{newpara}%
      Write $\dfadj_{\vphi_M} := \dfadj +\idxup{\diff\vphi_M} . \cdot$
      (formal adjoint of $\dbar$ with respect to $\vphi_F+\vphi_M$ and
      $\rs\omega$).  If $\dfadj_{\vphi_M} (su) = 0$ (hence $L^2$ on
      $X_c^\circ$) is known, the remaining claim
      $su \in \paren{\ker \dbadj_{\vphi_M}}^{n,q}_{(2),X_c} \cap
      \Harm'<X_c>{F\otimes M; \Phi},{\vphi_F+\vphi_M}$ can be proved
      as in the claim for $\sect_D u$, using the cut-off function
      $\eta_{c,\nu} \chi_\nu$ in the computation
      \eqref{eq:in-dom-dbadj} above with $s u$ in place of $u$,
      $\vphi_F+\vphi_M$ in place of $\vphi_F$ and $\dfadj_{\vphi_M}$
      in place of $\dfadj$.  The fact $\dfadj_{\vphi_M} (su) = 0$
      follows as in \cite{Chan&Choi_injectivity-I}*{Cor.~3.2.6} (with
      $D=0$) by applying the Bochner--Kodaira--Nakano formula with
      respect to $\vphi_F+\vphi_M$ and $\rs\omega$ on the form $su$,
      which yields \eqref{eq:BK-F+Phi} with $\vphi_F+\vphi_M$ in place
      of $\vphi_F+\Phi$, $\dfadj_{\vphi_M}$ in place of
      $\dfadj_{\Phi}$ and $s u\ps := s \eta_{c,\nu} \chi_{\nu'} u$ in
      place of $u\ps$.  The rest of the computation is sketched here
      for the sake of clarity.  Following the computation as above and
      making use of the identities $\dbar\paren{s u} = s \dbar u = 0$,
      $\nabla^{(0,1)}\paren{s u} = s \nabla^{(0,1)} u = 0$ and
      $\idxup{\diff\eta_{c,\nu}} . su = 0$, we obtain
      \begin{alignat*}{1}
        &~\int_{X_c^\circ} \abs{ \eta_{c,\nu} \chi_{\nu'}
          \dfadj_{\vphi_M}\paren{su} }_{\vphilist M}^2 -2\Re
        \int_{X_c^\circ} \inner{ \eta_{c,\nu} \chi_{\nu'}
          \dfadj_{\vphi_M}\paren{su} }{\eta_{c,\nu}
          \idxup{\diff\chi_{\nu'}} . su}_{\vphilist M}
        \\
        = &~\int_{X_c^\circ} \idxup{\ibddbar\paren{\vphi_F+\vphi_M}}
        \ptinner{su\ps}{su\ps}_{\vphilist M} \; .
      \end{alignat*}
      The positivity assumption
      $\ibddbar\vphi_M \leq C \ibddbar\vphi_F$ for some constant
      $C > 0$, together with the Cauchy-Schwarz and AM-GM inequalities
      on the term $2\Re \int \inner\dots\dots$, yields, for some
      constant $\alpha \in (0,1)$,
      \begin{alignat*}{1}
        &~(1 - \alpha) \int_{X_c^\circ} \abs{ \eta_{c,\nu} \chi_{\nu'}
          \dfadj_{\vphi_M}\paren{su} }_{\vphilist M}^2
        \\
        \leq &~(1+C) \int_{X_c^\circ} \abs{\eta_{c,\nu} \chi_{\nu'}
          s}_{\vphi_M}^2 \cancelto{0}{\idxup{\ibddbar\vphi_F}
          \ptinner{u}{u}}_{\vphilist} +\frac 1\alpha \int_{X_c^\circ}
        \abs{ \eta_{c,\nu} \idxup{\diff\chi_{\nu'}} . su }_{\vphilist
          M}^2 \; .
      \end{alignat*}
      The identity $\dfadj_{\vphi_M} \paren{su} = 0$ thus follows from
      the estimate
      $\abs{\diff \chi_{\nu'}}_{\rs\omega} \lesssim \frac 1{\nu'}$ and
      by taking the limits $\nu' \tendsto +\infty$ and then
      $\nu \tendsto +\infty$.  (Note in particular that, while
      $\abs{\diff\eta_{c,\nu}}_{\rs\omega}$ is not uniformly bounded
      on $X_c^\circ$, the Takegoshi property
      $\idxup{\diff\eta_{c,\nu}} .s u = 0$ makes it possible to take
      the limit $\nu \tendsto +\infty$.)  This completes the proof.\qedhere
    \end{newpara}
  \end{enumerate}
  \end{proof}

  The Takegoshi property ensures that the homomorphism
  % $\jmath$ in \eqref{eq:map-for-harmonic-representatives}, as well
  % as 
  $\jmath^c$ in \eqref{eq:map-for-harmonic-representatives-on-X_c} is
  injective.
  Recall that $X_\infty = X$, and write %$\jmath^\infty := \jmath$ and
  $\eta_{\infty, \nu} := \eta_\nu$ for convenience.

  \begin{prop}[cf.~\cite{Takegoshi_higher-direct-images}*{Thm.~4.3 (iv)}]
    \label{prop:injective-jmath}
    Suppose that $\ibddbar\vphi_F \geq 0$ on $X$. 
    Then, the homomorphisms
    \begin{equation*}
      \jmath^c \colon \Harm'<X_c>{F;\Phi} \to
      \cohgp q[X_c]{\logKX \otimes \mtidlof{\vphi_F}}
      \quad\text{ for } c \in (0,\infty]
    \end{equation*}
    defined as in %\eqref{eq:map-for-harmonic-representatives} and
    \eqref{eq:map-for-harmonic-representatives-on-X_c} are all injective. 
  \end{prop}

  \begin{proof}
    % The case where $c < \infty$ is proved in
    % \cite{Takegoshi_higher-direct-images}*{Thm.~4.3 (iv)}.
    % The proof below is valid for all $c > 0$, including the case $c
    % = \infty$.

    Let $u \in \Harm'(c)$ be a harmonic form such that $\jmath^c(u) = 0$.
    \newtext{Note that the $L^2_\tloc$ Dolbeault
      isomorphism \eqref{eq:L2-Dolbeault-isom} holds true also on $X_c$ for all $c \in (0, \infty]$.
    Therefore, we have} $u = \dbar\xi$ for some $\xi \in
    \Ltwo./n,q/<X_c>{F}_{\vphilist}$ on $X_c$.
   
    While $\xi$ may not be globally $L^2$ on $X_c^\circ$, the form $\eta_{c,\nu} \xi$ is globally $L^2$.
    The Takegoshi property $\idxup{\diff\Phi} . u = 0$ ensures that
    $\idxup{\diff\eta_{c,\nu}} . u = 0$ on $X_c^\circ$ 
    by the definitions of $\eta_\nu = \eta_{\infty,\nu}$ and $\eta_{c,\nu}$ (see \eqref{eq:cut-off-functions}). 
    The orthogonality between $\Harm'(c)$ and $\cl{\im\dbar}_{(2)}$ due to Proposition~\ref{prop:Takegoshi-harmonic-forms}
    then yields
    \begin{align*}
      0
      = \iinner{u}{\dbar\paren{\eta_{c,\nu} \xi}}_{\vphilist}
      &= \iinner{u}{\:\dbar\eta_{c,\nu} \wedge \xi + \eta_{c,\nu} \dbar\xi}_{\vphilist}
      \\
      &= \cancelto{0 \;\;(\because \text{Takegoshi property})}{
        \iinner{
        \idxup{\diff\eta_{c,\nu}} . u
        }{\:\xi}}_{\vphilist}
        + \iinner{u}{\eta_{c,\nu} u}_{\vphilist}
      \\
      &\xrightarrow{\nu \to +\infty} \norm u_{\vphilist}^2 \; .
    \end{align*}
    Therefore, we have $u = 0$ and $\jmath^c$ is injective.
  \end{proof}

  The vanishing of $\nabla^{(0,1)} u$ in Proposition
  \ref{prop:consequence-of-positivity} or
  \ref{prop:Takegoshi-harmonic-forms} helps to control the
  singularities of $u$ along $X \setminus X^\circ = P_F \cup P_M$. 
  Recall the function $\psi_{P_F\cup P_M}$ given in Section
  \ref{subsec:notation}.
  Also, on any admissible open set $V \subset X$, write
  \newtext{$D \cap V = \set{z_1 \dotsm z_{\sigma_V} = 0}$ and}
  $\paren{P_F\cup P_M} \cap V = \set{w_1\dotsm w_\mu = 0}$, where
  $z_1, \dots, z_{\sigma_V}, w_1,\dots,w_\mu$ are part of the
  holomorphic coordinates on $V$.

  \begin{prop}[cf.~\cite{Chan&Choi_injectivity-I}*{Thm.~2.5.1 and
      Prop.~3.3.1 and 3.3.2}]
    \label{prop:singularities-along-FM}
    
    Suppose that $\ibddbar\vphi_F \geq 0$ and $u \in
    \Harm'<X_c>{F;\Phi}$ for some given $c \in (0,\infty]$.
    Then, the form $*u$ is holomorphic on $X_c$, where $*$ is the
    Hodge $*$-operator with respect to $\rs\omega$.
    In particular, on any admissible open set $V \Subset X_c$, we have 
    \begin{equation*}
      \text{coef.~of } u \text{ and } \idxup{\diff\psi_D} . u \sect_D 
      \in \smooth_X \left[
        \newtext{\frac 1{\det\rs\omega} ,
        \frac 1{\abs{\psi_{P_F\cup P_M}}} ,
        \:\frac 1{\log\abs{\ell \psi_{P_F\cup P_M}}}} ,
        \: \frac 1{\abs{w_1}} , \dots , \frac 1{\abs{w_\mu}}
      \right] \quad\text{ on } V \; .
    \end{equation*}
    Moreover,  the restriction map $\jmath^{c'}_c$
    for $c' > c$ in
    \eqref{eq:restriction-maps-between-Takegoshi-harm-sp} is injective.
  \end{prop}

  \newtext*{
    \begin{remark} \label{rem:inverse-of-omega}
      It is shown in
      \cite{Chan&Choi_injectivity-I}*{Proof of Prop.~3.3.1} by a direct
      computation that
      coefficients of $\rs\omega$ lie in the algebra
      \begin{equation*}
        \smooth_X \left[
          \frac 1{\abs{\psi_{P_F\cup P_M}}} ,
          \:\frac 1{\log\abs{\ell \psi_{P_F\cup P_M}}} ,
          \: \frac 1{\abs{w_1}} , \dots , \frac 1{\abs{w_\mu}}
        \right] \; .
      \end{equation*}
      From the computation it can also be seen that
      \begin{equation*}
        \frac 1{\det \rs\omega}
        = \abs{w_1 \dotsm w_\mu}^{2}
        \abs{\psi_{P_F\cup P_M}}^{2\mu}
        \paren{\log\abs{\ell \psi_{P_F\cup P_M}}}^{2\mu} \:B
        \quad\text{ on } V
      \end{equation*}
      for some positive (nowhere zero) continuous function $B$ on $V$
      which is infinitely differentiable with respect to the variables
      $r_1, \dots, r_{\sigma_V}$, where $r_j := \abs{z_j}$ for $j = 1,
      \dots, \sigma_V$.
    \end{remark}
  }
  
  \begin{proof}
    A refined statement of the hard Lefschetz theorem
    \cite{Matsumura_injectivity-lc}*{Thm.~3.3} (see also
    \cite{Chan&Choi_injectivity-I}*{Thm.~2.5.1}) guarantees, in the
    case where $X$ is compact, the holomorphicity of $*u$ (more
    precisely, $*u \sect_D$). 
    The proof of that statement was complicated by the fact that one
    started by considering $\logKX<X>.$-valued harmonic forms with
    respect to $\phi_D +\vphi_F$ and $\rs\omega$.
    In the current situation (where $u$ is only $\logKX<X>$-valued), a
    simple proof, which works also on the 
    non-compact $X_c$, can be given as follows.
    Note that $\nabla^{(0,1)} u = 0$ on $X_c^\circ$ according to
    Proposition \ref{prop:consequence-of-positivity} or \ref{prop:Takegoshi-harmonic-forms}.
    From \cite{Chan&Choi_injectivity-I}*{Remark 2.4.3}, we have
    \begin{equation*}
      0
      = \norm{\nabla^{(0,1)} u}_{\vphilist}^2
      = \norm{* \dbar * u}_{\vphilist}^2
      = \norm{\dbar * u}_{\vphilist}^2 \; .
    \end{equation*}
    It follows that $*u$ is holomorphic on $X_c^\circ$.
    Furthermore, we have $\norm{* u}_{\vphilist,{\alert{\omega}}}^2 \leq
    \norm{* u}_{\vphilist}^2$, which is a consequence of $\rs\omega \geq \omega$.
    This inequality implies that $*u$ is locally $L^2$ everywhere in
    $X_c$ (not only $X_c^\circ$) with respect to the unweighted $L^2$ norm.
    Consequently, by \cite{Demailly_complete-Kahler}*{Lemma 6.9},
    $*u$ is holomorphic on the entire $X_c$.

    The smoothness of $* u$ on $X_c$, together with the pointwise identity
    \begin{equation*}
      u \wedge \conj{*u} = \abs{u}_{\rs\omega}^2
      = \abs{* u}_{\rs\omega}^2 \:\frac{\rs\omega^{\wedge n}}{n!} \; ,
    \end{equation*}
    shows that the singularities of $u$ along $X_c \setminus
    X_c^\circ$ are no worse than those of \newtext{the product of} $\rs\omega^{\wedge n}$
    \newtext{and the coefficients of the inverse of $\rs\omega$
      resulted from contracting the $(n-q)$-forms $* u$ and $\conj{* u}$,
      % and
      % thus of the form given in the claim (indeed $\abs{\psi_{P_F\cup P_M}}$
      % and $\log\abs{\ell\psi_{P_F\cup P_M}}$ are not needed in the generators of
      % the algebra in this case),
      hence the claim for the coefficients of $u$ as in
      \cite{Chan&Choi_injectivity-I}*{Prop.~3.3.1}\footnote{%
        \newtext{The algebra in \cite{Chan&Choi_injectivity-I}*{Prop.~3.3.1}
        is missing the generator $\frac 1{\det\rs\omega}$ by mistake.
        Fortunately it does not affect the arguments in the proof of
        \cite{Chan&Choi_injectivity-I}*{Prop.~3.3.2} and other parts
        of the paper.}
      }}.
    The singularities of $\idxup{\diff\psi_D} . u \sect_D$ along $X_c \setminus
    X_c^\circ$ are also of the form given in the claim
    % (where $\abs{\psi_{P_F\cup P_M}}$ and $\log\abs{\ell\psi_{P_F\cup P_M}}$ are
    % needed in the generators of the algebra in this case)
    as explained in
    the proof of \cite{Chan&Choi_injectivity-I}*{Prop.~3.3.2}.

    To see that $\jmath^{c'}_c$ is injective, suppose $u \in
    \Harm'(c')$ and $j^{c'}_c (u) = \res u_{X_c} \equiv
    0$.
    However, $*u$ is holomorphic on $X_{c'}$ and $X_c$ is open in
    $X_{c'}$, so $\res{*u}_{X_c} \equiv 0$ implies that $u =
    (-1)^{n-q} ** u \equiv 0$ on $X_{c'}$ by the identity theorem, as desired.
  \end{proof}
  
  The Takegoshi property with $\Phi$ also enables us to apply the
  twisted Bochner--Kodaira formula (see, for example,
  \cite{Chan&Choi_injectivity-I}*{Lemma 2.4.2}) to harmonic forms in
  $\Harm'(c)$ for all $c \in (0,\infty]$, which results in the
  following proposition. 

  \begin{prop}[cf.~\cite{Chan&Choi_injectivity-I}*{
      Prop.~3.2.3, 3.2.8 and 3.3.2}]
    \label{prop:Takegoshi-property-on-tBK}
    Suppose that $\ibddbar\vphi_F \geq 0$ and $u \in
    \Harm'<X_c>{F;\Phi}$.
    Then, we have
    \begin{equation*}
      \frac{
        \abs{\idxup{\diff\psi_D} . u}_{\vphilist}^2
      }{
        \abs{\psi_D}^{1+\eps}
      } \in \Lloc(\alert{X_c}) \quad(\text{not only on $X_c^\circ$!\:})
      \quad\text{ for any } \eps > 0 \; .
    \end{equation*}
    Furthermore, given the cut-off functions $\eta_{c,\nu}$ for $\nu \in\Nnum$ in
    \eqref{eq:cut-off-functions}, we have
    \begin{align*}
      \int_{X_c^\circ} \idxup{\ibddbar\sm\vphi_D}
      \ptinner{u}{u}_{\vphilist}
      &=\lim_{\nu \tendsto +\infty} \lim_{\eps \tendsto 0^+} \eps \int_{\mathrlap{X_c^\circ}} \;\;
        \frac{
        \abs{\idxup{\diff\psi_D} . \eta_{c,\nu} u}_{\vphilist}^2
        }{
        \abs{\psi_D}^{1+\eps}
        }
      \\
      &=\pi \sum_{b \in \Iset|1|} \int_{\lcS*|1|<c>[b]}
        \abs{\PRes[\lcS|1|[b]](\idxup{\diff\psi_D}. u)}_{\vphilist}^2
        \; ,
    \end{align*}
    where $D = \sum_{b \in \Iset|1|} \lcS|1|[b]$ with each
    $\lcS|1|[b]$ being an irreducible divisor, $\PRes[\lcS|1|[b]]$ is the
    Poincar\'e residue map from $X$ to $\lcS|1|[b]$, and $\lcS*|1|<c>[b] :=
    \lcS|1|[b] \cap X_c^\circ$.
    In particular, 
    $\PRes[\lcS|1|[b]](\idxup{\diff\psi_D}. u)$ is globally $L^2$ with
    respect to $\res{\vphi_F}_{\lcS|1|[b]}$ and
    $\res{\rs\omega}_{\lcS|1|[b]}$ on $\lcS*|1|<c>[b]$.
  \end{prop}

  \begin{proof}
    \setDefaultlcIndex{1}

    The assumptions on $\vphi_F$ and $u$ imply that
    \begin{enumerate}[label=(\emph{\roman*}), ref=\emph{\roman*}]
    \item \label{item:smooth-harmonic-u}
      $u$ is $\dbar$- and $\dfadj$-closed and is smooth in
      $X_c^\circ$ as a harmonic form (by the regularity of the
      $\dbar$-operator),
      
    \item \label{item:conseq-of-pos}
      $\nabla^{(0,1)} u = 0$ and $\idxup{\ibddbar\vphi_F} \ptinner
      u u_{\vphilist} = 0$ on $X_c^\circ$ by Proposition
      \ref{prop:consequence-of-positivity} or Proposition \ref{prop:Takegoshi-harmonic-forms}, and
      
    \item \label{item:Takegoshi-property}
      $\idxup{\diff\eta_{c,\nu}} . u = 0$ on $X_c^\circ$ (by the Takegoshi
      property with $\Phi$).
    \end{enumerate}
    Proposition \ref{prop:singularities-along-FM} provides crucial information about 
    the singularities of $u$ along $X_c \setminus X_c^\circ$. 
    With this knowledge, we can establish our claims by following the proofs 
    presented in \cite{Chan&Choi_injectivity-I}. 
    % \xb{Theorem numbers are not associated.}
    Specifically, we refer to \cite{Chan&Choi_injectivity-I}*{Prop.~3.2.3 and 3.3.2}, along with the 
    modification considered in \cite{Chan&Choi_injectivity-I}*{Remark
      3.2.4}.
    % {I think citing like this would be less confusing.
    % \xb{Than you!}

    We now present a brief outline of the key arguments.
    % Putting $u\ps_\nu := \eta_{c,\nu} u$ for $\nu \in \Nnum$.
    % , the
    % first equality in the claim is proved in
    % \cite{Chan&Choi_injectivity-I}*{Prop.~3.2.3 and 3.3.2 with Remark
    % 3.2.4} ($\eta_{c,\nu} u$ here corresponds to $\frac u{\sect_D}$ or
    % $\frac{\rs*u}{\sect_D}$ there)
    % via the residue computation.
    On any admissible set $V_i$ in the finite open cover $\cvr V =
    \set{V_i}_{i \in I}$ given in Section \ref{subsec:notation} on
    which
    \begin{equation*}
      D = \set{z_1 \dotsm z_{\sigma_{V_i}} = 0} 
      \quad\text{ and } \quad
      \lcS|1|[b] \cap V_i = \set{z_{b(1)} = 0}
      \;\;\forall~
      b \in \Iset|1|
      \text{ such that }
      \lcS|1|[b] \cap V_i \neq \emptyset 
    \end{equation*}
    (note also that $\sect\ps_b = z_{b(2)} \dotsm
    z_{b(\sigma_{V_i})}$), we have
    \begin{equation*}
      \idxup{\diff\psi_D} . \eta_{c,\nu} u
      = \sum_{b \in \Iset|1| \colon \lcS|1|[b] \cap V_i \neq
        \emptyset}
      \underbrace{
        \idxup{\frac{dz_{b(1)}}{z_{b(1)}}} . \eta_{c,\nu} u
      }_{ =:\: \frac{dz_{b(1)}}{z_{b(1)}} \wedge g_{b}} 
      \; - \idxup{\diff\sm\vphi_D} . \eta_{c,\nu} u  \quad\text{ on } V_i \; . 
    \end{equation*}
    Here $g_b$ is an $(n-1, q-1)$-form free from the forms
    $dz_{b(1)}$ and $d\conj{z_{b(1)}}$ such that
    \begin{equation*}
      \res{g_b}_{\lcS[b] \cap V_i}
      = \res{\PRes[\lcS[b]](\idxup{\frac{dz_{b(1)}}{z_{b(1)}}} . \eta_{c,\nu} u)}_{\lcS[b] \cap V_i}
      = \res{\PRes[\lcS[b]](\idxup{\diff\psi_D}
        . \eta_{c,\nu} u)}_{\lcS[b] \cap V_i} \; ,
    \end{equation*}
    which is a compactly supported smooth
    $\logKX<\lcS[b]>[\res{F}_{\lcS[b]}]$-valued $(0,q-1)$-form on
    $\lcS[b]$.
    The form of the singularities of $u$ along $X_c \setminus
    X_c^\circ$ shown in Proposition \ref{prop:singularities-along-FM}
    implies that the singularities there and the lc locus $D$ are
    ``independent'' of each other in view of Fubini's theorem
    \newtext{\label{page:sing-u-and-D-indep_explain}(i.e.~the variables defining the singularities of $u$ are
      independent from those defining $D$ when integration is concerned)}.
    A direct residue computation as in
    \cite{Chan&Choi_injectivity-I}*{Thm.~2.6.1} then ensures that $\frac{
      \abs{\idxup{\diff\psi_D} . u}_{\vphilist}^2
    }{
      \abs{\psi_D}^{1+\eps}
    } \in \Lloc(\alert{X_c})$ for all $\eps > 0$ and gives the
    equality 
    \begin{equation*}
      \lim_{\eps \tendsto 0^+} \eps \int_{\mathrlap{X_c^\circ}} \;\;
      \frac{
        \abs{\idxup{\diff\psi_D} . \eta_{c,\nu} u}_{\vphilist}^2
      }{
        \abs{\psi_D}^{1+\eps}
      }
      =\pi \sum_{b \in \Iset|1|} \int_{\lcS*|1|<c>[b]}
      \abs{\PRes[\lcS|1|[b]](\idxup{\diff\psi_D}. \eta_{c,\nu} u)}_{\vphilist}^2
      \; .
    \end{equation*}
    % (This computation actually guarantees that $\frac{
    % \abs{\idxup{\diff\psi_D} . u}_{\vphilist}^2
    % }{
    %   \abs{\psi_D}^{1+\eps}
    % } \in \Lloc[1](\alert{X_c^\circ})$ for all $\eps > 0$.)

    To show that these integrals converge when $\nu \ascendsto +\infty$,
    we can follow the proof in 
    \cite{Chan&Choi_injectivity-I}*{Prop.~3.2.8 with Lemma 3.3.3}
    by using the twisted Bochner--Kodaira formula in
    \cite{Chan&Choi_injectivity-I}*{Lemma 2.4.2}.

    \setDefaultPscript{\nu,\nu',\eps'}
    
    Write $u\ps := \eta_{c,\nu} \chi_{\nu'} \theta_{\eps'} u$ for
    $\nu, \nu' \in\Nnum$ and $\eps' > 0$
    (see \eqref{eq:cut-off-functions} for the cut-off functions
    $\eta_{c,\nu}$, $\chi_{\nu'}$ and $\theta_{\eps'}$).
    Then, the twisted Bochner--Kodaira formula in
    \cite{Chan&Choi_injectivity-I}*{Lemma 2.4.2} is valid for $u\ps$
    (by \eqref{item:smooth-harmonic-u}) and is read as
    \begin{align*}
      &~\alert[Gray]{\int_{X_c^\circ} \abs{\dbar u\ps}_{\vphilist}^2 \:\abs{\psi_D}^{1-\eps}
        +\int_{X_c^\circ} \abs{\dfadj u\ps}_{\vphilist}^2
        \:\abs{\psi_D}^{1-\eps}
        -\int_{X_c^\circ} \abs{\nabla^{(0,1)}  u\ps}_{\vphilist}^2
        \:\abs{\psi_D}^{1-\eps}}
      \\ % \tag*{$\tBK_{\eps,\vphilist}$}
      =&\!
         \begin{aligned}[t]
           &~\int_{X_c^\circ}
           \idxup{\cancelto{0 \;\text{ by \eqref{item:conseq-of-pos}}}{
               \alert[Gray]{\ibddbar\vphi_F}\vphantom{\frac{A}{1}}
             }
             +\frac{1-\eps
             }{\abs{\psi_D}} \ibddbar\psi_D }
           \ptinner{u\ps}{u\ps}_{\vphilist} \:\abs{\psi_D}^{1-\eps}
           \\
           &+
           \eps(1-\eps)
           \int_{X_c^\circ} 
           \abs{\idxup{\diff\psi_D} .  u\ps}_{\vphilist}^2
           \frac{\abs{\psi_D}^{1-\eps}}{\abs{\psi_D}^2}
           \\
           &~\alert[Gray]{\newtext{-}\, 2\paren{1-\eps} \Re
             \int_{X_c^\circ}
             \inner{\dfadj u\ps}{
               \frac{
                 \idxup{\diff\psi_D} .  u\ps
               }{
                 \abs{\psi_D}
               }
             }_{\mathrlap{\vphilist}} \;\;\abs{\psi_D}^{1-\eps}}
           \; .
         \end{aligned}
    \end{align*}
    From the equation $\dbar u\ps = \dbar\paren{\eta_{c,\nu}
      \chi_{\nu'} \theta_{\eps'}} \wedge u$, together with the analogous
    equations for $\dfadj u\ps$ and $\nabla^{(0,1)}
    u\ps$, and also the identity 
    \begin{equation*}
      \abs{\dbar \paren{\eta_{c,\nu} \chi_{\nu'} \theta_{\eps'}} \otimes u}_{\vphilist}^2
      =\abs{\dbar \paren{\eta_{c,\nu} \chi_{\nu'} \theta_{\eps'}} \wedge u}_{\vphilist}^2
      +\abs{\idxup{\diff \paren{\eta_{c,\nu} \chi_{\nu'} \theta_{\eps'}}} . u}_{\vphilist}^2
      \quad\text{ on } X_c^\circ
    \end{equation*}
    (see \cite{Chan&Choi_injectivity-I}*{footnote 9 on p.33 (arXiv
      version)} or \cite{Federer}*{1.5.3}), 
      the twisted Bochner--Kodaira formula is reduced to
    \begin{equation*} \tag{$\dagger$} \label{eq:pf:tBK-after-cancellation}
      \begin{aligned}
        0=
        &~
        -(1-\eps) \int_{\mathrlap{X_c^\circ}} \;\;
        \frac{
          \idxup{ \ibddbar\sm\vphi_D }
        }{\abs{\psi_D}^\eps}
        \ptinner{u\ps}{u\ps}_{\vphilist}
        \newtext{\qquad\paren{\because~\ibddbar\phi_D = 0 \;\text{ on
            } \supp \theta_{\eps'}}}
        \\
        &~+
        \eps(1-\eps)
        \int_{\mathrlap{X_c^\circ}} \;\;
        \frac{
          \abs{\idxup{\diff\psi_D} . u\ps}_{\vphilist}^2
        }{\abs{\psi_D}^{1+\eps}} 
        \\
        &~\alert[Gray]{~\newtext{+}\, 2\paren{1-\eps} \Re
          \int_{X_c^\circ}
          \inner{\idxup{\diff\paren{\eta_{c,\nu} \chi_{\nu'} \theta_{\eps'}}} . u}{
            \frac{
              \idxup{\diff\psi_D} .  u\ps
            }{
              \abs{\psi_D}^\eps
            }
          }_{\vphilist}} \; .
      \end{aligned}
    \end{equation*}
    % \newtext{Here we have $\ibddbar\psi_D= -\ibddbar\sm\vphi_D$, since 
    % $\ibddbar\phi_D=0$ on $X_c^\circ$.}
    The limits are taken in the order: $\nu' \ascendsto +\infty$,
    $\eps' \descendsto 0$, $\eps \descendsto 0$ and then $\nu
    \ascendsto +\infty$.
    First check that the last term (in \alert[Gray]{Gray}) converges to $0$
    in the limit.
    Note that
    \begin{equation*}
      \idxup{\diff\paren{\eta_{c,\nu} \chi_{\nu'} \theta_{\eps'}}} . u
      \overset{
        % \substack{\text{Takegoshi} \\ \text{property}}
        \text{\eqref{item:Takegoshi-property}}
      }=
      \eta_{c,\nu} \idxup{\diff\paren{\chi_{\nu'} \theta_{\eps'}}} . u
      =\alert[Cyan]{\eta_{c,\nu}  \theta_{\eps'}\idxup{\diff \chi_{\nu'}} . u}
      +\alert[RoyalBlue]{\eta_{c,\nu} \chi_{\nu'} \: \frac{\eps' \theta'_{\eps'}
          \idxup{\diff\psi_D} . u}{\abs{\psi_D}^{1+\eps'}}} 
    \end{equation*}
    (recall that $\theta'_{\eps'} := -\rho'\circ
    \frac{1}{\abs{\psi_D}^{\eps'}}$), and the term
    \begin{equation*}
      \int_{X_c^\circ} \inner{
        \alert[Cyan]{\eta_{c,\nu}  \theta_{\eps'}\idxup{\diff \chi_{\nu'}} . u}
      }{
        \frac{\idxup{\diff\psi_D} . u\ps}{\abs{\psi_D}^\eps}
      }_{\vphilist}
      \xrightarrow{\nu' \tendsto +\infty} 0
      \quad\text{ for any } \nu \in \Nnum \text{ and } \eps, \eps' > 0 \; ,
    \end{equation*}
    as $\abs{\diff\chi_{\nu'}}_{\rs\omega} \tendsto 0$ uniformly on $X_c^\circ$
    and $\eta_{c,\nu} \theta_{\eps'} \frac{\idxup{\diff\psi_D} . u}{\abs{\psi_D}^\eps}$ is
    $L^2$ in $X_c\setminus D$ with respect to $\vphi_F$ and $\rs\omega$.
    Furthermore,
    \begin{equation*}
      \begin{aligned}
        \int_{X_c^\circ} \inner{
          \alert[RoyalBlue]{
            \eta_{c,\nu} \chi_{\nu'} \frac{
              \eps' \theta'_{\eps'} \idxup{\diff\psi_D} . u
            }{\abs{\psi_D}^{1+\eps'}}
          }
        }{
          \frac{\idxup{\diff\psi_D} . u\ps}{\abs{\psi_D}^\eps}
        }_{\vphilist}
        &=\eps' \int_{\mathrlap{X_c^\circ}} \;\;
        \eta_{c,\nu}^2 \chi_{\nu'}^2 \theta_{\eps'} \theta'_{\eps'}
        \frac{
          \abs{\idxup{\diff\psi_D}. u}_{\vphilist}^2
        }{
          \abs{\psi_D}^{1+\eps +\eps'}
        }
        \\
        &\xrightarrow{\substack{\nu'
            \tendsto +\infty \\ \eps' \tendsto 0^+}} 0
      \end{aligned}
    \end{equation*}
    for any $\eps > 0$ and $\nu \in \Nnum$.
    Note that the fact that
    $\frac{\abs{\idxup{\diff\psi_D}.u}_{\vphilist}^2}{\abs{\psi_D}^{1+\eps}}
    \in \Lloc[1](X_c)$ for any $\eps > 0$ is used to ensure that the
    limit exists as $\nu' \ascendsto +\infty$.

    As a result, the term in \alert[Gray]{Gray} in
    \eqref{eq:pf:tBK-after-cancellation} goes to $0$ after the
    limits $\nu' \ascendsto +\infty$ and $\eps' \descendsto 0$.
    Notice that $\ibddbar\sm\vphi_D$ is bounded on $X$ and $u$ is
    $L^2$ with respect to $\vphi_F$ and $\rs\omega$ on $X_c^\circ$.
    Further taking the limits $\eps \descendsto 0$ and $\nu \ascendsto
    +\infty$ to \eqref{eq:pf:tBK-after-cancellation} yields the
    desired result.
  \end{proof}

} %%%%%% End of the redefinition of DefaultMetric %%%%%%

%%% Local Variables:
%%% mode: latex
%%% TeX-master: "Relative-Fujino"
%%% coding: utf-8
%%% End:

\section{Residue formula and harmonic residue}

% \xb{In the following subsections, we develop the residue computation and systematically study harmonic residues, which play a crucial role in carrying out the inductive argument using lc centers.}

% \mmark{}{Need some descriptive text for this section here.}

\mmark{
  The residue formula with respect to the lc centers of $(X,D)$ and the
  corresponding harmonic residue are established in
  \cite{Chan&Choi&Matsumura_injectivity} (although we coin the name
  ``harmonic residue'' only here).
  In this section, the residue statement is recalled and adapted to the
  current setup (with non-compact $X$ and singular $\vphi_F$ and $\vphi_M$). 
  The treatment to the singularities of $\vphi_F$ and $\vphi_M$ follows
  the one in \cite{Chan&Choi_injectivity-I}.
  While most statements and techniques used in the proofs in this
  section come from our previous works, the adjoint relation between the
  harmonic residue and the connecting morphism for the cohomology groups
  (Theorem \ref{thm:HRes-duality}) is new to us.
}{Sorry again for the change... Here I try to emphasise what
  are known and what is new in this section.}

\subsection{Adjoint ideal sheaves and the residue computation}
\label{sec:adjoint-ideal-n-residue}

{
  \setDefaultvphi{\vphi_L}

  We review below the basics of adjoint ideal sheaves in
  \cite{Chan_adjoint-ideal-nas}.
  Recall that $(L,\vphi_L)$ is a line bundle on $X$ equipped with a
  potential $\vphi_L$ with only \emph{neat analytic singularities}. 
  The adjoint ideal sheaf $\aidlof :=\aidlof<X>$ of index $\sigma$
  is given at each $x \in X$  by
  \begin{equation*}
    \aidlof,_x :=\setd{f \in \holo_{X,x}}{
      \begin{multlined}[c][0.4\textwidth]
        \exists~\text{open set } V_x \ni x \:, \: \forall~\eps > 0 \:,
        \\
        \eps\int_{V_x} \frac{\abs{f}^2 e^{-\vphi_L -\psi_D} \dvol_{V_x}}{\logpole} < +\infty
      \end{multlined}
    } \; .
  \end{equation*}
  Assume that $\lcdata<X>$ is in the snc configuration and that
  $\vphi_L^{-1}(-\infty)$ contains no lc centers of $(X,D)$
  (both of which hold true when $\vphi_L := \vphi_F$ or $\vphi_F +\vphi_M$
  according to the setup given in Section \ref{subsec:notation}).
  \newtext{\label{page:explain_m0m1}This implies that $\mtidlof{\vphi_L+\psi_D} \subsetneq
    \mtidlof{\vphi_L+m\psi_D} = \mtidlof{\vphi_L}$ for all $m \in
    [0,1)$ (so we can take $m_0=0$ and $m_1 =1$ in the notation in
    \cite{Chan_adjoint-ideal-nas}).}
  Then, \cite{Chan_adjoint-ideal-nas}*{Thm.~1.2.3} shows that the adjoint ideal
  sheaf can be written as 
  \begin{equation*}
    \aidlof = \mtidlof{\vphi_L} \cdot \defidlof{\lcc+1'}
    \quad\text{ for any } \sigma \geq 0 \; ,
  \end{equation*}
  where $\defidlof{\lcc+1'}$ is the defining ideal sheaf of $\lcc+1'$
  in $X$ (with the reduced structure).
  Furthermore, we have the residue short
  exact sequence
  \begin{equation} \label{eq:residue-exact-seq}
    \xymatrix@R-0.5cm@C+0.3cm{
      {0} \ar[r]
      & {\aidlof-1} \ar[r]
      & {\aidlof} \ar[r]^-{\Res^\sigma}
      & {\residlof} \ar[r]
      & {0 \; .}
    }
  \end{equation}
  Here the quotient sheaf ${\residlof}$, called the \emph{residue
    sheaf of index $\sigma$}, is supported on $\lcc' = \sum_{p \in
    \Iset} \lcS$ and given by 
  \begin{align*}
    \residlof
    &= \bigoplus_{p \in \Iset} \paren{\Diff_p D}^{-1}
      \otimes \mtidlof<\lcS>{\vphi_L} \;\;\text{ and thus }
    \\
    \logKX.[L] \otimes \residlof
    &=\bigoplus_{p \in\Iset} K_{\lcS}
      \otimes \res L_{\lcS} \otimes \mtidlof<\lcS>{\vphi_L} \; .
  \end{align*}

  The \emph{residue morphism $\Res^\sigma$} can be given in terms of
  the Poincar\'e residue map $\PRes[\lcS]$ given in
  \cite{Kollar_Sing-of-MMP}*{\S 4.18}.
  The Poincar\'e residue map $\PRes[\lcS]$ from $X$ to each $\lcS$ is
  uniquely determined after an orientation on the holomorphic conormal
  bundle of $\lcS$ in $X$ is fixed.
  On an admissible open set $V \subset X$ such that $\lcS
  \cap V =\set{z_{p(1)} =z_{p(2)} =\dotsm =z_{p(\sigma)}=0}$, 
  a section $f $ of  $\logKX.[L] \otimes
  \aidlof$ on $V \subset X$ can be written as
  \begin{equation*}
    f = \;\;\smashoperator{\sum_{p \in \Iset \colon \lcS \cap V
        \neq\emptyset}} \;\; dz_{p(1)} \wedge \dotsm \wedge dz_{p(\sigma)}
    \wedge g_p \:\sect_{(p)} 
    =\;\;\smashoperator[l]{\sum_{p \in \Iset \colon \lcS \cap V
        \neq\emptyset}}
    \frac{dz_{p(1)}}{z_{p(1)}} \wedge \dotsm
    \wedge \frac{dz_{p(\sigma)}}{z_{p(\sigma)}}
    \wedge g_p \:\sect_D \quad\text{ on } V. 
  \end{equation*}
  % Then, the Poincar\'e residue map $\PRes[\lcS]$ is given by
  Assuming that the orientation on the conormal bundle of
  $\lcS$ in $X$ on $V$ is given by $(dz_{p(1)}, \dots, dz_{p(\sigma)})$,
  we see that 
  \begin{equation*}
    \PRes[\lcS](\frac{f}{\sect_D})  =\res{g_p}_{\lcS} \in
    K_{\lcS} \otimes \res L_{\lcS} \otimes \mtidlof<\lcS>{\vphi_L}
    \quad\text{ on } \lcS \cap V \; .
  \end{equation*}
  Note that $\res{g_p}_{\lcS}$ takes values in
  $\mtidlof<\lcS>{\vphi_L}$ according to 
  \cite{Chan&Choi_ext-with-lcv-codim-1}*{Prop.~2.2.1}, which says
  that
  \begin{equation*} % \label{eq:residue-norm}
    \begin{aligned}
      \norm{g}_{\lcc' \cap V}^2
      :=\sum_{p\in\Iset} \norm{g_p}_{\lcS \cap V}^2
      := &~\sum_{p \in \Iset} \frac{\pi^\sigma}{(\sigma -1)!}
      \int_{\mathrlap{\lcS \cap V}} \quad \abs{g_p}^2 \:e^{-\vphi_L}
      \\
      % =&~\lim_{\rho \descendsto \charfct_{\cl V}}
      % \lim_{\eps \tendsto 0^+}
      % \eps \int_{V'} \frac{\rho \abs{f}^2
      %   e^{-\vphi_L-\phi_D}}{\logpole}
      \mathllap{\overset{\text{\alert[Gray]{$\vphi_L$ analytically}}}{\underset{\text{\alert[Gray]{singular}}}=}}
      &~\lim_{\rho \descendsto \charfct_{\cl V}} \lim_{\eps \tendsto
        0^+}
      \eps \int_{V'} \frac{
        \rho \abs f^2 \:e^{-\vphi_L -\phi_D}
      }{\abs{\psi_D}^{\sigma +\eps}}
      < +\infty \; ,\footnotemark
    \end{aligned}
  \end{equation*}%
  \footnotetext{\label{fn:explain_res-weights}%
    \newtext*{At the time of writing, it is still an ongoing project in
    investigating the behaviour of the residue functions when
    $\vphi_L$ has more general singularities.
    The assumption of analytic singularities on $\vphi_L$ is
    emphasised here for accuracy.
    Moreover, according to \cite{Chan_adjoint-ideal-nas}*{Thm.~4.1.2
      (2)} (or the computation in
    \cite{Chan_on-L2-ext-with-lc-measures}*{Prop.~2.2.1 and
      Cor.~2.3.3}), the same limit can also be achieved when the weight
    $\frac{e^{-\vphi_L-\phi_D}}{\abs{\psi_D}^{\sigma+\eps}}$ is
    replaced by $\frac{e^{-\vphi_L-\phi_D}}{\logpole|\psi_D|}$.
    It seems more natural to define the residue norm via the latter
    weight considering the given definition of adjoint ideal
    sheaves, but the current definition of residue norm is sufficient
    for the purpose of this paper.%
    %   It is shown in the residue computations in both
    % \cite{Chan_on-L2-ext-with-lc-measures}*{Prop.~2.2.1} (for
    % weights like ) and
    % \cite{Chan&Choi_ext-with-lcv-codim-1}*{Prop.~2.2.1} (for
    % weights like )
    % that, under the assumption $\vphi_L$ having only neat analytic
    % singularities, the limit gives the same residue norm.
    }
  }%
  where $f$ is assumed to be defined on a neighborhood $V'$ of the closure
  $\cl V$ of $V$ and $\rho \colon V' \to [0,1]$ is a compactly
  supported smooth function identically equal to $1$ on $V$.
  The limit $\lim_{\rho \descendsto \charfct_{\cl V}}$ refers to
  the pointwise limit as $\rho$ descends to the characteristic
  function $\charfct_{\cl V}$ of $\cl V$ on $X$.
  %%%%% \emph* is needed as the package embrac is used and \logKX[L]
  %%%%% appears inside \emph.
  Such a norm is referred to as the \emph*{residue norm on $\logKX.[L]
    \otimes \residlof$ on $V$}, or sometimes simply the residue norm
  on $\lcc<V>' := \lcc' \cap V$.
  The residue morphism $\Res^\sigma$ is then given in
  \cite{Chan_adjoint-ideal-nas}*{\S 4.2} by 
  \begin{equation*}
    \renewcommand{\objectstyle}{\displaystyle}
    \xymatrix@C+0.5cm@R-0.5cm{
      {\logKX.[L] \otimes \aidlof} \ar[r]^-{\Res^\sigma}
      \ar@{}[d]|*[left]+{\in} 
      & {\logKX.[L] \otimes \residlof}
      % \save +<4em,-1.3ex>*{\logKX[L] \otimes \residlof} \restore
      \ar@{}[d]|*[left]+{\in}
      % & *+<-2cm,-1cm>{}
      % \ar@{}[l]|(.41)*+{}
      & *+<-4.5em,0em>{\smash[b]{\bigoplus_{p \in\Iset} K_{\lcS} \otimes \res L_{\lcS}}\otimes
        \mtidlof<\lcS>{\vphi_L}}
      \ar@{}[l]|(.5)*+{=}
      \\
      *+<0.8cm,0cm>{f} \ar@{|->}[r]
      & {\paren{\res{g_p}_{\lcS}}_{\mathrlap{p\in\Iset}}
        \mathrlap{\hphantom{p\in\Iset} .}} 
    }
  \end{equation*}

  The above equation of the residue norm works also for $f$ with coefficients in
  $\smooth_{X}$.
  % , where
  % \begin{align*}
  %   \smooth_{X\, *}
  %   &:=\paren{\smooth_{X}\left[
  %     \frac{1}{\abs{\sect_i}} \colon i \in \Iset||
  %     \right]}_{\text{b}}
  %     \qquad\paren{\sect_i \text{ treated as a local defining function of }
  %     D_i} \\
  %   &:=\set{\text{locally bounded elements in the $\smooth_X$-algebra generated
  %     by } \frac{1}{\abs{\sect_i}} \text{ for all } i\in\Iset||} \;
  %     .\footnotemark
  % \end{align*}%
  % \footnotetext{See \cite{Chan&Choi&Matsumura_injectivity}*{footnote 2}.}%
  % % \footnotetext{
  % %   On an admissible open set $V$ under the holomorphic coordinate
  % %   system $(z_1,\dots, z_n)$ such that $D\cap V =\set{z_1 z_2 \dotsm
  % %     z_{\sigma_V} =0}$, one has
  % %   \begin{equation*}
  % %     \smooth_{X \,*}(V)
  % %     =\smooth_X(V)\left[e^{\pm \cplxi \theta_1}, \dots, e^{\pm \cplxi
  % %         \theta_{\sigma_V}} \right]
  % %   \end{equation*}
  % %   where $(r_j,\theta_j)$ is the polar coordinate system of the
  % %   $z_j$-plane for $j=1,\dots,\sigma_V$ in $V$, which is (almost) the
  % %   same as the ad hoc definition of $\smooth_{X\, *}(V)$ given in 
  % %   \cite{Chan&Choi_injectivity-I}*{\S 2.6} (in which
  % %   $e^{\pm\cplxi\theta_{k}}$ for $k \geq \sigma_V +1$ are also included
  % %   in the set of generators of the algebra).
  % %   The definition given here is independent of coordinates and its
  % %   sheaf structure can be seen easily.
  % % }%
  The coefficients of $\Res^\sigma$ (and hence $\PRes[\lcS]$ for any
  $p\in\Iset$) can be extended from $\holo_X$ to $\smooth_{X}$
  (also to $\smform/0,q/[X\,*]$, the sheaf of germs of smooth
  $(0,q)$-forms) accordingly.
  The residue norm is finite when the coefficients of $f$ belong to
  $\smooth_{X} \cdot \aidlof$ on $\cl V$.

}

%%%%%%%%%%%%%%%%%%%%%%%%%%%%%%%%%%%%%%%%%%%%%%%%%%%%%

{
  \setDefaultMetric{\rs\omega}

  We recall below the residue formula from
  \cite{Chan&Choi&Matsumura_injectivity}*{Prop.~2.3.3} 
  for $(L,\vphi_L)$ and adapt it to the current setting.
  For every $\sigma$-lc center $\lcS \subset \lcc'$, we write 
  \begin{equation*}
    \lcS* := \lcS \cap X^\circ \; , \quad
    \lcS<c> := \lcS \cap X_c \; , \;\text{ and }\quad
    \lcS*<c> := \lcS \cap X_c^\circ \; . 
  \end{equation*}
  and set
  \begin{equation} \label{eq:def-Takagoshi-harmonic-sp-on-lcS}
    \Harm<\lcS>,{\vphi_L}(c)
    := \Harm<\alert{\lcS<c>}>{\logKX<\lcS>[L];\Phi},{\vphi_L}
  \end{equation}
  as the Takegoshi harmonic space on $\lcS<c>$, for convenience.

  First, we recall a basic local residue computation.
  Let \newtext{$\smooth_{X\: \text{cpt}}$} denote the sheaf of germs of smooth functions
  on $X$ with compact support and let the ad hoc notation
  ``$\smooth_X[\rs\omega^{\pm}]$'' temporarily denote the algebra in
  Proposition \ref{prop:singularities-along-FM} for convenience.
  Further let ``$\smooth_X[\rs\omega^{\pm}] * \mtidlof{\vphi_L}$'' denote
  the multiplier ideal sheaf defined in $\smooth_X[\rs\omega^{\pm}]$
  in place of $\holo_X$, i.e.~the ideal sheaf of germs of functions in
  $\smooth_X[\rs\omega^{\pm}]$ which are locally $L^2$ with respect to
  $\vphi_L$.
  Define also ``$\smooth_X[\rs\omega^{\pm}] * \aidlof{\vphi_L}$'' similarly.

  \begin{prop}[\cite{Chan&Choi&Matsumura_injectivity}*{Prop.~2.3.2}]
    \label{prop:residue-formula-classical-kernel}
    
    Given any admissible open set $V \subset X$ and any compactly
    supported section $f \in \logKX.[L] \otimes \newtext{\smooth_{X
      \:\text{cpt}}}[\rs\omega^{\pm}] * \aidlof|1|{\vphi_L}\paren{V}$
    such that $\Res^1(f) = g =\paren{g_b}_{b\in\Iset|1|}$, we have,
    for any $\xi \in \logKX.[L] \otimes \smooth_{X}[\rs\omega^{\pm}]
    * \mtidlof{\vphi_L} \paren{V}$,
    \begin{align*}
      \lim_{\eps \tendsto 0^+} \eps \int_V
      \inner{\xi}{f} \:e^{-\phi_D-\vphi_L} e^{-\eps\abs{\psi_D}}
      &=\sum_{b \in \Iset|1|} \pi
        \int_{\lcS|1|[b] \cap V} \inner{\frac{\rs*\xi_b}{\sect_{(b)}}}{\: g_b}
        \:e^{-\vphi_L} \\
      &=\sum_{b \in \Iset|1|}
      % \underbrace{
        \pi
        \int_{\lcS|1|[b] \cap V} \inner{\rs*\xi_b}{\: g_b \sect_{(b)}}
        \:e^{-\phi_{(b)}-\vphi_L}
        % }_{\displaystyle =:
        %   \iinner{\rs*\xi_b}{g_b\sect_{(b)}}_{\mathrlap{\lcS|1|[b] \cap V,
        %   \phi_{(b)}}}}
    \end{align*}
    which is finite,
    where $\phi_{(b)} :=\log\abs{\sect_{(b)}}^2$ and
    \begin{equation*}
      \rs*\xi_b := \PRes[\lcS|1|[b]](\frac{\xi}{\sect_D}) \cdot \sect_{(b)}
      \in K_{\lcS|1|[b]} \otimes \Diff_b D \otimes \res L_{\lcS|1|[b]} \otimes
      \smooth_{\lcS|1|[b]\:c}[\rs\omega^{\pm}] * \mtidlof<\lcS|1|[b]>{\vphi_L} \paren{\lcS|1|[b] \cap V} \; .
    \end{equation*}    
  \end{prop}
  \begin{proof}[Remarks on the proof]
    The proof follows the same approach as that in
    \cite{Chan&Choi&Matsumura_injectivity}*{Prop.~2.3.2} with the
    relaxation of the coefficients of $f$ and $\xi$ from $\smooth_X$
    to $\smooth_X [\rs\omega^{\pm}]$.
    This change does not affect the residue computation due to Fubini's
    theorem.

    Note also that, in
    \cite{Chan&Choi&Matsumura_injectivity}*{Prop.~2.3.2}, the
    coefficients are allowed to be in $\smooth_{X\,*}$ (locally bounded germs in
    $\smooth_X\left[\frac 1{\abs{\sect_i}} \colon  i\in \Iset||
    \right]$, see \cite{Chan&Choi&Matsumura_injectivity}*{\S 2.3 and
      footnote 2}).
    Such coefficients are convenient when handling
    residues of $\dbar\psi_D \otimes u$ (see
    \cite{Chan&Choi_injectivity-I}*{proof of Prop.~3.2.3}), but are
    not necessary for dealing with residues of
    $\idxup{\diff\psi_D}. u$, as in this paper.
  \end{proof}

  All statements made above and in previous sections concerning
  $(X_c,D)$ for $c \in (0,\infty]$ (where $X_\infty = X$) are equally
  applicable to $(\lcS<c>, \Diff_p D)$ for any $p \in \Iset$.

  To adapt the global residue formula in
  \cite{Chan&Choi&Matsumura_injectivity}*{Prop.~2.3.3} to the
  non-compact setup, an additional assumption is required, namely, the
  involved harmonic form has to satisfy the Takegoshi property (see Theorem
  \ref{thm:Takegoshi-argument}).

  Recall from
  \cite{Chan&Choi&Matsumura_injectivity}*{Prop.~2.3.3} that, for any
  $\sigma$-lc center $\lcS$ and $(\sigma+1)$-lc center $\lcS+1[b]$ such
  that $\lcS+1[b] \subset \lcS$, the sign $\sgn{b:p}$ is defined by
  \begin{equation*}
    \PRes[\lcS+1[b]] =\sgn{b:p} \:\PRes[\lcS+1[b] | \lcS] \circ
    \PRes[\lcS] \; ,
  \end{equation*}
  where $\PRes[\lcS+1[b] | \lcS]$ denotes the Poincar\'e residue map
  from $\lcS$ to $\lcS+1[b]$.
  The ambient space in the following statement is assumed
  to be $X_c$ for any $c \in (0,\infty]$.
  % $\Phi$ is a smooth psh exhaustion function with
  % $\abs{\diff \Phi}_\omega$ bounded on $X$.
  % Define a sequence $\seq{\eta_\nu}_{\nu \in \Nnum}$ cut-off
  % functions by 
  % \begin{equation*}
  %   \eta_\nu := \rho\paren{\frac{\Phi}{\nu}} \; ,
  % \end{equation*}
  % where $\rho \colon [0,+\infty) \to [0,1]$ is a compactly supported
  % non-increasing smooth function such that $\res{\rho}_{[0,1]} \equiv 1$,
  % $\res\rho_{[2,+\infty)} \equiv 0$ and also $\abs{\rho'} \lesssim
  % 1$ on its domain ($\rho'$ denotes the derivative of $\rho$).
  % It can be seen that $\seq{\supp \eta_\nu}_{\nu\in\Nnum}$ is a
  % compact exhaustion of $X$ (which induces a compact exhaustion on
  % every lc center of $(X,D)$) and $\eta_\nu \ascendsto 1$ as $\nu
  % \ascendsto +\infty$.

  \begin{prop} \label{prop:res-formula-dbar-exact-dot-harmonic}
    Let $u_p \in \Harm<\lcS>(c)$
    % be a \emph{harmonic} $K_{\lcS} \otimes \res L_{\lcS}$-valued
    % $(0,q)$-form on $\lcS$ with respect to the
    (with norm $\norm\cdot_{\lcS*<c>} := \norm\cdot_{\lcS*<c>, \vphilist}$)
    % satisfying $\idxup{\diff\Phi} . u_p = 0$ on $\lcS$ (the Takegoshi property)
    for each $p \in \Iset$.
    With the finite cover $\cvr V$ and partition of unity
    $\set{\rho^i}_{i \in I}$ given in Section \ref{subsec:notation},
    let $\set{\gamma_{\idx 1.q}}_{\idx 1,q \in I}$ be a
    $\logKX. \newtext{{}\otimes \mtidlof<X>{\vphi_F}}$-valued \v Cech $(q-1)$-cochain with respect to $\cvr
    V \cap X_c$ and set, for each $p \in \Iset$ and $b \in \Iset+1$, 
    \begin{equation*}
      \rs\gamma_{\bullet; \:\idx 1.q} :=\PRes[\lcS||[\bullet]](\frac{\gamma_{\idx 1.q}}{\sect_D})
      \cdot \sect_{(\bullet)} \;\; \text{ and } \;\;
      v_{\bullet} := \sum_{\idx 1,q \in I} \underbrace{
        \dbar\rho^{i_q} \wedge \dotsm
        \wedge \dbar\rho^{i_2} \cdot \rho^{i_1}
      }_{=: \: \paren{\dbar\rho}^{\idx q.1}} \rs*\gamma_{\bullet;\:\idx 1.q}
      \quad\text{ on } \lcS||<c>[\bullet] \; ,
      % \\
      % \text{and }\quad
      % \rs\gamma_{b; \:\idx 1.q} :=\PRes[\lcS+1[b]](\frac{\gamma_{\idx 1.q}}{\sect_D})
      % \cdot \sect_{(b)} \; , \quad
      % v_{b} := \sum_{\idx 1,q \in I} \paren{\dbar\rho}^{\idx q.1}
      % \rs*\gamma_{b;\:\idx 1.q}
      % \quad\text{ on } \lcS+1[b] \; .
    \end{equation*}
    \newtext{where $\bullet$ denotes $p$ or $b$ and $\lcS||<c>[\bullet]$ denotes $\lcS<c> $, $\lcS+1<c>[b]$ respectively.} 
    Assume that $\frac{\dbar v_p}{\sect_{(p)}}$ (but not necessarily
    for $\frac{v_p}{\sect_{(p)}}$) is locally $L^2$ in
    $\lcS<c>$ (not only in $\lcS*<c>$) with respect to $\vphi_F$ and $\rs\omega$.
    Then,
    % after setting $\iinner{\cdot}{\cdot}_{\lcS||[\bullet], \phi_{(\bullet)}}
    % :=\iinner{\cdot}{\cdot \:e^{-\phi_{(\bullet)}}}_{\lcS||[\bullet]}$,
    we have, for every $\nu \in \Nnum$,
    \begin{equation*}
      \sum_{p\in\Iset} \iinner{\alert{\eta_{c,\nu}} \frac{\dbar v_{p}}{\sect_{(p)}} \:}{
        \: u_p
      }_{\lcS*} 
      =-\sigma_+ \smashoperator[l]{\sum_{b\in\Iset+1}}
      \iinner{\alert{\eta_{c,\nu}} \frac{v_{b}}{\sect_{(b)}} \:}{\quad\;\;
        \smash{\smashoperator{
            \sum_{p\in\Iset \colon \lcS+1[b] \subset \lcS}
          }} \;\;
        \sgn{b:p} \:
        \PRes[\lcS+1[b] | \lcS](\idxup{\diff\psi_{(p)}}. u_p)
      }_{\lcS*+1[b]} \; ,
    \end{equation*}
    where $\psi_{(p)} :=\phi_{(p)} -\sm\vphi_{(p)}$ and
    $\sm\vphi_{(p)}$ is some smooth potential on $\Diff_p D$, and
    \mmark{$\sigma_+ := \max\set{1, \sigma}$.}{
      $\sigma_+$ is used so that the statement is still valid when
      $\sigma = 0$.
    }
    Moreover, the equality also holds when each $u_p$ is replaced
    by $s u_p$ for all $p \in \Iset$ (and $\seq{\gamma_{\idx
        1.q}}_{\idx 1,q \in I}$ is taken as $\logKX. M
    \newtext{{}\otimes \mtidlof<X>{\vphi_F+\vphi_M}}$-valued and the norm
    is replaced by $\norm\cdot_{\lcS*||[\bullet],\vphilist M}$).
  \end{prop}

  \begin{proof}
    The chain rule and the Takegoshi property
    (i.e.~$\idxup{\diff\Phi}. u_p = 0$) ensure that
    $\idxup{\diff\eta_{c,\nu}} . u_p = 0$ on $\lcS<c>$.
    Apart from the introduction of the cut-off functions $\eta_{c,\nu}$
    into the formula, the proof is the same as that of
    \cite{Chan&Choi&Matsumura_injectivity}*{Prop.~2.3.3}.
    The treatment below focuses on the arguments involving $\eta_{c,\nu}$ and
    the description of the residue computation is kept brief.
    Readers are referred to \cite{Chan&Choi&Matsumura_injectivity} for details.
    
    Set $\iinner{\cdot}{\cdot}_{\lcS*||[\bullet], \phi_{(\bullet)}}
    :=\iinner{\cdot}{\cdot \:e^{-\phi_{(\bullet)}}}_{\lcS*||[\bullet]}$
    for $\bullet = p, b$.
    The smooth form $v_{p}$ on $\lcS<c>$ need not be locally $L^2$ with
    respect to the weight $e^{-\phi_{(p)}}$, so
    an integration by parts is done via the use of Proposition
    \ref{prop:residue-formula-classical-kernel} (with the knowledge of
    the singularities of $u_p$ along $\lcS \setminus \lcS*$ by
    Proposition \ref{prop:singularities-along-FM}), which yields \label{page:explain_limit-order}
    \begin{align*} 
      &~\sum_{p\in\Iset} \iinner{\alert{\eta_{c,\nu}} \frac{\dbar v_{p}}{\sect_{(p)}} \:}{
        \: u_p
        }_{\lcS*}
        =\sum_{p\in \Iset} \iinner{\alert{\eta_{c,\nu}} \dbar v_{p}}{ u_p
        \sect_{(p)}}_{\lcS*, \phi_{(p)}}
      \\ 
      \newtext{\xleftarrow{\eps \tendsto 0^+}}
      &~\newtext{\sum_{p \in \Iset} \iinner{
        e^{-\eps \abs{\psi_{(p)}}} \:\alert{\eta_{c,\nu}} \dbar v_{p}
        }{ u_p \sect_{(p)}}_{\lcS*, \phi_{(p)}}}
      \\
      \newtext{=}
      &~\newtext{\sum_{p\in \Iset} \paren{
        \iinner{
        e^{-\eps \abs{\psi_{(p)}}}\:
        \dbar\paren{\alert{\eta_{c,\nu}}  v_{p}}
        }{ u_p \sect_{(p)}}_{\lcS*,\phi_{(p)}}
        % -\cancelto{
        % \mathllap{(\because \text{ Takegoshi property})\quad}
        % 0 
        % }{\iinner{e^{-\eps \abs{\psi_{(p)}}} \:v_{p}}{
        % \idxup{\diff\alert{\eta_{c,\nu}}}. u_p \sect_{(p)}
        % }}_{\lcS*,\phi_{(p)}}
        -\iinner{e^{-\eps \abs{\psi_{(p)}}} \:v_{p}}{
        \smash{\cancelto{
        \mathllap{(\because \text{ Takegoshi property})\quad}
        0 
        }{\idxup{\diff\alert{\eta_{c,\nu}}}. u_p \sect_{(p)}}}
        }_{\lcS*,\phi_{(p)}}
        }}
      \\
      =&~\sum_{p \in \Iset} \paren{
         \cancelto{0 \;\;\;(\because~u_p \text{ harmonic, Lemma
         \ref{lem:su-harmonicity} or Proposition \ref{prop:Takegoshi-harmonic-forms}})}{\iinner{
         \dbar\paren{e^{-\eps \abs{\psi_{(p)}}} \:\alert{\eta_{c,\nu}} v_{p}}
         }{ u_p \sect_{(p)}}_{\mathrlap{\lcS*, \phi_{(p)}}}}
         \quad\;\; - \eps 
         \iinner{
         e^{-\eps \abs{\psi_{(p)}}} \:\alert{\eta_{c,\nu}} v_{p}
         }{\:\idxup{\diff\psi_{(p)}}.  u_p \sect_{(p)}}_{\lcS*,
         \phi_{(p)}}
         }
      \\
      =&~-\smashoperator[l]{\sum_{\substack{\idx 1,q \in I \, , \\ p \in \Iset}}}  \eps \:
         \iinner{
         e^{-\eps \abs{\psi_{(p)}}} \:\alert{\eta_{c,\nu}} % \paren{\dbar\rho}^{\idx q.1}
         \rs*\gamma_{p;\:\idx 1.q}
         }{\:
         \idxup{\diff\rho},[\idx 1.q] .
         \paren{\idxup{\diff\psi_{(p)}}.  u_p \sect_{(p)}}
         }_{\lcS*, \phi_{(p)}}
      \\
      \xrightarrow[
      \substack{\text{Prop.~\ref{prop:singularities-along-FM}}
      \\ \text{Prop.~\ref{prop:residue-formula-classical-kernel}} 
      }]{\eps \tendsto 0^+} 
      &~-\smashoperator[l]{\sum_{\substack{p \in \Iset \, ,\\ \idx 1,q \in I}}} 
        \sum_{k=\sigma +1}^{\mathclap{\sigma_{V_{\idx 1.q}}}} \sigma_+
        \iinner{ \alert{\eta_{c,\nu}}
        \PRes[p(k)](
        \frac{\rs*\gamma_{p;\:\idx 1.q}}{\sect_{(p)}}
        )
        }{\:
        \idxup{\diff\rho},[\idx 1.q] .
        \PRes[p(k)](\idxup{\diff\psi_{(p)}}.  u_p)
        }_{\lcS* \cap \set{z_{p(k)} =0}}
        \; ,
    \end{align*}
    where $\idxup{\diff\rho},[\idx 1.q] .\cdot$ is the adjoint
    of $\paren{\dbar\rho}^{\idx q.1} \cdot$, and $\PRes[p(k)]$ denotes
    the Poincar\'e residue map from $\lcS$ to $\lcS \cap
    \set{z_{p(k)}=0}$, in which $(z_1, \dots, z_n)$ is a holomorphic
    coordinate system such that $\lcS \cap V_{\idx 1.q}
    =\set{z_{p(1)} = \dotsm =z_{p(\sigma)} =0}$ and
    $\sect_{(p)} =z_{p(\sigma+1)} \dotsm z_{p(\sigma_{V_{\idx
          1.q}})}$.
    \newtext{\label{page:explain_assumption-met}Proposition \ref{prop:residue-formula-classical-kernel}
      is applicable in the last step as each $\rs\gamma_{p;\:\idx 1.q}$
      lies in $\mtidlof<\lcS>{\vphi_F}$ on $\lcS \cap V_{\idx 1.q}$ and the
      singularities of $\idxup{\diff\psi\ps} . u_p\sect\ps$ along
      $\lcS \setminus \lcS*$ guaranteed by Proposition
      \ref{prop:singularities-along-FM} fit in the assumption.}
    Notice that the introduction of $\eta_{c,\nu}$ helps to avoid the need to
    handle the boundary $\bdry X_c$ in the local residue computation.
    Note also that the appearance of the coefficient $\sigma_+$ comes from
    the different normalizations of the $L^2$ norms on various lc centers,
    namely, $\norm\cdot_X^2 := \int_X \dotsm$ and $\norm\cdot_{\lcS}^2
    := \frac{\pi^\sigma}{(\sigma -1)!} \int_{\lcS} \dotsm$ for every
    integer $\sigma \geq 1$.
    % With the assumption that $\ibddbar\vphi_F \geq 0$ on $X$, the
    % proof of \cite{Chan&Choi_injectivity-I}*{Prop.~3.3.2} is used 
    % in the last step to ensure that the singularities of
    % $\idxup{\diff\rho}[*, \idx 1.q]
    % \PRes[p(k)](\idxup{\diff\psi_{(p)}} . u_p)$ along $\lcS<c>
    % \setminus \lcS*<c>$ do not interfere with the residue
    % computation. 

    Following the argument in the proof of
    \cite{Chan&Choi&Matsumura_injectivity}*{Prop.~2.3.3}, the $(\sigma
    +1)$-lc centers $\lcS \cap \set{z_{p(k)} = 0}$ for $k = \sigma +1
    , \dots, \sigma_{V_{\idx 1.q}}$ in $V_{\idx 1.q}$
    can be re-indexed in terms of $b \in \Iset+1$ such that
    \begin{equation*}
      \lcS[p_{b,j}] \cap \set{z_{b(j)} = 0}
      = \lcS+1[b] \cap V_{\idx 1.q}
      \quad\text{ for } j = 1, \dots, \sigma +1 
    \end{equation*}
    and the summations transform as $\sum_{p \in \Iset} \sum_{k=\sigma
      +1}^{\sigma_{V}} \dotsm = \sum_{b \in
      \Iset+1} \sum_{j=1}^{\sigma +1} \dotsm$.
    With such a choice of indexing, we have
    \begin{equation*}
      \frac{\rs*\gamma_{b;\: \idx 1.q}}{\sect_{(b)}}
      :=\PRes[\lcS+1[b]](\frac{\gamma_{\idx 1.q}}{\sect_D})
      =\sgn{b:p_{b,j}} \:
      \PRes[b(j)](\frac{\rs*\gamma_{p_{b,j};\:\idx
          1.q}}{\sect_{(p_{b,j})}})
    \end{equation*}
    (noticing that % $\sect_{(b)} =\sect_{(\sigma+1 : b)}$,
    % $\sect_{(p_{b,j})} =\sect_{(\sigma : p_{b,j})}$ and
    $\sect_{(p_{b,j})} = z_{b(j)} \sect_{(b)}$).
    As a result, the expression in question becomes
    \begin{align*}
      &-\smashoperator[l]{\sum_{\substack{b \in \Iset+1  ,\\ \idx 1,q \in I}}}
        \sum_{j=1}^{\sigma +1} \sigma_+
        \iinner{ \sgn{b:p_{b,j}}\:\alert{\eta_{c,\nu}}
        \frac{\rs*\gamma_{b;\:\idx 1.q}}{\sect_{(b)}}
        }{\: 
        \idxup{\diff\rho},[\idx 1.q] .
        \PRes[b(j)](\idxup{\diff\psi_{(p_{b,j})}}.  u_{p_{b,j}})
        }_{\lcS*+1[b]}
      \\
      =&-\smashoperator[l]{\sum_{\substack{\idx 1,q \in I \, , \\ b \in \Iset+1}}}
         \sigma_+
         \iinner{\alert{\eta_{c,\nu}}
         \frac{\paren{\dbar\rho}^{\idx q.1} \rs*\gamma_{b;\:\idx 1.q}}{\sect_{(b)}}
         \:}{  \smash{\sum_{j=1}^{\sigma +1}} \sgn{b:p_{b,j}}\:
         \PRes[b(j)](\idxup{\diff\psi_{(p_{b,j})}}.  u_{p_{b,j}})
         }_{\lcS*+1[b]}
      \\
      =&-\sigma_+ \sum_{b \in \Iset+1} \iinner{
         \alert{\eta_{c,\nu}} \frac{v_b}{\sect_{(b)}}
         \:}{\quad\;
         \smash{\smashoperator{\sum_{p\in\Iset \colon \lcS+1[b] \subset
         \lcS}}} \;\;
         \sgn{b:p}\:
         \PRes[\lcS+1[b] | \lcS](\idxup{\diff\psi_{(p)}}.  u_{p})
         }_{\lcS*+1[b]} \; . 
    \end{align*}

    Note that the singularities of $s u_p$ along $\lcS<c> \setminus
    \lcS*<c>$ (and thus those singularities of $\idxup{\diff\rho},[\idx 1.q].
    \PRes[p(k)](\idxup{\diff\psi_{(p)}} . s u_p)$) do not
    interfere with the residue computation, so the above derivation,
    and the resulting equality remains valid with $s u_p$ in
    place of $u_p$.
  \end{proof}

}

\subsection{Harmonic residues} % $(w_b)_{b \in \Iset+1}$
\label{sec:residue-element}

%\input{residue-element}

%%%%%
%%%%% File name  : residue-element.tex
%%%%% Author     : Mario Chan
%%%%% Date       : 23rd May, 2024
%%%%% Description: The section on the harmonicity of w and the L2-ness
%%%%%              of its image under the connecting morphism.
%%%%%
%%
%%%

{
  \setDefaultMetric{\rs\omega}
  % \setDefaultAmbientSpace{X^\circ}

  In view of Proposition
  \ref{prop:res-formula-dbar-exact-dot-harmonic}, given any collection
  of $u :=(u_p)_{p \in \Iset}$ of harmonic forms $u_p \in
  \Harm<\lcS>(c)$ on $\lcS*<c> := \lcS \cap X_c^\circ$
  for each $p \in \Iset$, define
  \begin{equation} \label{eq:definition-of-w}
    w_b := \smashoperator[r]{
      \sum_{p\in\Iset \colon \lcS+1[b] \subset \lcS}
    } \;\;
    \sgn{b:p}\:
    \PRes[\lcS+1[b] | \lcS](\idxup{\diff\psi_{(p)}}.  u_{p})
    \quad\text{ on $\lcS*+1<c>[b]$ for each } b \in \Iset+1
  \end{equation}
  and set
  \begin{equation*}
    \HRes(u) := w := \paren{w_b}_{b \in \Iset+1}  \; ,
  \end{equation*}
  which is referred to as the \emph{harmonic residue} of $u$.
  The naming is justified by the following result.
  \begin{thm} \label{thm:Takegoshi-harmonicity-of-HRes}
    For any $c \in (0,\infty]$, the map $\HRes$ is a bounded linear
    operator between Takegoshi harmonic spaces
    \begin{equation*}
      \HRes \colon
      \bigoplus_{p \in \Iset} \Harm<\lcS>(c)
      \to
      \bigoplus_{b \in \Iset} \Harm/q-1/<\lcS+1[b]>(c) \; ,
    \end{equation*}
    that is, given any $(u_p)_{p\in\Iset} \in \bigoplus_{p \in \Iset}
    \Harm<\lcS>(c)$, each $w_b$ in \eqref{eq:definition-of-w} is a
    $\logKX<\lcS+1[b]>[\res F_{\lcS+1[b]}]$-valued harmonic $(0,q-1)$-form with respect
    to $\vphi_F$ and $\rs\omega$ on $\lcS*+1<c>[b]$ satisfying the
    Takegoshi property $\idxup{\diff\Phi}.w_b = 0$.
  \end{thm}

  \begin{proof}
    Let $u := (u_p)_{p \in \Iset} \in \bigoplus_{p\in\Iset}
    \Harm<\lcS>(c)$ and $w := (w_b)_{b\in \Iset+1} := \HRes(u)$.
    By the local computations in
    \cite{Chan&Choi&Matsumura_injectivity}*{Prop.~2.4.1 and Lemma
      2.4.2} (with $(\lcS*<c>, \lcS*+1<c>[b])$ in place of $(X, D_p)$
    there), under the assumption that $\ibddbar\vphi_F \geq 0$ (hence
    $\nabla^{(0,1)} u_p = 0$ by Proposition
    \ref{prop:consequence-of-positivity} or Proposition
    \ref{prop:Takegoshi-harmonic-forms}), we have 
    \begin{equation*}
      \dfadj w_b = 0 \quad\text{ and }\quad
      \dbar w_b = 0
      \quad\text{ on } \lcS*+1<c>[b]
    \end{equation*}
    (note that the computations work on non-compact spaces by their
    local nature).
    Since $\ibddbar\sm\vphi_{(p)}$ is bounded from above on $\lcS$, 
    Proposition \ref{prop:Takegoshi-property-on-tBK} with $(\lcS<c>,
    \lcS+1<c>[b], \psi_{(p)}, \PRes[\lcS+1[b] | \lcS], u_p)$ in place of $(X_c, \lcS|1|<c>,
    \psi_D, \PRes[\lcS|1|[b]], u)$ implies that
    $\PRes[\lcS+1[b] | \lcS](\idxup{\diff\psi_{(p)}}.u_p)$ is $L^2$ with
    respect to $\vphi_F$ and $\rs\omega$ on $\lcS*+1<c>[b]$ and is a
    bounded linear map in $u_p$.
    The definition \eqref{eq:definition-of-w} of $w_b$ then ensures
    that $w_b \in
    \Ltwo/0,q-1/<\lcS+1<c>[b]>{\logKX<\lcS+1<c>[b]>}_{\vphilist}$ and
    $u \mapsto w = \HRes(u)$ is a bounded linear map.
    % It remains to check that $w_b$ is globally $L^2$
    % on $\lcS*+1[b]$ with respect to $\res{\vphi_F}_{\lcS+1[b]}$ and
    % $\res{\rs\omega}_{\lcS+1[b]}$ and is also sitting in the domain
    % $\Dom\dbadj$ of the Hilbert space adjoint of $\dbar$ on
    % $\lcS*+1[b]$, therefore confirming that
    % $w_b \in \Harm/q-1/<\lcS+1[b]>{\logKX<\lcS+1[b]>}$.
    % This is a consequence of the following proposition.

    It remains to check that $\idxup{\diff\Phi}. w_b = 0$ on
    $\lcS*+1<c>[b]$.
    On any admissible open set $V \subset X$ with $\emptyset \neq \lcS+1[b] \cap V
    \subset \lcS \cap V$, let $z_1$ be a
    coordinate function on $V$ such that $dz_1$ generates the conormal
    bundle $N_{\lcS+1[b] | \lcS}^*$ on $\lcS*+1<c>[b] \cap V$ and write
    $u_p = dz_1 \wedge \rs u_p$ (where $\rs u_p$ contains no $dz_1$).
    Note that $\PRes[\lcS+1[b] | \lcS](\idxup{\diff\psi_{(p)}}.u_p) =
    \res{\idxup{dz_1}.\rs*u_p}_{\lcS+1[b]}$\footnote{\label{fn:explain_10-01-commut}%
      \newtext{% For readers who may get confused when they compare the
        % notation ``$\rs u_p$'' here with ``\rs u_{1,2}'' or ``$\rs
        % u_{p,k}$'' in \cite{Chan&Choi&Matsumura_injectivity}*{Step 3
        %   of pf.~of Thm.~1.2} or \cite{Chan&Choi&Matsumura_injectivity}*{Step III
        %   of pf.~of Thm.~3.4.1}:
        The differential form $u_p$ on $\lcS*<c>$ are treated
        as a section of $K_{\lcS<c>} \otimes \cTgt/0,q/_{\lcS<c>} \isom
        \cTgt/n-\sigma,0/_{\lcS<c>} \alert{{}\otimes{}}  
        \cTgt/0,q/_{\lcS<c>}$ instead of $\cTgt/n-\sigma, q/_{\lcS<c>} =\cTgt/n-\sigma,0/_{\lcS<c>}
        \alert{{}\wedge{}} \cTgt/0,q/_{\lcS<c>}$ so that, when
        contracted with $(0,1)$-vector fields like $\idxup{dz_1}$ (which
        interact only with the $(0,q)$-form), we can safely write
        $\idxup{dz_1} . u_p = \idxup{\diff\psi\ps} . dz_1 \wedge
        \rs u_p = dz_1 \wedge \idxup{dz_1} . \rs u_p$ without
        worrying about the sign when commuting $dz_1 \wedge \cdot$ and
        $\idxup{dz_1} . \cdot$.
        Other differential forms like $w_b$ on $\lcS*+1<c>[b]$ or $u$
        on $X$ are handled similarly.}
    }
    on $\lcS*+1<c>[b] \cap V$.
    The Takegoshi property $\idxup{\diff\Phi}.u_p = 0$ therefore
    implies that
    \begin{equation*}
      \begin{aligned}
        0 =\idxup{dz_1}. \idxup{\diff\Phi}. \paren{
          dz_1 \wedge \rs u_p
        }
        &=dz_1 \wedge \paren{
          \idxup{dz_1}. \idxup{\diff\Phi}. \rs* u_p
        }
        \\
        &=-dz_1 \wedge \paren{
          \idxup{\diff\Phi}. \idxup{dz_1}. \rs* u_p
        }
        \quad\text{ on } \lcS*<c> \cap V \; .
      \end{aligned}
    \end{equation*}
    Since $u_p$ is $K_{\lcS}$-valued, i.e.~it has a top holomorphic
    form on $\lcS \cap V$ as a basis, it follows that
    $\idxup{\diff\Phi}. \idxup{dz_1}. \rs* u_p = 0$ on $\lcS*<c> \cap
    V$, and thus $\idxup{\diff\Phi}. \PRes[\lcS+1[b] |
    \lcS](\idxup{\diff\psi_{(p)}}.u_p) = 0$ on $\lcS*+1<c>[b] \cap
    V$.
    After considering all admissible open sets $V$ and all $\sigma$-lc
    centers $\lcS$ intersecting $\lcS*+1<c>[b]$, we conclude that
    $\idxup{\diff\Phi}.w_b = 0$ on $\lcS*+1<c>[b]$.
  \end{proof}

  % As a result, the form $w_b$ for $b \in \Iset+1$ given in
  % \eqref{eq:definition-of-w} is $L^2$ on $\lcS*+1[b]$ with respect
  % to $\res{\vphi_F}_{\lcS+1[b]}$ and $\res{\rs\omega}_{\lcS+1[b]}$.
  % Since $\Dom\dbadj = \Dom\dfadj$ (the domain of the formal adjoint
  % of $\dbar$) on $\lcS*+1[b]$ by the completeness of
  % $\res{\rs\omega}_{\lcS+1[b]}$, this confirms that
  % % $w_b \in \Harm/q-1/<\lcS*+1[b]>{\logKX<\lcS+1[b]>}$.
  % % Considering the homomorphism
  % \begin{equation*}
  %   \xymatrix{
  %   {w_b \in \Harm/q-1/<\lcS+1[b]>{\logKX<\lcS+1[b]>}}
  %   \ar@{^(->}[r]^-{\jmath}
  %   & {\cohgp {q-1}[\lcS+1[b]]{\logKX<\lcS+1[b]> \otimes
  %   \mtidlof<\lcS+1[b]>{\vphi_F}}}
  % } \; ,
  % \end{equation*}
  % where the map $\jmath$ is defined as in
  % \eqref{eq:map-for-harmonic-representatives} (and is injective by
  % Theorem \ref{prop:injective-jmath}), that is, $\jmath(w_b)$
  % represents a locally $L^2$ Dolbeault class in
  % $\cohgp {q-1}[\lcS+1[b]]{\logKX<\lcS+1[b]> \otimes
  % \mtidlof<\lcS+1[b]>{\vphi_F}}$.

  % In order to simplify notation, $\jmath(w_b)$ is identified with
  % $w_b$ in what follows.

  % A collection of forms $w := \paren{w_b}_{b \in \Iset+1}$ described
  % above (identified as the image of $\jmath$) corresponds to a class
  % in
  Recall from %\eqref{eq:map-for-harmonic-representatives} or
  \eqref{eq:map-for-harmonic-representatives-on-X_c} the inclusion 
  \begin{align*}
    \spHarm/q-1/{\residlof+1*}(c)
    :=\smashoperator[r]{\bigoplus_{b \in \Iset+1}}
    \Harm/q-1/<\lcS+1[b]>(c)
    \xhookrightarrow{\;\jmath^c\;}
    &~\smashoperator[r]{\bigoplus_{b \in \Iset+1}}
      \cohgp{q-1}[\lcS+1<c>[b]]{\logKX<\lcS+1[b]> \otimes
      \mtidlof<\lcS+1<c>[b]>{\vphi_F}}
    \\
    =
    &~\cohgp{q-1}[X_c]{\logKX<X>. \otimes \residlof+1*}
      =: \spH/q-1/{\residlof+1*}_c 
  \end{align*}
  \newtext{for any $c \in (0,\infty]$,} where $\residlof+1* := \residlof+1 \isom \faidlof+1/$.
  Identifying the image and pre-image of $\jmath^c$, every $w \in
  \spHarm/q-1/{\residlof+1*}(c)$ can be viewed as
  a class in $\spH/q-1/{\residlof+1*}_c$.
  
  Write $\aidlof* := \aidlof$ for convenience.
  For $q \geq 1$, let
  \begin{equation*}
    \delta := \delta^{q-1} \colon
    \cohgp {q-1}[X_c]{\logKX<X>. \otimes \residlof+1*}
    \to \cohgp q[X_c]{\logKX<X>. \otimes \residlof*}
  \end{equation*}
  be a connecting morphism in the long exact sequence induced from the
  short exact sequence
  \begin{equation*}
    \renewcommand{\objectstyle}{\displaystyle}
    \xymatrix@R=3.5ex{
      {0} \ar[r]
      &{\faidlof/-1*} \ar[r] \ar[d]_-{\Res^{\sigma}}^-{\isom}
      &{\faidlof+1/-1*} \ar[r]
      &{\faidlof+1/*} \ar[r] \ar[d]^-{\Res^{\sigma+1}}_-{\isom}
      &{0 \; .}
      \\
      &{\residlof*} 
      &&{\residlof+1*}
    }
  \end{equation*}
  The image $\delta w$ plays an important role in the proof of the injectivity
  theorem (see \cite{Chan&Choi&Matsumura_injectivity}*{Proof of Thm.~3.4.1, Step
    IV} or Step \ref{item:pf:use_u-ortho-w} of the proof of
  Theorem \ref{T:prime_divisor_case}).
  It can be computed via the \v Cech representative, \mmark{where we take
    $\delta$ to be induced from the \v Cech coboundary operator.}{
    Added to fix a choice of $\delta$.
  }
  For the purpose of this paper, it suffices to consider only the
  images of $\spHarm/q-1/{\residlof+1*}(c)$.

  For each $w =\paren{w_b}_{b\in \Iset+1} \in
  \spHarm/q-1/{\residlof+1*}(c)$, write
  \begin{equation} \label{eq:Cech-Dolbeault-on-w_b}
    w_b \quad\;
    \overset{\mathclap{\text{\eqref{eq:Cech-Dolbeault-isom}}}}= \quad\;
    \dbar v_{b;(2)} + (-1)^{q-1} \frac{v_{b;(\infty)}}{\sect_{(b)}}
    \quad\text{ on } \lcS*+1<c>[b] \; ,
  \end{equation}
  where
  \begin{equation*}
    \frac{v_{b;(\infty)}}{\sect_{(b)}}
    := \paren{\dbar\rho}^{\idx q.1}
    \frac{\rs\gamma_{b;\:\idx 1.q}}{\sect_{(b)}}
    := \paren{\dbar\rho}^{\idx q.1} \alpha_{b;\:\idx 1.q}
    \quad\text{ on } \lcS+1<c>[b] 
  \end{equation*}
  such that $v_{b;(2)} \in
  \Ltwo/0,q-2/<\lcS+1<c>[b]>{\logKX<\lcS+1[b]>[F]}_{\vphilist}$ and
  $\set{\alpha_{b;\:\idx 1.q}}_{\idx 1,q \in I}$ is a \v Cech
  $(q-1)$-cocycle with respect to the cover $\cvr V
  \cap X_c$ representing the Dolbeault class of $w_b$ (each
  $\alpha_{b;\:\idx 1.q}$ is holomorphic and is 
  globally $L^2$ with respect to $\vphi_F$ and $\rs\omega$ on
  $\lcS+1<c>[b] \cap V_{\idx 1.q}$).
  (Write $\alpha_{b;\:\idx 1.q} =: \frac{\rs*\gamma_{b;\:\idx
      1.q}}{\sect_{(b)}}$ just to maintain consistency with Proposition
  \ref{prop:res-formula-dbar-exact-dot-harmonic}.)
  Take $\gamma_{\idx 1.q} \in \logKX. \otimes
  \aidlof+1<X>*\paren{V_{\idx 1.q}}$ such that
  \begin{equation*}
    \PRes[\lcS+1[b]](\frac{\gamma_{\idx 1.q}}{\sect_D})
    = \frac{\rs\gamma_{b;\:\idx 1.q}}{\sect_{(b)}} \; ,
    \;\;\text{ and then set }\quad
    \frac{\rs\gamma_{p;\:\idx 1.q}}{\sect_{(p)}}
    := \PRes[\lcS](\frac{\gamma_{\idx 1.q}}{\sect_D}) \; .
  \end{equation*}
  We see that the \v Cech cocycle $\set{\frac{\paren{\delta
        \rs*\gamma_p}_{\idx 0.q}}{\sect_{(p)}}}_{\idx 0,q \in I}$
  represents the component of the image $\delta w$ on $\lcS<c>$
  \newtext{(with values in $K_{\lcS<c>} \otimes \res F_{\lcS<c>}
    \otimes \mtidlof<\lcS<c>>{\vphi_F}$)}, where
  $\delta\rs\gamma_p$ denotes the image of the cochain
  $\set{\rs\gamma_{p;\:\idx 1.q}}_{\idx 1,q \in I}$ under the \v Cech
  coboundary operator.
  \newtext{\label{page:explain_loc-L2}(One can also see that $\set{\frac{\paren{\delta
          \rs*\gamma_p}_{\idx 0.q}}{\sect_{(p)}}}_{\idx 0,q \in I}$
    takes values in $K_{\lcS<c>} \otimes \res F_{\lcS<c>} \otimes
    \mtidlof<\lcS<c>>{\vphi_F}$, thus holomorphic in particular, by
    noticing that $(\delta\gamma)_{\idx 1.q} \in \ker
    \Res^{\sigma+1}$, i.e.~$(\delta\gamma)_{\idx 0.q} \in
    \logKX. \otimes \aidlof<X>* (V_{\idx 0.q})$, as
    $\set{\frac{\rs*\gamma_{b;\:\idx 1.q}}{\sect\ps_b}}_{\idx 1,q \in
      I}$ is a cocycle for all $b\in\Iset+1$.)}
  % In the sense of the $L^2_\tloc$ Dolbeault isomorphism
  In view of the \v Cech--Dolbeault map
  \eqref{eq:Cech-Dolbeault-isom}, the image $\delta
  w$ is then represented by $\paren{-\frac{\dbar
      v_{p;(\infty)}}{\sect_{(p)}}}_{p \in \Iset}$, where
  \begin{equation*}
    - \frac{\dbar v_{p;(\infty)}}{\sect_{(p)}}
    := -\frac{
      \dbar\paren{
        \paren{\dbar\rho}^{\idx q.1}\rs\gamma_{p;\:\idx 1.q}
      }
    }{\sect_{(p)}}
    = (-1)^q \paren{\dbar\rho}^{\idx q.0}
    \frac{\paren{\delta \rs*\gamma_p}_{\idx 0.q}}{\sect_{(p)}}
    \quad\text{ on } \lcS<c> \; ,
  \end{equation*}\mmark{}{The sign problem is fixed here.}%
  (see also \cite{Chan&Choi&Matsumura_injectivity}*{Proof of
    Thm.~3.4.1, Step IV}).
  % Abuse $\delta w$ to mean this representative in what follows.
  
  Note that \newtext{each $\frac{
      \paren{\delta\rs*\gamma_p}_{\idx 0.q}
    }{\sect\ps}$ has coefficients lying in $\mtidlof<\lcS<c>>{\vphi_F}$ on
    its domain and} $-\frac{\dbar v_{p;(\infty)}}{\sect_{(p)}}$ is
  therefore locally, but not necessarily globally, $L^2$ (with respect
  to $\vphi_F$ and $\rs\omega$) on $\lcS<c>$. 
  However, if $w$ is in the image of the restriction map
  \begin{equation*}
    \jmath^{c'}_c \colon
    \spHarm/q-1/{\residlof+1*}(c') \to \spHarm/q-1/{\residlof+1*}(c)
  \end{equation*}
  given in \eqref{eq:restriction-maps-between-Takegoshi-harm-sp} for
  some $c' > c$, then the corresponding $-\frac{\dbar
    v_{p;(\infty)}}{\sect_{(p)}}$ is locally $L^2$ in $\lcS<c'> = \lcS
  \cap X_{c'}$, thus globally $L^2$ in $\lcS<c> = \lcS \cap X_c$.
  In view of this observation, define
  \begin{equation} \label{eq:Gamma_c-L2-delta-image-subsp}
    \Gamma_c := \smashoperator[r]{\bigcup_{c' \in (c,\infty]}}
    \jmath^{c'}_c \spHarm/q-1/{\residlof+1*}(c')
    \;\;
    \subset \spHarm/q-1/{\residlof+1*}(c) \; .
  \end{equation}
  Notice that each $\Harm<\lcS>(c)$ forms a closed subspace in the
  Hilbert space $\Ltwo/0,q/<\lcS<c>>{\logKX<\lcS>}_{\vphilist}$.
  For every $w \in \Gamma_c$, \mmark{write $\deltaH w$ to mean the projection
    of $\paren{(-1)^q \frac{\dbar v_{p;(\infty)}}{\sect_{(p)}}}_{p \in
      \Iset}$, a representative of $(-1)^{q-1}\delta w$, to
    $\spHarm{\residlof*}(c)$.}{
    $\deltaH$ is redefined here.
  }\footnote{
    The image $\deltaH w$ of $\deltaH$ is well-defined because a different
    choice of the \v Cech representative of $w$ or a different choice of
    the lifting $\set{\gamma_{\idx 1.q}}_{\idx 1,q \in I}$ results in the
    representative of $\delta w$ being altered by an element in $\paren{\im\dbar}_\tloc
    \subset \bigoplus_p \Ltwo./0,q/<\lcS<\alert{c'}>>{\cdot}$ for
    some $c' > c$, hence in $\paren{\im\dbar}_{(2)} \subset \bigoplus_p
    \Ltwo/0,q/<\lcS<\alert{c}>>{\cdot}$ under the restriction map, which is
    orthogonal to $\spHarm{\residlof*}(c)$.
  } Proposition \ref{prop:res-formula-dbar-exact-dot-harmonic} can now
  be translated to the following adjoint relation between $\deltaH$ and the
  harmonic residue map $\HRes$.
  Note that the factor $(-1)^{q-1}$ in the definition of $\deltaH$ is
  there just to express this relation more neatly.

  \begin{thm} \label{thm:HRes-duality}
    For any $c \in (0, \infty)$ (excluding the case $c = \infty$), 
    any $w \in \Gamma_c \subset \spHarm/q-1/{\residlof+1*}(c)$ and
    any $u \in \spHarm{\residlof*}(c)$, we have 
    \begin{equation*}
      \iinner{\deltaH w}{u}_{\lcc' \cap X_c^\circ}
      =\sigma_+ \iinner{w}{\HRes(u)}_{\lcc+1' \cap X_c^\circ} \; , 
    \end{equation*}
    where the inner products are given by the residue norms on their
    respective sets of lc centers with respect to $\vphi_F$ and
    $\rs\omega$, and $\sigma_+ := \max\set{1, \sigma}$.
  \end{thm}

  \begin{proof}
    Write $\HRes(u) =: \paren{\HRes(u)_b}_{b\in\Iset+1}$ and abuse
    $\delta w$ to mean its representative $\paren{-\frac{\dbar
        v_{p;(\infty)}}{\sect_{(p)}}}_{p \in \Iset}$.
    The conclusion of Proposition
    \ref{prop:res-formula-dbar-exact-dot-harmonic} can then be written as
    \begin{align*}
      \iinner{
      \alert{\eta_{c,\nu}} (-1)^{q-1}\delta w
      }{u}_{\lcc' \cap X^\circ}
      &=
        \sum_{p\in\Iset} \iinner{\alert{\eta_{c,\nu}} (-1)^q\frac{\dbar v_{p;(\infty)}}{\sect_{(p)}} \:}{
        \: u_p
        }_{\lcS*}
      \\
      &\overset{\mathclap{\text{Prop.~\ref{prop:res-formula-dbar-exact-dot-harmonic}}}}= \quad\;
        \sigma_+ \smashoperator[l]{\sum_{b\in\Iset+1}}
        \iinner{\alert{\eta_{c,\nu}} (-1)^{q-1}\frac{v_{b;(\infty)}}{\sect_{(b)}}
        \:}{\:
        \HRes(u)_b
        }_{\lcS*+1[b]}
      \\
      &\overset{\mathclap{\text{\eqref{eq:Cech-Dolbeault-on-w_b}}}}= \quad\;
        \sigma_+ \smashoperator[l]{\sum_{b\in\Iset+1}}
        \iinner{\alert{\eta_{c,\nu}} \paren{w_b -\dbar v_{b;(2)}}
        \:}{\:
        \HRes(u)_b
        }_{\lcS*+1[b]}
      \\
      &\overset{\mathclap{\text{Thm.~\ref{thm:Takegoshi-harmonicity-of-HRes}}}}{
        \underset{\mathclap{\alert[Gray]{\idxup{\diff\Phi}. \HRes(u) = 0}}}=} \qquad
        \sigma_+ \smashoperator[l]{\sum_{b\in\Iset+1}}
        \iinner{\alert{\eta_{c,\nu}} w_b -\dbar \paren{\alert{\eta_{c,\nu}} v_{b;(2)}}
        \:}{\:
        \HRes(u)_b
        }_{\lcS*+1[b]}
      \\
      &\overset{\mathclap{\text{Thm.~\ref{thm:Takegoshi-harmonicity-of-HRes}}}}{
        \underset{\mathclap{\alert[Gray]{\HRes(u) \text{ harmonic}}}}=} \qquad
        \sigma_+ 
        \iinner{\alert{\eta_{c,\nu}} w}{\HRes(u)}_{\lcc+1' \cap X^\circ}
    \end{align*}\mmark{}{The sign problem is fixed here.}%
    for all $\nu \in \Nnum$.
    % It can be seen that all involved entities, in particular
    \newtext{The forms} $\delta
    w$ with $w \in \Gamma_c$ \newtext{and $\HRes(u)$ with $u \in
      \spHarm{\residlof*}(c)$} are $L^2$ (with respect to $\vphi_F$
    and $\rs\omega$) on $\lcc|\bullet|' \cap X_c^\circ$ for $\bullet =
    \sigma$ or $\sigma+1$ \newtext{(by the paragraph right before
      this theorem and by Theorem \ref{thm:Takegoshi-harmonicity-of-HRes})}, so
    it is legitimate to take the limit $\nu \to +\infty$ on both
    sides.
    The desired result follows by noting that $\iinner{(-1)^{q-1}\delta
      w}{u}_{\lcc' \cap X_c^\circ} =\iinner{\deltaH w}{u}_{\lcc' \cap X_c^\circ}$
    since $u \in \spHarm{\residlof*}(c)$.
  \end{proof}

  \begin{remark}
    When $X$ is compact, we have to consider only the case $c =
    \infty$ (as $X = X_\infty$).
    Since the representatives of $\delta w$ for any $w \in
    \spHarm/q-1/{\residlof+1*}(\infty)$ are $L^2$ on $\lcc'$, we can
    choose $\Gamma_\infty := \spHarm/q-1/{\residlof+1*}(\infty) \isom
    \spH/q-1/{\residlof+1*}$ and the adjoint relation still holds true.
  \end{remark}
  
}

\section{Proofs of main results}
\label{sec:main_results_proofs}

\subsection{Injectivity on $X$ with $D$ being prime}
\label{sec:proof-on-X}

In this section, we prove a special case of Theorem \ref{thm:main-log-smooth} under the assumption that $X$ is a K\"ahler manifold and $D$ is a prime divisor. 
This will help illustrate the main ideas of the proof in the general case. Specifically, we prove the following theorem.
\begin{thm}
\label{T:prime_divisor_case}
	Let $(X, D)$ be a log smooth pair where $D$ is a prime divisor and 
	let $\pi \colon X \to \Delta$ be a proper locally K\"ahler morphism to 
	an % (not necessarily irreducible or reduced)
        analytic space $\Delta$. 
	Let $F$ (resp.~$M$) be a line bundle on $X$ 
	with a smooth potential $\vphi_F$ (resp.~$\vphi_M$)
        % Hermitian metric $h_{F}$  (resp.~$h_{M}$)
	such that 
	\begin{equation*}
		\ibddbar\vphi_F \geq 0 \quad  \text{ and } \quad
		% \sqrt{-1}(\Theta_{h_F}(F)-t \Theta 
		% _{h_M}(M))\geq 0
		% -C\omega \leq
		% -C\omega \leq
                \ibddbar\vphi_M \leq C \ibddbar\vphi_F
		\quad \text{ for some } C>0 \; . 
	\end{equation*}
	Consider a section $s \in H^{0}(X, M)$ such that $s^{-1}(0)\nsupseteq D$.
	Then, the multiplication map  induced by the tensor product with $s$ between the higher direct image sheaves
	\begin{equation*}
		R^{q}\pi_{*}\paren{K_{X} \otimes D \otimes F}
		\xrightarrow{\otimes s} 
		R^{q}\pi_{*}\paren{K_{X} \otimes D \otimes F \otimes M}
	\end{equation*}
	is injective for every $q \geq 0$.
\end{thm}

\begin{proof}
	It suffices to prove that for a fixed $t\in\Delta$, the germ
        $\beta_t\in R^q\pi_*(K_{X}\otimes D\otimes F)_t$ vanishes if
        $s\beta_t = 0$ in $R^q\pi_*(K_{X}\otimes D\otimes F\otimes M)_t$.
	The proof consists of the following steps. 
	
	\begin{step}[Reduction step] \label{item:pf-reduction-simple}
		The standard exact sequence $0 \to K_{X} \to K_{X}\otimes D \to K_{D} \to 0$ induces the following diagram.
		\begin{equation}\label{h}
			\begin{aligned}
				\xymatrix@C=3.5em@R=3.5ex{
					\ar[d] & \ar[d]\\
					{\smash R^q\pi_*(K_X\otimes F)}
					\ar[d]^-{}\ar[r] ^-{\otimes s}
					\ar[d]^-{\tau}
					&{R^q\pi_*(K_X\otimes F\otimes M )}
					\ar[d]
					\\
					{R^q\pi_*(K_{X}\otimes D\otimes F)} \ar[d]^-{}\ar[r] ^-{\otimes s}
					&{R^q\pi_*(K_{X}\otimes D \otimes F\otimes M)} \ar[d]
					\\
					{R^q\pi_*(K_{D}\otimes F\vert_D)}
					\ar[d]\ar[r]^-{\otimes s|_D }
					&{R^q\pi_*(K_{D}\otimes(F \otimes M)\vert_D)} \ar[d]\\
					& }
			\end{aligned}
		\end{equation}
		The assumption on $s^{-1}(0)$ and the curvature assumption is still satisfied after restricting $F$ and $M$ to $D$.
		Hence, by \cite{Matsumura_injectivity-Kaehler}, the map $\otimes s\vert_D$ on the bottom row is injective. 
		An easy diagram chasing implies that there exists an element $\alpha_t\in R^q\pi_*(K_X\otimes F)_t$ such that $\tau_t(\alpha_t)=\beta_t$.
		Since the problem is local on $\Delta$, by shrinking
                $\Delta$ to a relatively compact Stein open set, we may assume that
		\begin{itemize}
			%\item $\alpha$ is defined up to the boundary $\partial\Delta$,
			\item $\varphi_F$ and $\varphi_M$ are smooth up to the boundary $\partial X$,
			\item the function $\abs{s}_{\varphi_M}$ is globally bounded on $X$.
		\end{itemize}
		As $X$ is holomorphically convex after the shrinking, we may further assume that there exists a representative $\alpha\in R^q\pi_*(K_X\otimes F)(\Delta)\isom H^q(X,K_X\otimes F)$ of $\alpha_t\in R^q\pi_*(K_X\otimes F)_t$ satisfying
        \begin{itemize}
            \item $0=s\tau(\alpha)\in R^q\pi_*(K_X\otimes D\otimes F\otimes M)(\Delta)$, and
            \item \mmark{$\alpha$ is the restriction of a class on a
                neighborhood of the compact $\cl X$ (the closure of $X$).}{$\alpha$ is not yet interpreted as a form on $X$,
              and even its representatives are $L^2_{\tloc}$, so
              nothing is smooth up to $\bdry X$, hence the change.}
        \end{itemize}
        In what follows, we will show that \mmark{$\tau\paren{\res\alpha_{X_c}}=0$ in
        $H^q(X_c,K_X\otimes D \otimes F)$}{Not to show $\res\alpha_{X_c} = 0$ here.} for some $c\in(0,\infty)$, which implies the desired conclusion.    
	\end{step}

    \begin{step}[{Takegoshi harmonic representative of $\alpha$ and
        its orthogonal projection $u^\perp$ to
        $\paren{\deltaH\Gamma_c}^\perp$}]
      \label{item:pf-Take-repr-simple}
		
      %   	In this step, we associate $\alpha$ with a suitable Takegoshi harmonic representative $u^\perp$ via
      %   	\begin{equation*}
      % \mathcal H^{q}_{X,\varphi_F,\omega}(c)
      % :=
      % \mathcal H^{n,q}_{\varphi_F,\omega}(c)
      %   		\xhookrightarrow{\jmath^{c}}
      %   		H^q(X_c,K_X\otimes F)
      %   	\end{equation*}
      %   as given in Sections \ref{sec:L2-theory} and \ref{sec:residue-element}.
    
    	The $L^2_{\tloc}$ Dolbeault isomorphism (see Section \ref{sec:L2-theory}) asserts that $\alpha$ can be represented by a $\dbar$-closed locally $L^2$ form on $X$.
		Since $\alpha$ is defined across the boundary $\bdry X$, the class $\alpha$ can be represented by a \emph{globally} $L^2$ form  (with respect to $\varphi_F$) on $X$.
		Fix a complete metric $\omega$ on $X$ as described in Section \ref{subsec:notation}.
		After taking a projection to the harmonic space (see
    \eqref{decom}) and taking into account Theorem
    \ref{thm:Takegoshi-argument} and Remark~\ref{rem-T-property}, we
    see that $\alpha$ is represented by a Takegoshi harmonic form
    $u\in\mathcal H^q_{X,\varphi_F,\omega}(\infty) := \Harm<X>{\logKX
      ; \Phi}$, 
    i.e.~$u$ is a harmonic form with respect to $\varphi_F,\omega$ on
    $X$ satisfying $\idxup{\diff\Phi}. u=0$ on $X$ (or see \eqref{T-space} and
    \eqref{eq:def-Takagoshi-harmonic-sp-on-lcS} for the definition).
		Notice that, for any $c \in (0, \infty]$, the restriction $\res u_{X_c}$ is still a Takegoshi harmonic form in $\mathcal H^q_{X,\varphi_F,\omega}(c)$ such that $\jmath^c\paren{\res u_{X_c}} =\res\alpha_{X_c}$ (see~\eqref{eq:restriction-maps-between-Takegoshi-harm-sp}). 

	  For any fixed $c\in(0,\infty)$, consider the subspace
    % if the subspace $\Gamma_c$ of $\mathcal H^{q-1}_{D,\varphi_F,\omega}$ is defined by
    \begin{equation} \label{eq:Gamma_c-simple}
      \Gamma_c := \smashoperator[r]{\bigcup_{c' \in (c,\infty]}}
      \jmath^{c'}_c \mathcal H^{q-1}_{D,\varphi_F,\omega}(c')
      \;\;
      \subset \newtext{\mathcal H^{q-1}_{D,\varphi_F,\omega}(c)}
    \end{equation}
    (cf.~\eqref{eq:Gamma_c-L2-delta-image-subsp}) and the map 
		\begin{equation*}
			\deltaH \colon \Gamma_c
			\to
			\mathcal H^q_{X,\varphi_F,\omega}(c)
    \end{equation*}
    induced from the connecting morphism $\delta$ in the long exact sequence 
		\begin{equation*}
			\xymatrix@R=0.3cm@C=1.5em{
				{\cdots\to \cohgp{q-1}[D_c]{K_D \otimes F}} \ar[r]^-{\delta}
				&
				\cohgp q[X_c]{K_X\otimes F} \ar[r]^-{\tau}  
				&
				\cohgp q[X_c]{K_X\otimes D \otimes F}  \to\cdots \; ,
			} 
		\end{equation*}
    where $D_c := D \cap X_c$, as discussed in Section \ref{sec:residue-element}.
    % which satisfies the duality in the sense of Theorem~\ref{thm:HRes-duality}.
		Then $u$ can be orthogonally decomposed as
		\begin{equation*}
			u|_{X_c} = u^\perp+\mu \in 
          \paren{\deltaH\Gamma_c}^\perp
          \oplus 
          \cl{\deltaH\Gamma_c}
        	=
        	\mathcal H^q_{X,\varphi_F,\omega}(c) \; .
		\end{equation*}
		From the long exact sequence above, we see that $\mu \in \cl{\deltaH\Gamma_c} \subset
		\ker\tau$ (the map $\jmath^c$ is made implicit).
		Therefore, to complete the proof, it suffices to show that
		$u^\perp = 0$ on $X_c$.
	\end{step}

	\begin{step}[$\HRes(u^\perp)$ as an obstruction of
      $\norm{su^\perp}_{X_c}^2 = 0$]\label{item:expression-of-su-simple}

		We make use of the assumption $\eqcls{s\tau(u)} =\eqcls{s\tau\paren{u^\perp}} =0$ in $H^q(X,K_X
    \otimes D \otimes F\otimes M)$ and the \v Cech--Dolbeault map in
    Section \ref{sec:Dolbeault} to re-express $\norm{s u^\perp}_{X_c}^2$ as follows.
			
		\begin{itemize}
			\item Given the finite Stein open cover $\cvr V = \set{V_i}_{i\in I}$ and the partition of unity $\set{\rho^i}_{i\in I}$ subordinate to $\cvr V$ given in Section \ref{subsec:notation}, it follows from the discussion in Section~\ref{sec:Dolbeault} that there exist a \v Cech cocycle $\{\alpha^\perp_{\idx 0.q}\}_{\idx 0.q \in I}$ representing the class of $u^\perp$ via the $L^2_{\tloc}$ Dolbeault isomorphism on $X_c$ and a globally $L^2$ section $v_{(2)}$ of $K_{X}\otimes F$ on $X_c$ with respect to $\norm\cdot_{X_c} := \norm\cdot_{X_c, \vphilist}$ such that
			\begin{equation*} \tag{$*$} \label{E:CechDolbeault_isom}
				u^\perp
				=
				\:\dbar v_{(2)}
				+(-1)^q \:\underbrace{\dbar \rho^{i_{q}} \wedge \dotsm \wedge
					\dbar\rho^{i_1} \cdot \rho^{i_0} }_{=: \:
					\paren{\dbar\rho}^{\idx q.0}} \alpha^\perp_{\idx 0.q} \; .
			\end{equation*}
      Note that each $\alpha^\perp_{\idx 0.q} \in H^{0}(V_{\idx 0.q} \cap X_c , K_{X}\otimes F)$ is globally $L^2$ on the open set $V_{\idx 0.q} \cap X_c$.

			\item The fact $\eqcls{s\tau\paren{u^\perp}}=0$ in $\cohgp
        q[X_c]{K_X\otimes D \otimes F \otimes M}$ guarantees the
        existence of $\lambda_{\idx 1.{q}} \in H^{0}(V_{\idx 1.{q}}
        \cap X_c , K_X \otimes D\otimes F \otimes M)$ for $\idx 1,{q}
        \in I$ such that (note that $\tau$ is given by $\otimes \sect_D$)
			\begin{equation*}
				s\alpha^\perp_{\idx 0.q}\sect_D = (\delta\lambda)_{\idx 0.q}
				\quad\text{ on } V_{\idx 0.q} \cap X_c \; ,
			\end{equation*}
			where $\delta$ is \newtext{the} \v Cech coboundary operator.
			Thanks to this identity, the second term on the right-hand side in
      \eqref{E:CechDolbeault_isom}, after applying $\otimes s
      \sect_D$, can be expressed as \label{page:explain_cobdry}
			\begin{alignat*}{1}
				(-1)^q
				\paren{\dbar\rho}^{\idx q.0}
				s\alpha^\perp_{\idx 0.q}\sect_D
				&=
				(-1)^q
				\dbar\rho^{i_{q}} \wedge
				\dotsm\wedge
				\dbar\rho^{i_1} \cdot \rho^{i_0}
				\paren{\delta\lambda}_{\idx 0.q}
        \\
        &\newtext{%
          {}=(-1)^q
          \dbar\rho^{i_{q}} \wedge
          \dotsm\wedge
          \dbar\rho^{i_1} \cdot \rho^{i_0}
          \sum_{\nu=0}^q (-1)^\nu \lambda_{\idx 0.[\widehat{i_\nu} \dotsm]q}
        }
				\\
				&=
				(-1)^q\dbar\rho^{i_{q}} \wedge
				\dotsm\wedge
				\dbar\rho^{i_1}
				\cdot\lambda_{\idx 1.q}
        \qquad\newtext{\alert[Gray]{(\because \sum_{i_\nu} \dbar\rho^{i_\nu} \equiv 0)}}
				\\
				&=
				-
				\dbar
				\paren{
					\dbar \rho^{i_{q}} \wedge
					\dotsm\wedge
					\dbar\rho^{i_2}\rho^{i_1}
					\cdot\lambda_{\idx 1.q}
				}
				\\
				&=
				-
				\dbar\paren{\paren{\dbar\rho}^{\idx q.1}\lambda_{\idx1,q}}
				=:-\dbar v_{(\infty)} \; ,
			\end{alignat*}
			where $v_{(\infty)}$ is a smooth section of $K_X\otimes D\otimes
      F\otimes M$, which is globally $L^2$ with respect to
      $\sm\vphi_D +\vphi_F +\vphi_M$ (not to $\phi_D$) on $X_c$.
      Thus, we have 
      \begin{equation*}
        su^\perp=\dbar \paren{s v_{(2)}} - \frac{\dbar v_{(\infty)}}{\sect_D}.
      \end{equation*}
    \end{itemize}
		
		Since $u^\perp$ is harmonic with respect to $\varphi_F$ on $X_c$ and
		we also have $\ibddbar\vphi_F \geq 0$ and $%-C\omega \leq
		\ibddbar\vphi_M \leq C\ibddbar\vphi_F$ on $X$ for some constant
		$C > 0$ by assumption,
		Proposition~\ref{prop:consequence-of-positivity} and Lemma~\ref{lem:su-harmonicity} guarantee that  
		$su^\perp$ and $su^\perp\sect_D$ are harmonic with respect to $\vphi_F+\vphi_M$ and $\phi_D+\vphi_F+\vphi_M$, respectively, on $X_c$.
		Then $\norm{su^\perp}^2_{X_c}$ can be computed 
		using the Takegoshi property and the residue computation in Proposition~\ref{prop:residue-formula-classical-kernel} 
		\newtext{(the same argument as in the proof of
      Proposition~\ref{prop:res-formula-dbar-exact-dot-harmonic},
      repeated here for clarifying ideas in this simple case)}, which
    gives \label{page:explain_limit-order2}
		\begin{align*}
			\norm{su^\perp}_{X_c}^2
			=
			&~\iinner{
				\dbar \paren{s v\ps_{2}} -\frac{\dbar v\ps_\infty
				}{\sect_D}}{su^\perp}_{X_c}
       \quad
       \overset{\alert[Gray]{\hphantom{\im\dbar}\mathllap{\paren{\im\dbar}_{(2)}} \perp s \newtext{u^\perp}}}
			=
			-\iinner{
				\frac{\dbar v\ps_\infty}{\sect_D}
			}{su^\perp}_{X_c}
			\\
			\xleftarrow{\nu \tendsto +\infty}
			&-
			\iinner{\eta_{c,\nu}\frac{\dbar v\ps_\infty}{\sect_D}}
			{su^\perp}_{X_c}
            =
            -
            \iinner{\eta_{c,\nu}\dbar v\ps_\infty}
			{su^\perp\sect_D}_{X_c,\phi_D}
            \\ 
            \newtext{\xleftarrow{\eps\tendsto0+}}
			      &\newtext{~-
            \iinner{e^{-\eps\abs{\psi_D}}\eta_{c,\nu} \dbar v\ps_\infty}
			{su^\perp\sect_D}_{X_c,\phi_D}}
            \\
            \newtext{=}&
            \newtext{~-
            \iinner{e^{-\eps\abs{\psi_D}} \dbar(\eta_{c,\nu} v\ps_\infty)}
			{su^\perp\sect_D}_{X_c,\phi_D}
            +
            \iinner{e^{-\eps\abs{\psi_D}} v\ps_\infty}
			      {\smash{\cancelto{
            \mathllap{(\because \text{ Takegoshi property})\quad}
            0 
            }{\idxup{\partial\eta_{c,\nu}}.su^\perp\sect_D}}}_{X_c,\phi_D}
            % \quad\;\;\; \vphantom{\sum_{p\in \Iset}}
            }
            \\
            =&
            -
            \cancelto{0 \;\;\;(\because~u^\perp \text{ harmonic, Lemma~\ref{lem:su-harmonicity} or Proposition \ref{prop:Takegoshi-harmonic-forms}})}
            {\iinner{\dbar\paren{e^{-\eps\abs{\psi_D}}\eta_{c,\nu} v\ps_\infty}}
            {su^\perp\sect_D}_{\mathrlap{X_c,\phi_D}}} \quad\;\;\; \vphantom{\sum_{p\in \Iset}}
            +
            \eps\iinner{e^{-\eps\abs{\psi_D}}\eta_{c,\nu} v\ps_\infty}
			{\idxup{\partial\psi_D}.su^\perp\sect_D}_{X_c,\phi_D}
            \\
            =&
            \sum_{\idx 1.q\in I}
            \eps\iinner{e^{-\eps\abs{\psi_D}}\eta_{c,\nu}\lambda_{\idx1.q}}
			      {\idxup{\partial\rho},[\idx 1.q].\idxup{\partial\psi_D}.su^\perp\sect_D}_{X_c,\phi_D}
            \\
            \xrightarrow[
            \substack{
            \text{Prop.~\ref{prop:residue-formula-classical-kernel}}}
            ]{\eps \tendsto 0^+}
			&
            \sum_{\idx 1.q\in I}
            \iinner{\eta_{c,\nu}
            \mathcal R_D\paren{\frac{\lambda_{\idx1.q}}{\sect_D}}}
			{\idxup{\partial\rho},[\idx 1.q].
            \mathcal R_D\paren{\idxup{\partial\psi_D}.su^\perp}}_{D_c}
            \\
            =&
            \;\;\iinner{\eta_{c,\nu}
               \:v_{b;(\infty)}}{s\:\HRes(u^\perp)}_{D_c} \; ,
        \end{align*}
        where
		\begin{equation*}
			v_{b;(\infty)}
			=\sum_{\idx1,q\in I}
			\paren{\dbar\rho}^{\idx q.1}
			\rs*\lambda_{b;\idx 1.q}
			\quad\text{ with }\quad
			\rs*\lambda_{b;\idx 1.q}
			:=
			\mathcal R_D\paren{\frac{\lambda_{\idx 1.q}}{\sect_D}}
		\end{equation*}
		and
		\begin{equation*}
			\HRes(u^\perp)=\PRes[D](\idxup{\diff\psi_D}. u^\perp) \; .
		\end{equation*}
		Note that the notation here follows those used in
    Proposition~\ref{prop:res-formula-dbar-exact-dot-harmonic} (with
    $I_D^1=\set{b}$, a singleton, here) and Section~\ref{sec:residue-element}.
    It is shown below that $\HRes(u^\perp)$ is actually $0$ on $D_c$, which will then conclude the proof.
	\end{step}

    \begin{step} [{Vanishing of $\HRes(u^\perp)$ from $u^\perp \in
      \paren{\deltaH \Gamma_c}^{\perp}$ and the adjoint relation between
      $\deltaH$ and $\HRes$}]
		\label{item:pf:use_u-ortho-w}
		% In this step, we apply the assumption $u^\perp\in\paren{\deltaH \Gamma_c}^{\perp}$ and Theorem~\ref{thm:HRes-duality} to conclude the proof.
		  
		It follows from Theorem \ref{thm:Takegoshi-harmonicity-of-HRes} that    $\HRes(u^\perp) \in\mathcal H^{q-1}_{D,\varphi_F,\omega}(c)$ and from    Theorem \ref{thm:HRes-duality} that
		\begin{equation*}
			0 \overset{\alert[Gray]{\deltaH\Gamma_c \perp u^\perp}}
			=
			\iinner{\deltaH w}{u^\perp}_{X_c}
			=
			\iinner{w}{\HRes(u^\perp)}_{D_c}
			\quad\text{ for all } w \in \Gamma_c \; .
    	\end{equation*}
		It thus suffices to show that $\HRes(u^\perp) \in \cl{\Gamma_c}$ (the closure of $\Gamma_c$ in $\mathcal H^{q-1}_{D,\varphi_F,\omega}(c)$),    which will imply that $\HRes(u^\perp) = 0$.
		
		Recall that $u^\perp = \res u_{X_c} - \mu$.
		We have $\HRes(u) \in \mathcal H^{q-1}_{D,\varphi_F,\omega}(\infty)$ by    Theorem \ref{thm:Takegoshi-harmonicity-of-HRes}, so $\HRes(\res    u_{X_c}) = \res{\HRes(u)}_{X_c} \in \Gamma_c$ by the definition    \eqref{eq:Gamma_c-L2-delta-image-subsp} of $\Gamma_c$.
		Furthermore, by $\mu \in \cl{\deltaH\Gamma_c}$, there is a sequence $\seq{w_\nu}_{\nu \in \Nnum}$ such that
		\begin{equation*}
    		w_\nu \in \mathcal H^{q-1}_{D,\varphi_F,\omega}(c_\nu)
			\;\;\text{ for some } c_\nu > c
			\quad\text{ and }\quad
			\deltaH w_\nu
			\xrightarrow{\nu \tendsto \infty} \mu
			\;\;\text{ in }
			\mathcal H^q_{X,\varphi_F,\omega}(c)\; .
		\end{equation*}
    	Recall from the discussion in Section \ref{sec:residue-element} that, although the representatives of $\delta w_\nu$ may not be globally $L^2$ on $X_{\alert{c_\nu}}$, they are globally $L^2$ on the smaller space $X_{\alert{c_\nu^-}}$ for any $c_\nu^- \in(c, c_\nu)$.
		This means that we have $\deltaH w_\nu \in\mathcal H^q_{X,\varphi_F,\omega}(c_\nu^-)$.
		(All the maps $\jmath^c$ and $\jmath^{c'}_c$ are made implicit here for clarity.) 
		Theorem \ref{thm:Takegoshi-harmonicity-of-HRes} again guarantees that 
		$\HRes(\deltaH w_\nu) \in\mathcal H^{q-1}_{D,\varphi_F,\omega}(c_\nu^-)    \subset \Gamma_c$.
		Since $\HRes$ is a bounded linear operator, this clearly implies that $\HRes(\mu) \in \cl{\Gamma_c}$.
		As a result, we obtain $\HRes(u^\perp) = \res{\HRes(u)}_{X_c} -\HRes(\mu)\in \cl{\Gamma_c}$, which completes the proof. \qedhere
	\end{step}

\end{proof}

\subsection{Injectivity for $(X,D)$ and $(Y,0)$ in general}
\label{sec:proof-on-D}

In this section, we prove a special case of Theorem \ref{main-thm} on
$Y$ (with $D_Y=0$) and Theorem \ref{thm:main-log-smooth} on $X$ with a
generalization such that the potentials $\vphi_F$ and $\vphi_M$ are
allowed to be singular.
These will be used to prove the full version of Theorem \ref{main-thm}
in Section \ref{sec:reduction-to-log-smooth}.

\begin{thm} % [Theorem \ref{thm:main-log-smooth} with singular metrics]
  \label{thm:main-thm-in-section}
  \mmark{}{
    Since other details related to $Y$ are added, and Section
    \ref{sec:intro} seems far away from this section, I keep the full
    assumptions here. When I read a paper, I feel better if I don't
    have to flip through all the pages to check the assumptions for
    a technical statement.
  }%
    \mmark{}{\xb{Thank you. This version is better.}
   }

  Let $(X, D)$ be a K\"ahler log smooth lc pair and $\pi \colon X \to
  \Delta$ be a proper locally K\"ahler morphism to an analytic space
  $\Delta$.
  Let $F$ (resp.~$M$) be a line bundle on $X$ 
  equipped with a potential $\vphi_F$  (resp.~$\vphi_M$) such
  that it has the analytic singularities described as in Section
  \ref{subsec:notation} and satisfies
  \begin{equation*}
    \ibddbar\vphi_F \geq 0 \quad  \text{ and } \quad
    % \sqrt{-1}(\Theta_{h_F}(F)-t \Theta 
    % _{h_M}(M))\geq 0
    % -C\omega \leq
    -C\omega \leq \ibddbar\vphi_M \leq C \ibddbar\vphi_F
    \quad \text{ for some } C>0 \; . 
  \end{equation*}
  \newtext*{Also let $Y$ be a reduced snc divisor satisfying the same
  assumptions as $D$ (in particular, $Y$ has only no common components
  with the polar divisors $P_F$ and $P_M$ and $Y+P_F+P_M$ is an snc
  divisor) and $\defidlof{Y}$ be its defining ideal sheaf on $X$.}
  Consider a section $s \in H^{0}(X, M)$ such that its zero locus
  $s^{-1}(0)$ contains no lc centers of the pairs $(X,D)$ and $(X,Y)$, and
  that the function $\abs s_{\vphi_M}$ is locally bounded on
  $X$.
  Write $\pi_Y$, $F_Y$, $M_Y$ and $s_Y$ as the pullback of $\pi$, $F$,
  $M$ and $s$ to $Y$ respectively.
  Then, the multiplication maps induced by $\otimes s$ and
  $\otimes s_Y$ between the higher direct image sheaves
  \begin{gather*}
    R^{q}\pi_{*}\paren{K_{X} \otimes D \otimes F \otimes \mtidlof{\vphi_F}}
    \xrightarrow{\otimes s} 
    R^{q}\pi_{*}(K_{X} \otimes D \otimes F \otimes M
    \otimes \mtidlof{\vphi_F+\vphi_M})
    \\
    \drR q[{\pi_Y}](K_Y \otimes F_Y \otimes
      \frac{\mtidlof{\vphi_F}}{\mtidlof{\vphi_F} \cdot \defidlof{Y}})
    \xrightarrow{\otimes s_Y}
    \drR q[\pi_Y](K_Y \otimes F_Y \otimes M_Y \otimes \frac{\mtidlof{\vphi_F
        +\vphi_M}}{\mtidlof{\vphi_F +\vphi_M} \cdot \defidlof{Y}})
  \end{gather*}
  are injective for every $q \geq 0$, where $K_Y := \parres{K_X
    \otimes Y}_Y = K_X \otimes Y \otimes \frac{\holo_X}{\defidlof{Y}}$.
\end{thm}

\begin{remark} \label{remark:holo_Y-structure-on-quotient}
  \newtext*{For $\vphi_L := \vphi_F$ or $\vphi_F +\vphi_M$, the quotient
  $\frac{\mtidlof{\vphi_L}}{\mtidlof{\vphi_L} \cdot \defidlof{Y}}$ has
  the structure of an $\holo_Y$-sheaf (as $\defidlof{Y}$ is in its
  annihilator).
  It follows that $K_X \otimes Y \otimes L \otimes_{\holo_X}
  \frac{\mtidlof{\vphi_L}}{\mtidlof{\vphi_L} \cdot \defidlof{Y}} = K_Y
  \otimes L_Y \otimes_{\holo_Y} \frac{\mtidlof{\vphi_L}}{\mtidlof{\vphi_L} \cdot
    \defidlof{Y}}$.
  By the snc assumptions on $\vphi_L$ and on $Y$, we have
  $\mtidlof{\vphi_L} \cdot \defidlof Y = \mtidlof{\vphi_L+\psi_Y}
  =\aidlof|0|{\vphi_L}[\psi_Y]$ and thus
  $\frac{\mtidlof{\vphi_L}}{\mtidlof{\vphi_L} \cdot \defidlof{Y}} =
  \faidlof|n|/|0|{\vphi_L}[\psi_Y]$.}
  % is a locally free $\holo_X$-sheaf,
\end{remark}

\begin{proof}
  \setDefaultMetric{\rs*\omega}

  The outline of the proof is the same as the one in
  \cite{Chan&Choi&Matsumura_injectivity}*{\S 3.4} with some
  adjustments.

  \begin{enumerate}[label=\textbf{Step \arabic*}, ref=\arabic*,
    leftmargin=0pt, labelsep=*, widest=I,
    itemindent=*, align=left,
    itemsep=1.5ex]

    \newcommand{\tpln}[1]{(#1)\textbf{.}} %% topic line

  \item
    \tpln{Reduction via induction on lc centers and
      the local nature of the problem}
    Set
    \begin{equation*}
      \aidlof|-1| := \aidlof|-1|{\vphi_F+\vphi_M} := 0 
      \quad\text{ and }\quad
      \residlof|0| := D^{-1} \otimes \mtidlof{\vphi_F} \; .
    \end{equation*}
    % By the diagram chaing argument using the commutative
    % diagram~\eqref{eq:commut-diagram_sing-Fujino-conj}, the problem can
    % be reduced to the injectivity of the multiplication map
    % \begin{equation}
    %   \otimes s:H^q(D,K_D\otimes F)\rightarrow H^q(D,K_D\otimes F\otimes M).
    % \end{equation}
    % For the detailed arguemnt, see \cite{CCM}*{Corollary 1.3}.  \medskip
    % 
    % The proof is slight variant of the one of Theorem~1.2 in \cite{CCM}.
    % As in the proof of \cite{CCM}, we 
    Write
    \begin{align*}
      \aidlof* &:= \aidlof =\mtidlof{\vphi_F} \cdot \defidlof{\lcc+1'} \; , 
      \\
      \residlof* &:= \residlof \isom \faidlof/-1* \; ,
      \\
      \aidlof*^M &:= \aidlof{\vphi_F+\vphi_M} \quad \text{ and }
      \\
      \spR{\sheaf F} &:=\drR q(\logKX. \otimes
        \sheaf F) \quad\text{ for any sheaf } \sheaf F \text{ on } X
                       \; . 
    \end{align*}
    % \mmark{}{I keep using \texttt{enumerate} instead of the
    %   \texttt{step} environment which is a properer
    %   environment here. It also provides better control of various
    %   lengths and distances. The ``(...)'' description is not a good
    %   style in this situation, yet I'm trying my best to mimic it.}%
    Recall that
    % \begin{equation*}
    %   K_D = K_X \otimes D \otimes \faidlof|\sigma_{\mlc}|/|0|* \; ,
    % \end{equation*}
    % and
    the inclusions between adjoint ideal sheaves induce the short
    exact sequences
    \begin{equation} \label{eq-ex2}
      \xymatrix{
        0 \ar[r]
        & {\faidlof-1/|\sigma_0-1|*} \ar[r]
        & {\faidlof|n|/|\sigma_0-1|*} \ar[r]
        & {\faidlof|n|/-1*} \ar[r]
        & 0
      } \quad\text{ for } 0 \leq \sigma_0 \leq \sigma \leq n \; .
    \end{equation}
    % \mmark{}{\xb{I want to change the subsections $\sigma_0, \sigma, n$
    %     to $a, b, c$...} Sorry that I can't change the indices to
    %   $a,b,c$, as they really correspond to the indices in the
    %   following commutative diagram. I've adjusted a bit so that the
    %   correspondence looks more apparent.}% 
    They, together with the multiplication map $\otimes s$, induce the commutative diagrams
    %%%%%%%%%%%%%%
    
% \begin{equation} \label{eq:commut-diagram_sing-Fujino-conj}
%   \begin{gathered}
%     \xymatrix@R=0.85cm@C+0.75cm{
%       {\vdots} \ar[d]
%       & {\vdots} \ar[d]
%       & {\vdots} \ar[d]
%       \\
%       {\spH{\faidlof-1/|\rho|*}} \ar[d] \ar@{=}[r]
%       & {\spH{\faidlof-1/|\rho|*}} \ar[r]^-{\otimes s}
%       \ar[d]_-{\iota_{\sigma-1}} \ar[dr]|-*+{\mu_{\sigma-1}}
%       & {\spH M{\faidlof-1/|\rho|*M}} \ar[d]
%       \\
%       {\spH{\faidlof/|\rho|*}} \ar[d] \ar[r]^-{\iota_{\sigma}}
%       \ar@/_1.99pc/[rr]|(.65)*+{\mu_{\sigma}}
%       & {\spH{\faidlof|\sigma_{\mlc}|/|\rho|*}}
%       \ar[d]|(.38)*+<3pt>{ }
%       \ar[r]^-{\otimes s}
%       & {\spH M{\faidlof|\sigma_{\mlc}|/|\rho|*M}}
%       \ar[d]
%       \\
%       {\spH{\residlof*}} \ar[d] \ar[r]^-{\tau_\sigma}
%       \ar@/_1.68pc/[rr]+<-39pt,-15pt>|(.67)*+{\nu_\sigma}
%       & {\spH{\faidlof|\sigma_{\mlc}|/-1*}} \ar[d]|(.52)*+<3pt>{}
%       \ar[r]^-{\otimes s}
%       & {\spH M{\faidlof|\sigma_{\mlc}|/-1*M}} \ar[d] \\
%       {\vdots} & {\vdots} & {\vdots} }
%   \end{gathered}
% \end{equation}
\begin{equation} \label{eq:commut-diagram_sing-Fujino-conj}
  \begin{gathered}
    \xymatrix@R=0.85cm@C+0.75cm{
      {\vdots} \ar[d]
      & {\vdots} \ar[d]
      & {\vdots} \ar[d]
      \\
      {\spR{\faidlof-1/|\sigma_0-1|*}} \ar[d] \ar@{=}[r]
      & {\spR{\faidlof-1/|\sigma_0-1|*}} \ar[r]^-{\otimes s}
      \ar[d]_-{\iota_{\sigma-1}} \ar[dr]|-*+{\mu_{\sigma-1}}
      & {\spR M{\faidlof-1/|\sigma_0-1|*^M}} \ar[d]
      \\
      {\spR{\faidlof/|\sigma_0-1|*}} \ar[d] \ar[r]^-{\iota_{\sigma}}
      \ar@/_1.95pc/[rr]|(.65)*+{\mu_{\sigma}}
      & {\spR{\faidlof|n|/|\sigma_0-1|*}}
      \ar[d]|(.38)*+<3pt>{ }
      \ar[r]^-{\otimes s}
      & {\spR M{\faidlof|n|/|\sigma_0-1|*^M}}
      \ar[d]
      \\
      {\spR{\residlof*}} \ar[d] \ar[r]^-{\tau_\sigma}
      \ar@/_1.68pc/[rr]+<-39pt,-15pt>|(.67)*+{\nu_\sigma}
      & {\spR{\faidlof|n|/-1*}} \ar[d]|(.52)*+<3pt>{}
      \ar[r]^-{\otimes s}
      & {\spR M{\faidlof|n|/-1*^M}} \ar[d] \\
      {\vdots} & {\vdots} & {\vdots} }
  \end{gathered}
\end{equation}
%
    %%%%%%%%%%%%%%
    for $\sigma =\sigma_0, \sigma_0 +1,\dots,n$.
    Here the columns are exact, and $\mu_\sigma$ (resp.~$\nu_\sigma$)
    is the composition of $\iota_\sigma$ (resp.~$\tau_\sigma$) with the
    map induced from $\otimes s$.
    Note that $\mu_{\sigma_0} =\nu_{\sigma_0}$ and $\iota_{\sigma_0}
    =\tau_{\sigma_0}$.

    Since $\newtext{\iota_{n}} = \id$ (the identity map) at $t$
    % for any
    % % $\sigma$ and
    % \newtext{fixed} $\sigma'$
    % % such that $\sigma' \geq \sigma \geq
    % % \sigma_{\mlc}$
    (hence $\newtext{\mu_{n}}
    % = \mu_{\sigma_{\mlc}}
    = \otimes s$ at $t$), both injectivity statements in the claim at $t$ are
    proved if we put in the diagram \eqref{eq:commut-diagram_sing-Fujino-conj} 
    \begin{itemize}
    \item  $\sigma_0 := 0$ (for the map $\otimes s$)
    \item  $\sigma_0 := 1$ and $D := Y$ (for the map
      $\otimes s_Y$, note that $\faidlof/|0|*$ has an $\holo_Y$-sheaf
      structure as explained in Remark \ref{remark:holo_Y-structure-on-quotient})
    \end{itemize}
    and show that $\paren{\ker\mu_{n}}_t = \paren{\ker\iota_{n}}_t \;\paren{= 0}$.
    % ($\iota_{n}$ is the identity map at $t$).
    Following the argument in
    \cite{Chan&Choi_injectivity-I}*{Thm.~1.3.2}, since
    $\paren{\ker\mu_{\sigma-1}}_t =\paren{\ker\iota_{\sigma-1}}_t$ and
    $\paren{\ker\nu_\sigma}_t =\paren{\ker\tau_\sigma}_t$ together
    imply $\paren{\ker\mu_{\sigma}}_t =\paren{\ker\iota_{\sigma}}_t$
    via a diagram-chasing argument, to prove the injectivity of the
    map in the claim, it suffices to show that
    \begin{equation*} % \label{eq:ker-nu=ker-tau}
      \paren{\ker\nu_\sigma}_t =\paren{\ker\tau_\sigma}_t \quad\text{ for all }
      \sigma = \newtext{\sigma_0, \dots, n}
      % \begin{cases}
      %   0, 1, \dots, n & (\text{for } \otimes s) \\
      %   1, \dots, n & (\text{for } \otimes s_Y) \\
      % \end{cases}
      \:\text{ and for all } t \in \Delta \; .
    \end{equation*}
    Note that we obviously have $\ker\tau_\sigma \subset
    \ker\nu_\sigma$ for all $\sigma \geq 0$ (on $\Delta$).
    The remainder of the proof is devoted to proving the reverse
    inclusions for all $t \in \Delta$.
    \newtext*{
      \begin{remark} \label{remark:mlc-index}
        For any $t \in \Delta$ and any neighborhood $U_t \Subset \Delta$
        of $t$, there exists the smallest positive integer $\sigma_{\mlc}
        \leq n$ (\newtext{possibly} dependent on $U_t$\newtext{, but only decreases as $U_t$ shrinks}) such that
        \begin{equation*}
          \aidlof|\sigma_{\mlc} -1|* \subsetneq \aidlof|\sigma_{\mlc}|* =
          \aidlof* = \mtidlof{\vphi_F}
          \;\;\text{ on $\pi^{-1}(U_t)$ }\quad
          \text{ for all integers $\sigma \geq \sigma_{\mlc}$}
        \end{equation*}
        by the strong Noetherian property of increasing sequences of
        coherent sheaves.\footnote{
          Indeed $\sigma_{\mlc}$ is the codimension of the \newtext{minimal log-canonical centers} of
          $(\pi^{-1}(U_t), D)$, as $\aidlof|\sigma_{\mlc} -1|* =
          \mtidlof{\vphi_F} \cdot \defidlof{\lcc|\sigma_{\mlc}|'} \neq
          \mtidlof{\vphi_F}$. 
        }
        \newtext{Set it to be the one given on $\pi^{-1}(\Delta) =
          X$.}
        Then we see that the proof till this point is the same as in
        \cite{Chan&Choi&Matsumura_injectivity}*{\S 3.4}, where 
        the index $n$ is replaced by $\sigma_{\mlc}$ everywhere.
      \end{remark}
    }
    \newtext{%
      \begin{remark} \label{remark:general_quotients}
        \newtext*{Notice that, by replacing the index $n$ by some integer
        $\sigma' \in [\sigma_0, n]$ in this proof,}
        % when the proof is done,
        it is actually proved that the injectivity in
        the theorem (for $\otimes s$) holds true even when the
        multiplier ideal sheaves $\mtidlof{\vphi_F}$ and
        $\mtidlof{\vphi_F+\vphi_M}$ are replaced by
        $\faidlof|\sigma'|/|\sigma_0-1|*$ and
        $\faidlof|\sigma'|/|\sigma_0-1|*^M$ for any choices of
        $\sigma_0$ and $\sigma'$ such that $0\leq \sigma_0 \leq
        \sigma' \;\;(\leq n)$.
      \end{remark}
    }
    At this point, fix any point $t \in \Delta$ and any integer
    $\sigma = % (0, )\, 1, \dots, n
    \newtext{\sigma_0, \dots, n}$.
    % (set $\sigma_{\mlc} :=\sigma_{\mlc}^t$ for convenience).
    Pick any germ
    \begin{equation*}
      \alpha_t \in \paren{\ker\nu_\sigma}_t \subset \spR{\residlof*}_t
      \; .
    \end{equation*}
    It suffices to prove that $\alpha_t \in \paren{\ker\tau_\sigma}_t$.

    As the problem is local on $\Delta$, we can shrink $\Delta$ to a
    (sufficiently small) relatively compact Stein open neighborhood
    of $t \in \Delta$ such that $\alpha_t$ is lifted to a section
    $\alpha \in \ker\nu_\sigma \subset \spR{\residlof*}$ on $\Delta$
    which is defined even across the boundary $\bdry\Delta$.
    \newtext*{(This may reduce $\sigma_{\mlc}$, but does not affect the
      current proof.)}
    The manifold $X$ is shrunk accordingly and becomes a
    \emph{holomorphically convex manifold}.
    % and there exists a smooth psh exhaustion function $\Phi
    % \geq 0$ on $X$ with \mmark{$\sup_X \Phi=\infty$}{
    %   \xb{memo: I want to assume that $\sup_X \Phi=\infty$ as in
    %     \ref{subsec:notation}}
    %   Sure! I mistakenly thought that exhaustion functions must satisfy
    %   $\sup_X \Phi =\infty$ and thus missed out this condition.
    % } which is lifted from the
    % one on the Stein space $\Delta$,
    % such that $\sup_X
    % \abs{d\Phi}_\omega < \infty$ can also be chosen accordingly.
    % \mmark{
    \mmark{}{
      I keep the more extended description of the assumptions after
      shrinking $\Delta$ instead
      of referring them to Sec.~\ref{subsec:notation} as this is
      the place where we explain why the assumptions in
      Sec.~\ref{subsec:notation} make sense.
    }%
    Thanks to such shrinking, we can therefore assume that 
    \begin{itemize}
    \item the potentials $\vphi_F$ and $\vphi_M$ are
      smooth on their regular loci across the boundary $\bdry X$,
      
    \item the function $\abs s_{\vphi_M}$ is globally bounded on $X$,
      
      % \item the differential form $d\Phi$ is smooth up to $\bdry X$ and
      %   thus $\abs{d\Phi}_\omega$ is bounded on $X$ (where $\omega$ is
      %   a \emph{complete} K\"ahler metric on $X$), 
      %   \xb{memo: I want to assume that $\sup_X \Phi=\infty$ as in \ref{subsec:notation}} and
      
    \item $X$ admits a \emph{finite} Stein open covering $\cvr V :=
      \set{V_i}_{i \in I}$
    \end{itemize}
    as in the assumptions stated in Section \ref{subsec:notation}.
    % }{
    %   Will shorten this and refer back to the corresponding step in
    %   the special case if that remains.
    % }
    Furthermore, construct the smooth exhaustion psh function $\Phi$ and
    the \emph{complete} K\"ahler metric on $X$ described in Section
    \ref{subsec:notation}.
    For any coherent sheaf $\sheaf F$ on $X_c := \set{\Phi < c}$ for
    $c \in (0,\infty]$, we have the isomorphism 
    \begin{equation*}
      \spH{\sheaf F}_c :=\cohgp q[X_c]{\logKX. \otimes \sheaf F}
      \isom
      \spR{\sheaf F}\paren{\Delta_c} \; ,
      % \quad\text{ for any sheaf } \sheaf F \text{ on } X_c :=
      % \set{\Phi < c}
    \end{equation*}
    thanks to the Leray spectral sequence and Cartan's Theorem B on
    the Stein space $\Delta_c := \set{\Phi_\Delta < c}$. 
    % and , from the isomorphism $\cohgp q[X]{\sheaf F}$  and for the Stein space $\Delta$, we have
    Therefore, the section
    \begin{equation*}
      \alpha \in \ker\nu_\sigma \subset \spR{\residlof*}(\Delta) \isom
      \spH{\residlof*}_\infty \; ,
    \end{equation*}
    is abused to mean the corresponding class in $\spH{\residlof*}_\infty$ on $X$
    (supported on $\lcc'$) which is \emph{the restriction of a class
      on a neighborhood of $\cl X$} (the closure of $X$).
    Abusing also $\ker\nu_\sigma$ and $\ker\tau_\sigma$ to mean their
    corresponding subspaces in $\spH{\residlof*}_\infty$, we are going
    to prove that $\res\alpha_{X_c} \in \res{\ker\tau_\sigma}_{X_c}
    \subset \spH{\residlof*}_c$ for some $c > 0$, which will complete
    the proof.

    % we are left to prove the injectivity of the
    % multiplication map 
    % \begin{equation} \label{eq:injectivity_D}
    %   H^q(X,K_X\otimes D\otimes F \otimes \mtidlof{\vphi_F})
    %   \xrightarrow{\:\otimes s\:}
    %   H^q(X,K_X\otimes D\otimes F\otimes M \otimes \mtidlof{\vphi_F +\vphi_M})
    % \end{equation}
    % for all $q \geq 0$.

    % \renewcommand{\objectstyle}{\displaystyle} As in the proof of
    % Theorem~\ref{T:prime_divisor_case}, the problem can be reduced
    % to the injectivity of the multiplication map
    % \begin{equation*}
    %   \otimes s:H^q(X,K_X\otimes D\otimes F) \rightarrow
    %     % 	 H^q(X,K_X\otimes D\otimes F\otimes M)
    % \end{equation*}
    % under the assumption $\pi\colon X\rightarrow\Delta$ is a surjective
    % holomorphic submersion to a subvariety $\Delta$ in $\bC^N$
    % together with conditions (\romannumeral1),(\romannumeral2),
    % (\romannumeral3), and (\romannumeral4) in Step 1.

  \item \label{item:pf:step-2}
    \tpln{Takegoshi harmonic representative of $\res\alpha_{X_c}$ and its
      orthogonal projection $u^\perp$ to $\paren{\deltaH\Gamma_c}^\perp$}
    We first show that $\res\alpha_{X_c}$ can be associated to a Takegoshi
    harmonic form via the map
    $\spHarm{\residlof*}(c)
    % := \bigoplus_{p \in \Iset} \Harm<\lcS>(\infty)
    \xhookrightarrow{\jmath^{c}} \spH{\residlof*}_c$ given in
    Sections \ref{sec:L2-theory} and \ref{sec:residue-element}.
    % 
    % For a cohomology class $\alpha\in\ker\nu_\sigma$, we will prove that $\tau_\sigma(\alpha)$ vanishes.
    % 
    Write $\lcc' =\bigcup_{p \in \Iset} \lcS$ as the union of
    $\sigma$-lc centers $\lcS$ of $(X,D)$ and let $\rs\omega$ be the
    complete K\"ahler metric on $X^\circ := X \setminus \paren{P_F
      \cup P_M}$ as described in Section \ref{subsec:notation}.
    The \newtext{map} $\jmath^\infty$ in
    \eqref{eq:map-for-harmonic-representatives-on-X_c}
    % \eqref{eq:map-for-harmonic-representatives}
    induces the monomorphism
    \begin{equation*}
      \bigoplus_{p \in \Iset} \Harm<\lcS>{\logKX<\lcS>}
      \xhookrightarrow{\:\jmath^\infty\:}
      \bigoplus_{p \in \Iset} \cohgp q[\lcS]{K_{\lcS}
        \otimes F \otimes \mtidlof<\lcS>{\vphi_F}}
      =
      \spH{\residlof*}_\infty \; .
    \end{equation*}
    The $L^2_{\tloc}$ Dolbeault isomorphism (see Section \ref{sec:L2-theory}) 
    asserts that the component of $\alpha$ on each
    $\sigma$-lc center $\lcS$ can be represented 
    by a $\dbar$-closed locally $L^2$
    form with respect to $\vphi_F$ on $\lcS$.
    Since $\alpha$ is defined across the boundary $\bdry X$, 
    the component of $\alpha$ on $\lcS$ can be represented by a
    \emph{globally} $L^2$ form on $\lcS$.
    After taking a projection to the harmonic space (see \eqref{decom}) and taking into account
    Theorem \ref{thm:Takegoshi-argument} and Remark
    \ref{rem-T-property}, we see that $\alpha$ is represented by 
    \begin{align*}
      u :=\paren{u_p}_{p \in \Iset}
      &\in \bigoplus_{p\in \Iset} \Harm<\lcS>{\logKX<\lcS>}
      \\
      &\overset{\mathclap{\text{Thm.~\ref{thm:Takegoshi-argument}}}}=
        \quad\;
        \bigoplus_{p\in \Iset} \Harm<\lcS>(\infty)
        =: \spHarm{\residlof*}(\infty) \; ,
    \end{align*}
    where $u =\paren{u_p}_{p \in \Iset}$ is a collection of Takegoshi harmonic forms (see \eqref{T-space} for the definition),  
    i.e.~each $u_p$ is a harmonic form with respect to $\vphilist$
    % on $\lcS$
    satisfying $\idxup{\diff\Phi}. u_p = 0$ on $\lcS* := \lcS
    \cap X^\circ$.
    Notice that, for any $c \in (0, \infty]$, the restriction $\res
    u_{X_c}$ is still a collection of Takegoshi harmonic
    forms in $\spHarm{\residlof*}(c)$ such that
    $\jmath^c\paren{\res u_{X_c}} =\res\alpha_{X_c}$ (see
    \eqref{eq:restriction-maps-between-Takegoshi-harm-sp}).

    % such that each summand are isomorphic to the locally $L^2$ cohomology,
    % \begin{equation*}
    %   \cohgp q[\lcS]{K_{\lcS} \otimes F}
    %   \xrightarrow{\quad \cong \quad }
    %   \dfrac{\ker\dbar}{\im\dbar} 
    %   \text{ of } \Ltwoloc^{n-1,q}\paren{\lcS,F}_{\varphi_F},
    % \end{equation*}		
    % by \cite{Matsumura_injectivity-Kaehler}*{Proposition 2.16} as before.
    % Hence $\alpha$ can be written as
    % \begin{equation*}
    %   \alpha=\paren{\alpha_p}_{p\in\Iset}=:\paren{[u_p]}_{p\in\Iset}=:[u]
    % \end{equation*}
    % for some $u_p\in\ker\dbar\in\Ltwoloc^{0,q}\paren{K_{\lcS}\otimes\lcS,F}_{\varphi_F}$.
    Now fix any positive $c < \infty$.
    Recall from \eqref{eq:Gamma_c-L2-delta-image-subsp} the subspace
    $\Gamma_c$ of $\spHarm/q-1/{\residlof+1*}(c)$ and the map
    \begin{equation*}
      \deltaH \colon \Gamma_c \to \spHarm{\residlof*}(c)
    \end{equation*}
    induced from the connecting homomorphism $\delta \colon
    \spH/q-1/{\residlof+1*}_c \to \spH{\residlof*}_c$, as discussed in
    Section \ref{sec:residue-element}.
    Recall also that we take the (squared) residue norm
    $\norm\cdot_{\lcc' \cap X_c^\circ}^2 = \sum_{p \in \Iset}
    \norm\cdot_{\lcS*<c>, \vphilist}^2$ as the $L^2$ norm on $\spHarm{\residlof*}(c)$ 
    (see Section \ref{sec:adjoint-ideal-n-residue}).
    Then, we decompose $u$ orthogonally as 
    \begin{equation*}
      \res u_{X_c} = u^\perp + \mu
      \:\in\: \paren{\deltaH\Gamma_c}^\perp \oplus \cl{\deltaH\Gamma_c}
      = \spHarm{\residlof*}(c) .
    \end{equation*}
    \mmark{}{Order of terms on the RHS changed, as I just feel that we
      care more the entities on the left; think
    about the Taylor series...}%
    %Let $u^\perp$ be the orthogonal projection of $\res u_{X_c}$ to
    %$\paren{\deltaH\Gamma_c}^\perp$ in $\spHarm{\residlof*}(c)$ such
    %that
    %\begin{equation*}
     % \mu := \res u_{X_c} -u^\perp \in \cl{\deltaH\Gamma_c} \subset
     % \spHarm{\residlof*}(c) \; .
    %\end{equation*}
    The short exact sequences in \eqref{eq-ex2} induce the commutative diagram
    \begin{equation} \label{eq:commut-diagram_delta-tau}
      \begin{aligned}
        \xymatrix@R-0.3cm{
          {\dotsm} \ar[r]
          & {\spH/q-1/{\residlof+1*}_c} \ar[r]^-{\delta}
          \ar[d]^-{\tau_{\sigma+1}}
          & {\spH{\residlof*}_c} \ar[r]
          \ar@{=}[d]
          & {\spH{\faidlof+1/-1*}_c} \ar[r]
          \ar[d]
          & {\dotsm}
          \\
          {\dotsm} \ar[r]
          & {\spH/q-1/{\faidlof|\newtext{n}|*}_c} \ar[r]
          & {\spH{\residlof*}_c} \ar[r]^-{\tau_\sigma}
          & {\spH{\faidlof|\newtext{n}|/-1*}_c} \ar[r]
          & {\dotsm}
        }
      \end{aligned}
    \end{equation}
    in which the rows are exact.
    This shows that $\mu \in \cl{\deltaH\Gamma_c} \subset
    \ker\tau_{\sigma}% \subset \ker\nu_\sigma
    $ (the map $\jmath^c$ is made implicit).
    Therefore, to complete the proof, it suffices to show that
    $u^\perp = 0$ on $\lcc' \cap X_c^\circ$.

    % Using the same argument of shrinking the space $\Delta$ and
    % applying the orthogonal decompositions
    % \eqref{E:orthogonal_decompostion1} and
    % \eqref{E:orthogonal_decompostion2}, it suffices to consider the
    % problem under the following conditions. 
    
    % \begin{itemize}
    % \item each $u_p$ is a harmonic form on $\lcS$ with respect to the
    %   given norm $\norm\cdot_{\lcS}$ and
    
    % \item $u \in \ker\nu_\sigma \cap \paren{\ker\tau_\sigma}^\perp$,
    %   where the orthogonal complement $\paren{\ker\tau_\sigma}^\perp$
    %   of $\ker\tau_\sigma$ is taken with respect to the residue norm
    %   $\norm\cdot_{\lcc'}$.
    % \end{itemize}
    % The theorem is proved if it is shown that $u_p = 0$ for all
    % $p \in \Iset$.

  \item \label{item:express-su-in-residue-norm}
    \tpln{$\HRes(u^\perp)$ as an obstruction of
      $\norm{su^\perp}_{\lcc' \cap X_c^\circ}^2 = 0$}
    We show below that the harmonic residue $\HRes(u^\perp)$
    is an obstruction to our desired vanishing $u^\perp = 0$ on $\lcc'
    \cap X_c^\circ$ by rewriting $\norm{su^\perp}_{\lcc' \cap X_c^\circ}^2$ using the
    assumption $u^\perp \in \ker\nu_\sigma$ and the \v Cech--Dolbeault
    map \eqref{eq:Cech-Dolbeault-isom}.
    Note that $u^\perp \in \ker\nu_\sigma$ as both $\res u_{X_c}$ and $\mu$ belong
    to $\ker\nu_\sigma$ on $X_c$.
    
    Recall that $\cvr V :=\set{V_i}_{i \in I}$ is the finite Stein cover of $X$
    % by admissible open sets with respect to
    % $(\vphi_F,\vphi_M,\psi_D)$ 
    and $\set{\rho^i}_{i \in I}$ is the
    partition of unity subordinate to $\cvr V$ as
    described in Section \ref{subsec:notation}.
    % Their notations are abused to mean also their induced cover and
    % partition of unity on $\lcc'$ for any $\sigma \geq 0$.
    % For any choice of indices $\idx 0,q \in I$, write $V_{\idx 0.q}
    % :=V_{i_0} \cap V_{i_1} \dotsm \cap V_{i_q}$ as usual.
    % 
    Through the \v Cech--Dolbeault map
    \eqref{eq:Cech-Dolbeault-isom}, the cohomology 
    class of the component $u^\perp_p$ of $u^\perp$ on
    each $\lcS<c>$ is represented by a \v Cech $q$-cocycle
    $\set{\alpha^\perp_{p; \:\idx 0.q}}_{\idx 0,q \in I}$
    % Under the Einstein summation convention on the indices $\idx 0,q$,
    % we see that
    such that
    \begin{equation} % \tag{$*$}
      \label{eq:pf:Cech-Dolbeault-u_p}
      u^\perp_p
      % \overset{\text{\eqref{eq:Cech-Dolbeault-isom}}}
      =
      \:\dbar v_{p;(2)}
      +(-1)^q \:\underbrace{\dbar \rho^{i_{q}} \wedge \dotsm \wedge
        \dbar\rho^{i_1} \cdot \rho^{i_0} }_{=: \:
        \paren{\dbar\rho}^{\idx q.0}} \alpha^\perp_{p; \:\idx 0.q} \; ,
    \end{equation}
    where $v_{p; (2)}$ is a $K_{\lcS} \otimes \res{F}_{\lcS}$-valued $(0,q-1)$-form
    on $\lcS*<c>$ with $L^2$ coefficients with respect to
    $\norm\cdot_{\lcS*<c>}$ and
    $\alpha^\perp_{p; \:\idx 0.q} \in K_{\lcS} \otimes \res F_{\lcS} \otimes
    \mtidlof<\lcS>{\vphi_F}$ on
    $\lcS<c> \cap V_{\idx 0.q} :=\lcS<c> \cap V_{i_0} \cap \dotsm \cap V_{i_q}$.
    By the residue exact sequence \eqref{eq:residue-exact-seq}, 
    for each choice of the multi-indices $(\idx 0,q)$, there exists a section
    $f_{\idx 0,q} \in \logKX. \otimes \aidlof*$ on the Stein open set
    $V_{\idx 0.q} \cap X_c$ such that
    \begin{equation*}
      \Res^\sigma(f_{\idx 0.q})
      =\paren{\alpha^\perp_{p; \:\idx 0.q}}_{p \in \Iset} \; .
    \end{equation*}
    Considering the inclusion $\aidlof* \subset
    \aidlof|\newtext{n}|*$ and multiplying the cochain by the
    section $s$, we see that the image $\nu_\sigma\paren{u^\perp}$ is
    represented by the \v Cech $q$-cocycle $\set{\eqcls{s f_{\idx
          0.q}}}_{\idx 0,q \in I}$, in which
    \begin{equation*}
      \eqcls{s f_{\idx 0.q}} := \paren{s f_{\idx 0.q} \bmod \aidlof-1*^M}
      \;\in \logKX. M \otimes\faidlof|\newtext{n}|/-1*^M
      \quad\text{ on } V_{\idx 0.q} \cap X_c \; .
    \end{equation*}
    % the collection is then representing $\nu_\sigma(u)$ in $\spH M{\faidlof|n|/-1*}$.
    %(Note that $\set{\eqcls{sf_{\idx 0.q}}}_{\idx 0,q \in I}$ is a
    %cocycle since $\set{\alpha^\perp_{p;\:\idx 0.q}}_{\idx 0,q \in I}$
    %is for each $p \in \Iset$.)
    The assumption $u^\perp \in \ker\nu_\sigma$ implies that this cocycle is
    a coboundary, that is, there exists $\lambda_{\idx 1.q} \in \logKX. M
    \otimes\aidlof|\newtext{n}|*^M$ on $V_{\idx 1.q} \cap X_c$ for each $(\idx 1,q)$ such that 
    \begin{equation*}
      \set{\eqcls{sf_{\idx 0.q}}}_{\idx 0,q \in I}
      =\delta\set{\eqcls{\lambda_{\idx 1.q}}}_{\idx 1,q \in I}
      =\set{\eqcls{\paren{\delta\lambda}_{\idx 0.q} }}_{\idx 0,q \in
        I} \; ,
    \end{equation*}
    where $\paren{\delta\lambda}_{\idx 0.q}$ is given by the usual formula of \v Cech coboundary operator $\paren{\delta\lambda}_{\idx 0.q}:=\sum_{k =0}^q (-1)^k \lambda_{\idx 0[\dotsm \widehat{i_k}].q}$.
    Note that $\lambda_{\idx 1.q}$ need \emph{not} take values in $\aidlof*^M$.    
    Since $sf_{\idx 0.q}$ and $\paren{\delta\lambda}_{\idx 0.q}$
    differ by an element in $\logKX. M \otimes \aidlof-1*^M$ on
    $V_{\idx0.q} \cap X_c$, we see that 
    \begin{equation*}
      \newtext{\paren{\PRes[\lcS](
          \frac{\paren{\delta\lambda}_{\idx 0.q}}{\sect_D}
          )}_{\mathrlap{p \in \Iset}}} \quad
      =\Res^\sigma\paren{\paren{\delta\lambda}_{\idx 0.q}}
      =\Res^\sigma\paren{s f_{\idx 0.q}}
      =\paren{s \alpha^\perp_{p;\:\idx 0.q}}_{p \in \Iset} \; .
    \end{equation*}
    
    Following the notation in Proposition
    \ref{prop:res-formula-dbar-exact-dot-harmonic},
    % \newtext{and
    % recalling that $\Res^\sigma\paren{\paren{\delta\lambda}_{\idx
    % 0.q}} = \paren{\PRes[\lcS](\frac{\paren{\delta\lambda}_{\idx
    % 0.q}}{\sect_D})}_{p \in \Iset}$},
    set
    \begin{alignat*}{2}
      \rs\lambda_{p;\:\idx 1.q}
      &:= \PRes[\lcS](\frac{\lambda_{\idx 1.q}}{\sect_D}) \cdot \sect_{(p)}
      &\quad&\text{ for each } p \in \Iset \;\text{ and}
      \\
      \rs\lambda_{b;\:\idx 1.q}
      &:= \PRes[\lcS+1[b]](\frac{\lambda_{\idx 1.q}}{\sect_D}) \cdot \sect_{(b)}
      &&\text{ for each } b \in \Iset+1 \; .
    \end{alignat*}
    We then have $\frac{\paren{\delta\rs*\lambda_p}_{\idx
        0.q}}{\sect_{(p)}} = s \alpha^\perp_{p;\:\idx 0.q}$ on $\lcS<c>$ for
    each $p\in \Iset$.
    Together with \eqref{eq:pf:Cech-Dolbeault-u_p}, we obtain that
    \begin{align*}
      s u^\perp_p - \dbar\paren{s v_{p;(2)}}
      &= (-1)^q \paren{\dbar\rho}^{\idx q.0}
        \frac{\paren{\delta\rs*\lambda_p}_{\idx 0.q}}{\sect_{(p)}}
      \\
      &=(-1)^q\:
        \dbar\rho^{i_q} \wedge \dotsm \wedge \dbar\rho^{i_1}
        \cdot\frac{\rs*\lambda_{p;\:\idx 1.q}}{\sect_{(p)}}
      \\
      &=-\frac{
        \dbar\paren{\paren{\dbar\rho}^{\idx q.1}
          \rs*\lambda_{p;\:\idx 1.q}}
        }{\sect_{(p)}}
        =: -\frac{\dbar v_{p;(\infty)}}{\sect_{(p)}} \; .
    \end{align*}
    Setting also $v_{b;(\infty)} := \sum_{\idx 1,q \in I}
    \paren{\dbar\rho}^{\idx q.1} \rs*\lambda_{b; \:\idx 1.q}$ on
    $\lcS+1<c>[b]$, Proposition
    \ref{prop:res-formula-dbar-exact-dot-harmonic} then yields
    \begin{align*}
      \norm{su^\perp}_{\lcc' \cap X_c^\circ}^2
      \qquad\;\;
      &\overset{\mathclap{\text{Prop.~\ref{prop:consequence-of-positivity}}}}{
      \underset{\mathclap{\alert[Gray]{\hphantom{\im\dbar}\mathllap{\paren{\im\dbar}_{(2)}} \perp s u_p^\perp}}}=}
      \quad\;\;
      \sum_{p\in\Iset}
      \iinner{s u_p^\perp - \dbar\paren{s v_{p;(2)}}}{ s u_p^\perp}_{\lcS*<c>}
      =~-\sum_{p\in \Iset}
         \iinner{\frac{\dbar v_{p;(\infty)}}{\sect_{(p)}}}{\: s
        u_p^\perp}_{\lcS*<c>}
      \\
      &\xleftarrow{\nu \tendsto \infty}
        ~-\sum_{p\in \Iset}
        \iinner{\alert{\eta_{c,\nu}} \frac{\dbar v_{p;(\infty)}}{\sect_{(p)}}}{\: s u_p^\perp}_{\lcS*<c>}
      \\
      &\overset{
        \mathclap{
        \text{Prop.~\ref{prop:res-formula-dbar-exact-dot-harmonic}}
        }
        }
      = \quad\;\;
      %    \sigma \smashoperator[l]{\sum_{b\in\Iset+1}}
      %    \iinner{\alert{\eta_\nu} \frac{v_{b}}{\sect_{(b)}} \:}{\quad\;\;
      %    s\;\smash{\smashoperator{
      %    \sum_{p\in\Iset \colon \lcS+1[b] \subset \lcS}
      %    }} \;\;
      %    \sgn{b:p} \:
      %    \PRes[\lcS+1[b] | \lcS](\idxup{\diff\psi_{(p)}} u_p)
      %    }_{\lcS+1[b]}
      % \\
      % =&~\quad
         \sigma_+ \smashoperator[l]{\sum_{b\in\Iset+1}}
         \iinner{\alert{\eta_{c,\nu}} \frac{v_{b;(\infty)}}{\sect_{(b)}} \:}{s
           \:\HRes(u^\perp)_b}_{\lcS*+1<c>[b]}
      \\
      &\xrightarrow{\nu \tendsto \infty}
        \lim_{\nu \tendsto \infty}
        \sigma_+ \smashoperator[l]{\sum_{b\in\Iset+1}}
        \iinner{\alert{\eta_{c,\nu}} \frac{v_{b;(\infty)}}{\sect_{(b)}} \:}
        {s \:\HRes(u^\perp)_b}_{\lcS*+1<c>[b]} \; ,
    \end{align*}%
    where $\HRes(u^\perp) = \paren{\HRes(u^\perp)_b}_{b\in \Iset+1}$
    is the harmonic residue of $u^\perp$ defined in Section
    \ref{sec:residue-element}. 
    % $v_{b;(\infty)} := \sum_{\idx 1,q \in I}
    % \paren{\dbar\rho}^{\idx q.1} \rs*\lambda_{b; \:\idx 1.q}$,
    % $\rs*\lambda_{b; \:\idx 1.q}
    % :=\PRes[\lcS+1[b]](\frac{\lambda_{\idx 1.q}}{\sect_D}) \cdot \sect_{(b)}$, and
    % \begin{equation}\label{eq-def-w}
    %   w_b := \smashoperator[r]{\sum_{p \in \Iset \colon \lcS+1[b] \subset
    %       \lcS}} \;\; \sgn{b:p}\:
    %   \PRes[\lcS+1[b] | \lcS](\idxup{\diff\psi_{(p)}}  u_{p})
    %   % \in \Harm/q-1/<\lcS*+1[b]>{\logKX<\lcS+1[b]>}
    % \end{equation}
    Therefore, if we show that $\HRes(u^\perp) = 0$ on $\lcc+1' \cap
    X_c^\circ$, 
    then $s u^\perp = 0$ and hence $u^\perp = 0$, which will conclude
    the proof.
    
  % \item Show that $w_b$ is $L^2$ harmonic with respect to
  %   $\res{\vphi_F}_{\lcS+1[b]}$ (and $\res{\omega}_{\lcS+1[b]}$) on
  %   $\lcS+1[b]$ for all $b \in \Iset+1$ and thus $\paren{w_b}_{b
  %     \in\Iset+1}$ represents a class in $\spH/q-1/{\residlof+1*}$.
    
  %   This is follows from the argument in Section~\ref{sec:residue-element}.

  \item % \label{item:pf:use_u-ortho-w}
    \tpln{Vanishing of $\HRes(u^\perp)$ from $u^\perp \in
    \paren{\deltaH \Gamma_c}^{\perp}$ and the adjoint relation between
    $\deltaH$ and $\HRes$}
    % Apply the assumption  and Theorem
    % \ref{thm:HRes-duality} to conclude the proof.
    % via the use of $w		:=\paren{w_b}_{b \in
    % \Iset+1} \in \spH/q-1/{\residlof+1*}=\bigoplus_{b
    % \in\Iset+1}
    % \cohgp{q-1}[\lcS+1[b]]{K_{\lcS+1[b]}\otimes F}$ in
    % view of the commutative diagram 
    % \begin{equation*}
    % \xymatrix@R-0.3cm{
    % {\dotsm} \ar[r]
    % & {\spH/q-1/{\residlof+1*}} \ar[r]^-{\delta}
    % \ar[d]^-{\tau_{\sigma+1}}
    % & {\spH{\residlof*}} \ar[r]
    % \ar@{=}[d]
    % & {\spH{\faidlof+1/-1*}} \ar[r]
    % \ar[d]
    % & {\dotsm}
    % \\
    % {\dotsm} \ar[r]
    % & {\spH/q-1/{\faidlof|n|*}} \ar[r]
    % & {\spH{\residlof*}} \ar[r]^-{\tau_\sigma}
    % & {\spH{\faidlof|n|/-1*}} \ar[r]
    % & {\dotsm}
    % }
    % \end{equation*}
    % and conclude that $u_p = 0$ on $\lcS$ for each $p\in\Iset$.
    % 
    It follows from Theorem \ref{thm:Takegoshi-harmonicity-of-HRes} that
    $\HRes(u^\perp) \in \spHarm/q-1/{\residlof+1*}(c)$ and from
    Theorem \ref{thm:HRes-duality} that
    \begin{equation*}
      0 \overset{\alert[Gray]{\deltaH\Gamma_c \perp u^\perp}}=
      \iinner{\deltaH w}{u^\perp}_{\lcc' \cap X_c^\circ}
      =\sigma_+ \iinner{w}{\HRes(u^\perp)}_{\lcc+1' \cap X_c^\circ}
      \quad\text{ for all } w \in \Gamma_c \; .
    \end{equation*}
    It thus suffices to show that $\HRes(u^\perp) \in \cl{\Gamma_c}$
    (the closure of $\Gamma_c$ in $\spHarm/q-1/{\residlof+1*}(c)$),
    which will imply that $\HRes(u^\perp) = 0$.

    The rest of the proof is the same as in Step
    \ref{item:pf:use_u-ortho-w} of the proof of Theorem
    \ref{T:prime_divisor_case}, mutatis mutandis (with
    $\Harm/q-1/<D>(\bullet)$ there replaced by
    $\spHarm/q-1/{\residlof+1*}(\bullet)$ and $\Harm<X>(\bullet)$ by
    $\spHarm{\residlof*}(\bullet)$).
    The proof is then completed. \qedhere
  \end{enumerate}
\end{proof}

\subsection{Injectivity for $(Y,D_Y)$}
\label{sec:reduction-to-log-smooth}

The full version of Theorem \ref{main-thm}, together with a
generalization allowing singularities on $\vphi_F$ and $\vphi_M$, is
given below.

\begin{thm}[Theorem \ref{main-thm} with singular Hermitian metrics] \label{thm:reduction-to-log-smooth}
  % Let $(Y,D_Y)$ be an snc pair globally embedded into a log smooth lc
  % pair $(X,D)$ in the sense that there is an embedding $i \colon Y
  % \hookrightarrow X$ such that both $i\paren{D_Y} =D$ and $Y$
  % (identified with its image in $X$) are both snc divisors in $X$ with
  % no common irreducible components and that the divisor $Y+D$ is also
  % snc in $X$.
  % Suppose that $\pi \colon X \to \Delta$ is a proper locally K\"ahler morphism
  % from $X$ to an analytic space $\Delta$, and $(F, e^{-\vphi_F})$ and
  % $(M, e^{-\vphi_M})$ are Hermitian line bundles on $X$ (in which
  % $\vphi_F$ and $\vphi_M$ are smooth) such that
  % \begin{equation*}
  %   \ibddbar\vphi_F \geq 0 \quad  \text{ and } \quad
  %   \ibddbar\vphi_M \leq C \ibddbar\vphi_F
  %   \quad \text{ for some } C>0 \; .
  % \end{equation*}
  % Let $s \in \cohgp 0[X]{M}$ be a holomorphic section of $M$ on $X$
  % such that $s^{-1}(0)$ contains no lc centers of $(X, Y+D)$.
  % Write $\pi_Y$, $F_Y$, $M_Y$ and $s_Y$ as the pullback of $\pi$, $F$,
  % $M$ and $s$ to $Y$ respectively.
  Using the notation and assumptions in Theorem
  \ref{thm:main-thm-in-section}, assume further that the divisor $Y$
  and $D$ has no common irreducible components, the divisor $Y+D$
  is snc and $s^{-1}(0)$ contains no lc centers of $(X,Y+D)$.
  % Assume also that $s^{-1}(0)$ contains no lc centres of $(X,Y+D)$.
  Let $D_Y := D \cap Y$.
  Then, Theorem \ref{thm:main-thm-in-section} implies that the
  multiplication map induced by $\otimes s_Y$ between
  the higher direct image sheaves
  \begin{equation*}
    R^{q}{\pi_Y}_{*}\paren{K_{Y} \otimes D_Y \otimes F_Y \otimes
    \frac{\mtidlof{\vphi_F}}{\mtidlof{\vphi_F} \cdot \defidlof{Y}}}
    \xrightarrow{\otimes s_Y} 
    \drR q[{\pi_Y}](K_{Y} \otimes D_Y \otimes F_Y \otimes M_Y
    \otimes
    \frac{\mtidlof{\vphi_F+\vphi_M}}{\mtidlof{\vphi_F+\vphi_M} \cdot \defidlof{Y}})
  \end{equation*}
  is injective for every $q \geq 0$.
  % Here, $K_Y := \parres{K_X \otimes Y}_Y$.
\end{thm}

\begin{proof}
  \newtext{The proof goes by applying
    Theorem~\ref{thm:main-thm-in-section} and its arguments to
    suitably constructed diagrams.}
  Let $\phi_Y := \log\abs{\sect_Y}^2$ be the potential on (the line
  bundle associated with) $Y$ induced
  from a canonical section $\sect_Y$.
  Define $\psi_Y := \phi_Y -\sm\vphi_Y \leq -1$, where $\sm\vphi_Y$ is some
  smooth potential on $Y$, and set $\psi_{Y+D} := \psi_Y +\psi_D$.
  
  Let $\vphi_L := \vphi_F$ or $\vphi_F +\vphi_M$. 
  Recall that the divisors $P_L$ (the polar locus of $\vphi_L$), $D$
  and $Y$ have only snc against each other, and there are no common
  components among any two of them, so we have
  \begin{align*}
    \aidlof{\vphi_L +\phi_Y}
    &=\mtidlof{\vphi_L+\phi_Y} \cdot \defidlof{\lcc+1'}
      \; ,
    \\
    \aidlof{\vphi_L}[\psi_{Y+D}]
    &=\mtidlof{\vphi_L} \cdot \defidlof{\lcc+1'(Y+D)}
      \quad\text{ and }
    \\
    \aidlof{\vphi_L}[\psi_{Y}]
    &=\mtidlof{\vphi_L} \cdot \defidlof{\lcc+1'(Y)} 
  \end{align*}
  for all $\sigma \geq 0$.
  It also follows from 
  $\mtidlof{\vphi_L +\phi_Y} \cdot \defidlof{D} = \mtidlof{\vphi_L} \cdot
  \defidlof{Y} \cdot \defidlof{D} = \mtidlof{\vphi_L} \cdot
  \defidlof{Y+D}$ that
  \begin{equation*}
    \aidlof|0|{\vphi_L+\phi_Y}
    % = \mtidlof{\vphi_L+\phi_Y+\phi_D}
    =\aidlof|0|{\vphi_L}[\psi_{Y+D}] \; .
  \end{equation*}
  Moreover,
  % for $\newtext{\sigma'} \gg 0$ (say, when $\sigma' = n$), all
  % increasing sequences of adjoint ideal sheaves above stabilize such
  % that
  \newtext*{we have}
  \begin{align*}
    \aidlof|\newtext{n}|{\vphi_L+\phi_Y} =\mtidlof{\vphi_L+\phi_Y}
    &=\mtidlof{\vphi_L}\cdot \defidlof{Y}
      \;\;\paren{= \aidlof|0|{\vphi_L}[\psi_Y]}
    \quad\text{ and}
    \\
    \aidlof|\newtext{n}|{\vphi_L}[\psi_{Y+D}] = \aidlof|\newtext{n}|{\vphi_L}[\psi_Y]
    &=\mtidlof{\vphi_L} \; .
  \end{align*}
  % There exist integers $\sigma_D, \sigma_Y, \sigma_{Y+D} \in [1, n]$
  % such that (after shrinking $\Delta$)
  % \begin{align*}
  %   \aidlof|\sigma_D-1|{\vphi_L+\phi_Y}
  %   &\subsetneq
  %     \aidlof{\vphi_L+\phi_Y}
  %     % = \mtidlof{\vphi_L+\phi_Y}
  %     =\mtidlof{\vphi_L} \cdot \defidlof{Y}
  %   &&\text{ for all } \sigma \geq \sigma_D \; ,
  %   \\
  %   \aidlof|\sigma_{Y+D}-1|{\vphi_L}[\psi_{Y+D}]
  %   &\subsetneq
  %   \aidlof{\vphi_L}[\psi_{Y+D}] = \mtidlof{\vphi_L}
  %   &&\text{ for all } \sigma \geq \sigma_{Y+D} \; ,
  %   \\
  %   \aidlof|\sigma_Y-1|{\vphi_L}[\psi_Y]
  %   &\subsetneq
  %   \aidlof{\vphi_L}[\psi_Y] = \mtidlof{\vphi_L}
  %   &&\text{ for all } \sigma \geq \sigma_Y 
  % \end{align*}
  % on $X$.
  Note also that
  \begin{equation*} \tag{$\dagger$} \label{eq:Y+D_system-n-Y_system}
    \aidlof{\vphi_L}[\psi_{Y+D}] \subset \aidlof{\vphi_L}[\psi_Y]
    \quad\text{ on $X$ for all } \sigma \geq 0 \; .%\footnotemark
  \end{equation*}%
  % \footnotetext{
    This can be seen by noticing that $\abs{\psi_{Y+D}} = \abs{\psi_Y}
    +\abs{\psi_D} \geq \abs{\psi_Y} \geq 1$ and the expression
    $e^{-\abs{\psi}} \logpole|\psi|$ is decreasing in $\abs\psi$
    ($\geq 1$) as soon as, say, $\abs\psi \geq \sigma + 1 +\eps$.
  % }%
  % This shows in particular that $\sigma_Y \leq \sigma_{Y+D}$.
  
  Write
  % \begin{alignat*}{2}
  %   \aidlof*^Y_{D} 
  %   &:=\aidlof{\vphi_F+\phi_Y}[\psi_{D}] \; , \quad
  %   &\aidlof*^{M,Y}_{D}
  %   &:=\aidlof{\vphi_F+\vphi_M+\phi_Y}[\psi_{D}] \; ,
  %   \\
  %   \aidlof*
  %   &:=\aidlof{\vphi_F}[\psi_{Y+D}] \; ,
  %   &\aidlof*^M
  %   &:=\aidlof{\vphi_F+\vphi_M}[\psi_{Y+D}] \; ,
  %   \\
  %   \aidlof*_{Y}
  %   &:=\aidlof{\vphi_F}[\psi_{Y}] \; ,
  %   &\aidlof*^M_{Y}
  %   &:=\aidlof{\vphi_F+\vphi_M}[\psi_{Y}]
  % \end{alignat*}
  \begin{gather*}
    \aidlof*^Y_{D} :=\aidlof{\vphi_F+\phi_Y}[\psi_{D}] \; ,
    \quad
    \aidlof*:=\aidlof{\vphi_F}[\psi_{Y+D}] \; ,
    \quad
    \aidlof*_{Y}:=\aidlof{\vphi_F}[\psi_{Y}] \; ,
    \\
    \newtext*{
      \mathclap{
        \aidlof*^{M,Y}_{D} :=\aidlof{\vphi_F+\vphi_M+\phi_Y}[\psi_{D}] \; ,
        \;\;
        \aidlof*^M :=\aidlof{\vphi_F+\vphi_M}[\psi_{Y+D}] \; ,
        \;\;
        \aidlof*^M_{Y} :=\aidlof{\vphi_F+\vphi_M}[\psi_{Y}]
      }
    }
  \end{gather*}
  and
  % \newtext{define $\aidlof*^{M,Y}_{D}$, $\aidlof*^M$,
  % $\aidlof*^M_{Y}$ similarly with $\vphi_F+\vphi_M$ in place of
  % $\vphi_F$.
  % Also set}
  \begin{equation*}
    \spR{\sheaf F} := \drR q(K_X \otimes \alert{Y \otimes D} \otimes F \otimes
    \sheaf F) \quad\text{ for any sheaf } \sheaf F \text{ on } X 
  \end{equation*}
  for convenience.
  \newtext{Recalling that $\aidlof|n|*^Y_D = \aidlof|0|*_Y$ and
    $\aidlof|0|*^Y_D = \aidlof|0|*$,} the above consideration shows
  that we have the short exact sequence 
  \begin{equation*}
    \xymatrix@R=0.5em{
      {0} \ar[r]
      &{\faidlof|\newtext{n}|/|0|*^Y_{D}} \ar[r]
      &{\faidlof|\newtext{n}|/|0|*} \ar[r]
      &{\faidlof|\newtext{n}|/|0|*_Y} \ar[r]
      &{0}
      \\
      &&&{\frac{\mtidlof{\vphi_F}}{\mtidlof{\vphi_F} \cdot
          \defidlof{Y}}}
      % \ar@{}[u]|(0.45)*[left]+{=}
      \ar@{=}[u]
    } \; ,\footnotemark
  \end{equation*}%
  \footnotetext{\label{fn:general_sigma-prime_explain}\newtext{The
      statement of the current theorem still holds true when 
      $\frac{\mtidlof{\vphi_F}}{\mtidlof{\vphi_F} \cdot \defidlof{Y}}$
      is replaced by $\faidlof|\sigma'|/|0|*_Y$ \newtext*{and the
        index $n$ everywhere in the proof is replaced by $\sigma'$},
      but not for arbitrary $\sigma' \geq 0$ as in Remark
      \ref{remark:general_quotients}.
      The integer $\sigma'$ has to be at least sufficiently large such that
      $\aidlof|\sigma'|*^Y_D = \aidlof|0|*_Y$ for this exact sequence
      to hold.}}%
  which induces the commutative diagram
  \begin{equation*}
    \xymatrix@R=1.2em{
      {\dotsm} \ar[r]
      &{\spR{\faidlof|\newtext{n}|/|0|*}} \ar[r]^-{\iota_{\newtext{n}}}
      \ar[d]_-{\otimes s_{Y+D}}
      \ar[dr]|-*+{\mu_{\newtext{n}}}
      &{\spR{\faidlof|\newtext{n}|/|0|*_Y}} \ar[r]
      \ar[d]^-{\otimes s_Y}
      &{\spR/q+1/{\faidlof|\newtext{n}|/|0|*^Y_{D}}} \ar[r]
      \ar[d]^-{\otimes s_D}
      &{\dotsm}
      \\
      {\dotsm} \ar[r]
      &{\spR M{\faidlof|\newtext{n}|/|0|*^M}} \ar[r]
      &{\spR M{\faidlof|\newtext{n}|/|0|*^M_Y}} \ar[r]
      &{\spR/q+1/ M{\faidlof|\newtext{n}|/|0|*^{M,Y}_{D}}} \ar[r]
      &{\dotsm}
    }
  \end{equation*}
  where the rows are exact, and all columns are induced from $\otimes s$.
  The proof is finished if we show that the map $\otimes s_Y$ is
  injective.
  
  By applying Theorem \ref{thm:main-thm-in-section} (with $(F\otimes
  Y, \vphi_F+\phi_Y)$ in place of $(F, \vphi_F)$ there) to the pair
  $(X,D)$ and using the fact that $\vphi_F+\phi_Y$ is psh, the map
  $\otimes s_D$ is injective.
  By a diagram-chasing argument, to prove that the map $\otimes s_Y$
  is injective, it suffices to show that $\ker\mu_{\newtext{n}} =
  \ker\iota_{\newtext{n}}$.

  % Set $\sigma_{\mlc} := \sigma_{Y+D}$ below for convenience.
  The short exact sequence of adjoint ideal sheaves of the form
  \eqref{eq-ex2} and the inclusion \eqref{eq:Y+D_system-n-Y_system}
  induce the commutative diagram
  %%%%%%%%%%%%%%%%%%%%%%%%%%%
  
\begin{equation*} % \label{eq:commut-diagram_sing-Fujino-lc-conj}
  \begin{gathered}
    \xymatrix@R=0.85cm@C+0.75cm{
      {\vdots} \ar[d]
      & {\vdots} \ar[d]
      & {\vdots} \ar[d]
      \\
      {\spR{\faidlof-1/|0|*}} \ar[d] \ar[r]
      \ar[dr]^-{\iota_{\sigma-1}}
      \ar@/^0.75pc/[drr]|(.6)*+{\mu_{\sigma-1}} 
      & {\spR{\alert{\frac{\aidlof-1*_Y}{\aidlof|0|*_Y}}}}
        \ar[r]^-{\otimes s_Y}
        \ar@{}[d]
        \ar[d]|(.325)*+<3pt>{}
      & {\spR M{\alert{\frac{\aidlof-1*^M_Y}{\aidlof|0|*^M_Y}}}} \ar[d]
      \\
      {\spR{\faidlof/|0|*}} \ar[d] \ar[r]^-{\iota_{\sigma}}
      \ar@/_2.15pc/[rr]|(.65)*+{\mu_{\sigma}}
      & {\spR{\alert{\faidlof|n|/|0|*_Y}}}
      \ar[d]|(.38)*+<3pt>{ }
      \ar[r]^-{\otimes s_Y}
      & {\spR M{\alert{\faidlof|n|/|0|*^M_Y}}}
      \ar[d]
      \\
      {\spR{\residlof*}} \ar[d] \ar[r]^-{\tau_\sigma}
      \ar@/_1.79pc/[rr]+<-39pt,-15pt>|(.67)*+{\nu_\sigma}
      & {\spR{\alert{\frac{\aidlof|n|*_Y}{\aidlof-1*_Y}}}}
      \ar[d]|(.52)*+<3pt>{}
      \ar[r]^-{\otimes s_Y}
      & {\spR M{\alert{\frac{\aidlof|n|*^M_Y}{\aidlof-1*^M_Y}}}} \ar[d] \\
      {\vdots} & {\vdots} & {\vdots} }
  \end{gathered}
\end{equation*}
%
  %%%%%%%%%%%%%%%%%%%%%%%%%%%
  for $\sigma = 1, \dots, \newtext{n}$, analogous to
  \eqref{eq:commut-diagram_sing-Fujino-conj}.
  From the exactness of the columns and the commutativity of the
  diagram, it follows from a diagram-chasing argument that
  $\ker\mu_{\sigma-1} = \ker\iota_{\sigma -1}$ and $\ker\nu_\sigma =
  \ker\tau_\sigma$ implies $\ker\mu_\sigma = \ker\iota_\sigma$ for any
  $\sigma \geq 1$.
  Since $\mu_1 = \nu_1$ and $\iota_1 = \tau_1$, we again have to show
  that
  \begin{equation*} % \label{eq:ker-nu=ker-tau}
    \paren{\ker\nu_\sigma}_t =\paren{\ker\tau_\sigma}_t
    \quad\text{ for all } \sigma = 1, \dots, \newtext{n}
    \:\text{ and for all } t \in \Delta \; .
  \end{equation*}
  \newtext{Notice that the result of the analogous equalities in
    the proof of Theorem \ref{thm:main-thm-in-section} cannot be
    applied directly (with $Y+D$ in place of $D$ there), as the maps
    $\nu_\sigma$ and $\tau_\sigma$ here involve the restriction from
    $Y+D$ to $Y$ while those maps in the proof of Theorem
    \ref{thm:main-thm-in-section} remain on $Y+D$.
    However,}
  the proof of Theorem
  \ref{thm:main-thm-in-section}
  % mutatis mutandis (with, at worst, a
  % replacement of ``$\aidlof|n|*^M$'' by
  % ``$\aidlof|n|*^M_Y$'' in Step
  % \ref{item:express-su-in-residue-norm} of the proof).
  \begin{newpara}\label{page:explain_modified-proof}
    can still be applied to obtain the above equalities, with two modifications:
    \begin{itemize}
    \item the commutative diagram in
      \eqref{eq:commut-diagram_delta-tau} (in Step
      \ref{item:pf:step-2} of the proof) is replaced by
      \begin{equation*}
        \begin{aligned}
          \xymatrix@R-0.3cm{
            {\dotsm} \ar[r]
            & {\spH/q-1/{\residlof+1*}_c} \ar[r]^-{\delta}
            \ar[d]^-{\tau_{\sigma+1}}
            & {\spH{\residlof*}_c} \ar[r]
            \ar[d]
            \ar[dr]^-{\tau_\sigma}
            & {\spH{\faidlof+1/-1*}_c} \ar[r]
            \ar[d]
            & {\dotsm}
            \\
            {\dotsm} \ar[r]
            & {\spH/q-1/{\faidlof|\newtext{n}|*_Y}_c} \ar[r]
            & {\spH{\residlof*_Y}_c} \ar[r]
            & {\spH{\faidlof|\newtext{n}|/-1*_Y}_c} \ar[r]
            & {\dotsm}
          }
        \end{aligned}
      \end{equation*}
      with exact rows, where $\residlof*_Y \isom \faidlof/-1*_Y$ and
      the vertical maps are induced by the inclusions $\aidlof*
      \subset \aidlof*_Y \subset \aidlof|\newtext{n}|*_Y$ for any integers $\sigma \in
      [0,\newtext{n}]$;

    \item ``$\aidlof*^M$'' in Step
      \ref{item:express-su-in-residue-norm} of the proof is replaced by
      ``$\aidlof*^M_Y$'' for any integers $\sigma \geq 0$.
    \end{itemize}
  \end{newpara}
  This completes the proof.
\end{proof}

% \subsection{Extension theorem on $X$ (?)}
% \label{sec:proof-extension_thm}

%%%%% Reference list %%%%%

\end{document}